\documentclass[10pt]{article}

\RequirePackage{amsmath}
\RequirePackage{amsfonts}
\RequirePackage{amssymb}
\RequirePackage{hyperref}
\RequirePackage{multirow}
\RequirePackage{algorithm}
\RequirePackage{algpseudocode}
\RequirePackage{graphicx}
\RequirePackage{xcolor}
\RequirePackage{amsthm}
\RequirePackage{amsmath}
\RequirePackage{amsfonts}
\RequirePackage{amssymb}
\RequirePackage{hyperref}
\RequirePackage{multirow}
\RequirePackage{algorithm}
\RequirePackage{algpseudocode}
\RequirePackage{graphicx}
\RequirePackage{xcolor}
\usepackage{booktabs} 
\usepackage{array}
\usepackage{arydshln}
\usepackage{dpr_fec,ifthen,dsfont}

\newtheorem{definition}{Definition}
\newtheorem{lemma}{Lemma}
\newtheorem{proposition}{Proposition}

\newcommand{\Halmos}{\qed}

\usepackage{isetechreport}
\def\reportyear{21T}
\def\reportno{013}
\def\revisionno{1}
\def\originaldate{June 24, 2021}
\def\revisiondate{\today}
\coraltrue
\cvcrfalse

\begin{document}

\title {A Stochastic Sequential Quadratic Optimization Algorithm for Nonlinear Equality Constrained Optimization with Rank-Deficient Jacobians}

\author{Albert S.~Berahas\thanks{E-mail: albertberahas@gmail.com}}
\affil{Department of Industrial and Operations Engineering, University of Michigan}
\author{Frank E.~Curtis\thanks{E-mail: frank.e.curtis@lehigh.edu}}
\author{Michael J.~O'Neill\thanks{E-mail: moneill@lehigh.edu}}
\author{Daniel P.~Robinson\thanks{E-mail: daniel.p.robinson@lehigh.edu}}
\affil{Department of Industrial and Systems Engineering, Lehigh University}

\titlepage

\maketitle

%**********
% Abstract
%**********
\begin{abstract}
  A sequential quadratic optimization algorithm is proposed for solving smooth nonlinear equality constrained optimization problems in which the objective function is defined by an expectation of a stochastic function.  The algorithmic structure of the proposed method is based on a step decomposition strategy that is known in the literature to be widely effective in practice, wherein each search direction is computed as the sum of a normal step (toward linearized feasibility) and a tangential step (toward objective decrease in the null space of the constraint Jacobian).  However, the proposed method is unique from others in the literature in that it both allows the use of stochastic objective gradient estimates and possesses convergence guarantees even in the setting in which the constraint Jacobians may be rank deficient.  The results of numerical experiments demonstrate that the algorithm offers superior performance when compared to popular alternatives.
\end{abstract}

%**************
% Body of paper
%**************
%% Macros
\newcommand{\dtrue}{d^{\rm true}}
\newcommand{\utrue}{u^{\rm true}}
\newcommand{\wtrue}{w^{\rm true}}
\newcommand{\ytrue}{y^{\rm true}}
\newcommand{\tautrial}{\tau^{\rm trial}}
\newcommand{\tautruetrial}{\tau^{\rm trial,true}}
\newcommand{\xitrial}{\xi^{\rm trial}}
\newcommand{\chitrial}{\chi^{\rm trial}}
\newcommand{\alphasuff}{\alpha^{\rm suff}}
\newcommand{\alphamin}{\alpha^{\rm min}}
\newcommand{\alphatrial}{\alpha^{\rm trial}}
\newcommand{\toedit}[1]{\textcolor{blue}{#1}}
\renewcommand{\big}{{\rm big}}
\newcommand{\low}{{\rm low}}
\newcommand{\zero}{{\rm zero}}

%*********
% Section
%*********
\section{Introduction}\label{sec.introduction}

We propose an algorithm for solving equality constrained optimization problems in which the objective function is defined by an expectation of a stochastic function.  Formulations of this type arise throughout science and engineering in important applications such as data-fitting problems, where one aims to determine a model that minimizes the discrepancy between values yielded by the model and corresponding known outputs.

Our algorithm is designed for solving such problems when the decision variables are restricted to the solution set of a (potentially nonlinear) set of equations.  We are particularly interested in such problems when the constraint Jacobian---i.e., the matrix of first-order derivatives of the constraint function---may be rank deficient in some or even all iterations during the run of an algorithm, since this can be an unavoidable occurrence in practice that would ruin the convergence properties of any algorithm that is not specifically designed for this setting.  The structure of our algorithm follows a step decomposition strategy that is common in the constrained optimization literature; in particular, our algorithm has roots in the Byrd-Omojokun approach~\cite{Omoj89}.  However, our algorithm is unique from previously proposed algorithms in that it offers convergence guarantees while allowing for the use of stochastic objective gradient information in each iteration.  We prove that our algorithm converges to stationarity (in expectation), both in nice cases when the constraints are feasible and convergence to the feasible region can be guaranteed (in expectation), and in more challenging cases, such as when the constraints are infeasible and one can only guarantee convergence to an infeasible stationary point.  To the best of our knowledge, there exist no other algorithms in the literature that have been designed specifically for this setting, namely, stochastic optimization with equality constraints that may exhibit rank deficiency.

The step decomposition strategy employed by our algorithm makes it similar to the method proposed in \cite{CurtNoceWaec09}, although that method is designed for deterministic optimization only and employs a line search, whereas our approach is designed for stochastic optimization and requires no line searches.  Our algorithm builds upon the method for solving equality constrained stochastic optimization problems proposed in~\cite{BeraCurtRobiZhou21}.  The method proposed in that article assumes that the singular values of the constraint Jacobians are bounded below by a positive constant throughout the optimization process, which implies that the linear independence constraint qualification (LICQ) holds at all iterates.  By contrast, the algorithm proposed in this paper makes no such assumption.  Handling the potential lack of full-rank Jacobians necessitates a different algorithmic structure and a distinct approach to proving convergence guarantees; e.g., one needs to account for the fact that primal-dual stationarity conditions may not be necessary and/or the constraints may be infeasible.

Similar to the context in \cite{BeraCurtRobiZhou21}, our algorithm is intended for the highly stochastic regime in which the stochastic gradient estimates might only be unbiased estimators of the gradients of the objective at the algorithm iterates that satisfy a loose variance condition.  Indeed, we show that in nice cases---in particular, when the adaptive merit parameter employed in our algorithm eventually settles at a value that is sufficiently small---our algorithm has convergence properties in expectation that match those of the algorithm in \cite{BeraCurtRobiZhou21}.  These results parallel those for the stochastic gradient method in the context of unconstrained optimization \cite{BottCurtNoce18,RobbMonr51,RobbSieg71}.  However, for cases not considered in \cite{BeraCurtRobiZhou21} when the merit parameter sequence may vanish, we require the stronger assumption that the difference between each stochastic gradient estimate and the corresponding true gradient of the objective eventually is bounded \emph{deterministically} in each iteration.  This is appropriate in many ways since in such a scenario the algorithm aims to transition from solving a \emph{stochastic} optimization problem to the \emph{deterministic} one of minimizing constraint violation.  Finally, we show under reasonable assumptions the total probability is zero that the merit parameter settles at too large of a value.

Our algorithm has some similarities, but many differences with another recently proposed algorithm, namely, that in \cite{NaAnitKola21}.  That algorithm is also designed for equality constrained stochastic optimization, but: $(i)$ like for the algorithm in \cite{BeraCurtRobiZhou21}, for the algorithm in \cite{NaAnitKola21} the LICQ is assumed to hold at all algorithm iterates, and $(ii)$ the algorithm in \cite{NaAnitKola21} employs an adaptive line search that may require the algorithm to compute relatively accurate stochastic gradient estimates throughout the optimization process.  Our algorithm, on the other hand, does not require the LICQ to hold and is meant for a more stochastic regime, meaning that it does not require a procedure for refining the stochastic gradient estimate within an iteration.  Consequently, the convergence guarantees that can be proved for our method, and the expectations that one should have about the practical performance of our method, are quite distinct from those for the algorithm in \cite{NaAnitKola21}.

Besides the methods in \cite{BeraCurtRobiZhou21,NaAnitKola21}, there have been few proposed algorithms that might be used to solve problem of the form~\eqref{prob.opt}.  Some methods have been proposed that employ stochastic (proximal) gradient strategies applied to minimizing penalty functions derived from constrained problems \cite{ChenTungVeduMori18,KumaSoumMhamHara18,NandPathAbhiSing19}, but these do not offer convergence guarantees to stationarity with respect to the original constrained problem.  On the other hand, stochastic Frank-Wolfe methods have been proposed \cite{hazan2016variance,locatello2019stochastic,lu2020generalized,RaviDinhLokhSing19,reddi2016stochastic,zhang2020one}, but these can only be applied in the context of convex feasible regions.  Our algorithm, by contrast, is designed for nonlinear equality constrained stochastic optimization.

%************
% Subsection
%************
\subsection{Notation}

The set of real numbers is denoted as $\R{}$, the set of real numbers greater than (respectively,~greater than or equal to) $r \in \R{}$ is denoted as $\R{}_{> r}$ (respectively,~$\R{}_{\geq r}$), the set of $n$-dimensional real vectors is denoted as $\R{n}$, the set of $m$-by-$n$-dimensional real matrices is denoted as $\R{m \times n}$, and the set of $n$-by-$n$-dimensional real symmetric matrices is denoted as $\mathbb{S}^n$.  Given $J \in \R{m \times n}$, the range space of $J^T$ is denoted as $\Range(J^T)$ and the null space of $J$ is denoted as $\Null(J)$.  (By the Fundamental Theorem of Linear Algebra, for any $J \in \R{m \times n}$, the spaces $\Range(J^T)$ and $\Null(J)$ are orthogonal and $\Range(J^T) + \Null(J) = \R{n}$, where in this instance `$+$' denotes the Minkowski sum operator.)  The set of nonnegative integers is denoted as $\N{} := \{0,1,2,\dots\}$.  For any $m \in \N{}$, let~$[m]$ denote the set of integers $\{0,1,\dots,m\}$. Correspondingly, to represent a set of vectors $\{v_0,\dots,v_k\}$, we define $v_{[k]} := \{v_0,\dots,v_k\}$.

The algorithm that we propose is iterative in the sense that, given a starting point $x_0 \in \R{n}$, it generates a sequence of iterates $\{x_k\}$ with $x_k \in \R{n}$ for all $k \in \N{}$.  For simplicity of notation, the iteration number is appended as a subscript to other quantities corresponding to each iteration; e.g., with a function $c : \R{n} \to \R{}$, its value at $x_k$ is denoted as $c_k := c(x_k)$ for all $k \in \N{}$.  Given $J_k \in \R{m \times n}$, we use $Z_k$ to denote a matrix whose columns form an orthonormal basis for $\Null(J_k)$.

%************
% Subsection
%************
\subsection{Organization}

Our problem of interest and basic assumptions about the problem and the behavior of our algorithm are presented in Section~\ref{sec.problem}.  Our algorithm is motivated and presented in Section~\ref{sec.algorithm}.  Convergence guarantees for our algorithm are presented in Section~\ref{sec.convergence}.  The results of numerical experiments are provided in Section~\ref{sec.numerical} and concluding remarks are provided in Section~\ref{sec.conclusion}.

%*********
% Section
%*********
\section{Problem Statement}\label{sec.problem}

Our algorithm is designed for solving (potentially nonlinear and/or nonconvex) equality constrained optimization problems of the form
\bequation\label{prob.opt}
  \min_{x\in\R{n}}\ f(x)\ \ \st\ \ c(x) = 0,\ \ \text{with}\ \ f(x) = \E[F(x,\iota)],
\eequation
where the functions $f : \R{n} \to \R{}$ and $c : \R{n} \to \R{m}$ are smooth, $\iota$ is a random variable with associated probability space $(\Omega,\Fcal,P)$, $F : \R{n} \times \Omega \to \R{}$, and $\E[\cdot]$ denotes expectation taken with respect to~$P$.  We assume that values and first-order derivatives of the constraint functions can be computed, but that the objective and its associated first-order derivatives are intractable to compute, and one must instead employ stochastic estimates.  (We formalize our assumptions about such stochastic estimates starting with Assumption~\ref{ass.g} on page~\pageref{ass.g}.)  Formally, we make the following assumption with respect to \eqref{prob.opt} and our proposed algorithm, which generates a sequence of iterates $\{x_k\}$.

\bassumption\label{ass.main}
  Let $\Xcal \subseteq \R{n}$ be an open convex set containing the sequence~$\{x_k\}$ generated by any run of the algorithm.  The objective function $f : \R{n} \to \R{}$ is continuously differentiable and bounded over~$\Xcal$ and its gradient function $\nabla f : \R{n} \to \R{n}$ is Lipschitz continuous with constant $L \in \R{}_{>0}$ $($with respect to $\|\cdot\|_2$$)$ and bounded over $\Xcal$.  The constraint function $c : \R{n} \to \R{m}$ $($with $m \leq n$$)$ is continuously differentiable and bounded over $\Xcal$ and its Jacobian function $J := \nabla c^T : \R{n} \to \R{m \times n}$ is Lipschitz continuous with constant $\Gamma \in \R{}_{>0}$ $($with respect to $\|\cdot\|_2$$)$ and bounded over~$\Xcal$.
\eassumption

The aspects of Assumption~\ref{ass.main} that pertain to the objective function~$f$ and constraint function~$c$ are typical for the equality constrained optimization literature.  Notice that we do not assume that the iterate sequence itself is bounded.  Under Assumption~\ref{ass.main}, it follows that there exist positive real numbers $(f_{\inf},f_{\sup},\kappa_{\nabla f},\kappa_c,\kappa_J) \in \R{}_{>0} \times \R{}_{>0} \times \R{}_{>0} \times \R{}_{>0} \times \R{}_{>0}$ such that
\bequation\label{eq.bounds}
  f_{\inf} \leq f_k \leq f_{\sup},\ \ \|\nabla f(x_k)\|_2 \leq \kappa_{\nabla f},\ \ \|c_k\|_2 \leq \kappa_c,\ \ \text{and}\ \ \|J_k\|_2 \leq \kappa_J\ \ \text{for all}\ \ k \in \N{}.
\eequation

Given that our proposed algorithm is stochastic, it is admittedly not ideal to have to assume that the objective value, objective gradient, constraint value, and constraint Jacobian are bounded over the set $\Xcal$ containing the iterates.  This is a common assumption in the deterministic optimization literature, where it may be justified in the context of an algorithm that is guaranteed to make progress in each iteration, say with respect to a merit function.  However, for a stochastic algorithm such as ours, such a claim may be seen as less than ideal since a stochastic algorithm may only be guaranteed to make progress \emph{in expectation} in each iteration, meaning that it is possible for the iterates to drift far from desirable regions of the search space during the optimization process.

Our justification for Assumption~\ref{ass.main} is two-fold.  First, any reader who is familiar with analyses of stochastic algorithms for unconstrained optimization---in particular, those analyses that do not require that the objective gradient is bounded over a set containing the iterates---should appreciate that additional challenges present themselves in the context of constrained optimization.  For example, whereas in unconstrained optimization one naturally considers the objective $f$ as a measure of progress, in (nonconvex) constrained optimization one needs to employ a merit function for measuring progress, and for practical purposes such a function typically needs to involve a parameter (or parameters) that must be adjusted dynamically by the algorithm.  One finds that it is the adaptivity of our merit parameter (see \eqref{eq.merit} later on) that necessitates the aforementioned boundedness assumptions that we use in our analysis.  (Certain exact merit functions, such as that employed in \cite{NaAnitKola21}, might not lead to the same issues as the merit function that we employ.  However, we remark that the merit function employed in \cite{NaAnitKola21} is not a viable option unless the LICQ holds at all algorithm iterates.)  Our second justification is that we know of no other algorithm that offers convergence guarantees that are as comprehensive as ours (in terms of handling feasible, degenerate, and infeasible settings) under an assumption that is at least as loose as Assumption~\ref{ass.main}.

Let the Lagrangian $\ell : \R{n} \times \R{m} \to \R{}$ corresponding to \eqref{prob.opt} be given by $\ell(x,y) = f(x) + c(x)^Ty$, where $y \in \R{m}$ represents a vector of Lagrange multipliers.  Under a constraint qualification (such as the LICQ), necessary conditions for first-order stationarity with respect to \eqref{prob.opt} are given by
\bequation\label{eq.KKT}
  0 = \bbmatrix \nabla_x \ell(x,y) \\ \nabla_y \ell(x,y) \ebmatrix = \bbmatrix \nabla f(x) + J(x)^Ty \\ c(x) \ebmatrix;
\eequation
see, e.g., \cite{NoceWrig06}.  However, under only Assumption~\ref{ass.main}, it is possible for \eqref{prob.opt} to be degenerate---in which case \eqref{eq.KKT} might not be necessary at a solution of \eqref{prob.opt}---or \eqref{prob.opt} may be infeasible.  In the latter case, one aims to design an algorithm that transitions automatically from seeking stationarity with respect to \eqref{prob.opt} to seeking stationarity with respect to a measure of infeasibility of the constraints.  For our purposes, we employ the infeasibility measure $\varphi : \R{n} \to \R{}$ defined by $\varphi(x) = \|c(x)\|_2$.  A point $x \in \R{n}$ is stationary with respect to $\varphi$ if and only if either $c(x) = 0$ or both $c(x) \neq 0$ and
\bequation\label{eq.infeasible_stationary}
  0 = \nabla \varphi(x) = \frac{J(x)^T c(x)}{\|c(x)\|_2}.
\eequation

%*********
% Section
%*********
\section{Algorithm Description}\label{sec.algorithm}

Our algorithm can be characterized as a sequential quadratic optimization (commonly known as SQP) method that employs a step decomposition strategy and chooses step sizes that attempt to ensure sufficient decrease in a merit function in each iteration.  We present our complete algorithm in this section, which builds upon this basic characterization to involve various unique aspects that are designed for handling the combination of $(i)$ stochastic gradient estimates and $(ii)$ potential rank deficiency of the constraint Jacobians.

In each iteration $k \in \N{}$, the algorithm first computes the \emph{normal component} of the search direction toward reducing linearized constraint violation.  Conditioned on the event that $x_k$ is reached as the $k$th iterate, the problem defining this computation, namely,
\bequation\label{prob.normal}
  \min_{v \in \R{n}}\ \thalf \|c_k + J_kv\|_2^2\ \ \st\ \ \|v\|_2 \leq \omega \|J_k^Tc_k\|_2
\eequation
where $\omega \in \R{}_{>0}$ is a user-defined parameter, is deterministic since the constraint function value~$c_k$ and constraint Jacobian $J_k$ are available.  If $J_k$ has full row rank, $\omega$ is sufficiently large, and \eqref{prob.normal} is solved to optimality, then one obtains $v_k$ such that $c_k + J_kv_k = 0$.  However, an exact solution of~\eqref{prob.normal} may be expensive to obtain, and---as has been shown for various step decomposition strategies, such as the Byrd-Omojokun approach~\cite{Omoj89}---the consideration of \eqref{prob.normal} is viable when $J_k$ might not have full row rank.  Fortunately, our algorithm merely requires that the normal component $v_k \in \R{n}$ is feasible for problem~\eqref{prob.normal}, lies in $\Range(J_k^T)$, and satisfies the Cauchy decrease condition
\bequation\label{eq.Cauchy}
  \|c_k\|_2 - \|c_k + J_kv_k\|_2 \geq \epsilon_v (\|c_k\|_2 - \|c_k + \alpha_k^C J_k v_k^C\|_2)
\eequation
for some user-defined parameter $\epsilon_v \in (0,1]$.  Here, $v_k^C := -J_k^Tc_k$ is the steepest descent direction for the objective of problem~\eqref{prob.normal} at $v=0$ and the step size $\alpha_k^C \in \R{}$ is the unique solution to the problem to minimize $\thalf \|c_k + \alpha^C J_k v_k^C\|_2^2$ over $\alpha^C \in \R{}_{\geq0}$ subject to $\alpha^C \leq \omega$ (see, e.g., \cite[Equations~(4.11)--(4.12)]{NoceWrig06}).  Since this allows one to choose $v_k \gets v_k^C$, the normal component can be computed at low computational cost.  For a more accurate solution to \eqref{prob.normal}, one can employ a so-called matrix-free iterative algorithm such as the linear conjugate gradient (CG) method with Steihaug stopping conditions~\cite{Stei83} or GLTR \cite{GoulLuciRomaToin99}, each of which is guaranteed to yield a solution satisfying the aforementioned conditions no matter how many iterations (greater than or equal to one) are performed.

After the computation of the normal component, our algorithm computes the \emph{tangential component} of the search direction by minimizing a model of the objective function subject to remaining in the null space of the constraint Jacobian.  This ensures that the progress toward linearized feasibility offered by the normal component is not undone by the tangential component when the components are added together.  The problem defining the computation of the tangential component is
\bequation\label{prob.tangential}
  \min_{u\in\R{n}}\ (g_k + H_kv_k)^Tu + \thalf u^TH_ku\ \ \st\ \ J_ku = 0,
\eequation
where $g_k \in \R{n}$ is a stochastic gradient estimate \emph{at least} satisfying Assumption~\ref{ass.g} below and the real symmetric matrix $H_k \in \mathbb{S}^n$ satisfies Assumption~\ref{ass.H} below.  (Specific additional requirements for $\{g_k\}$ are stated separately for each case in our convergence analysis.)

\bassumption\label{ass.g}
  For all $k \in \N{}$, the stochastic gradient estimate $g_k \in \R{n}$ is an unbiased estimator of $\nabla f(x_k)$, i.e., $\E_k[g_k] = \nabla f(x_k)$, where $\E_k[\cdot]$ denotes expectation conditioned on the event that the algorithm has reached $x_k$ as the $k$th iterate.  In addition, there exists a positive real number $M \in \R{}_{>0}$ such that, for all $k \in \N{}$, one has $\E_k[\|g_k - \nabla f(x_k)\|_2^2] \leq M$.
\eassumption

\bassumption\label{ass.H}
  The matrix $H_k \in \mathbb{S}^n$ is chosen independently from $g_k$ for all $k \in \N{}$, the sequence~$\{H_k\}$ is bounded in norm by $\kappa_H \in \R{}_{>0}$, and there exists $\zeta \in \R{}_{>0}$ such that, for all $k \in \N{}$, one has $u^TH_ku \geq \zeta \|u\|_2^2$ for all $u \in \Null(J_k)$.
\eassumption

In our context, one can generate $g_k$ in iteration $k \in \N{}$ by independently drawing $b_k$ realizations of the random variable $\iota$, denoting the \emph{mini-batch} as $\Bcal_k := \{\iota_{k,1},\dots,\iota_{k,b_k}\}$, and setting
\bequation\label{eq.mini-batch}
  g_k \gets \frac{1}{b_k} \sum_{\iota \in \Bcal_k} \nabla f(x_k,\iota).
\eequation
It is a modest assumption about the function $f$ and the sample sizes $\{b_k\}$ to say that $\{g_k\}$ generated in this manner satisfies Assumption~\ref{ass.g}.  As for Assumption~\ref{ass.H}, the assumptions that the elements of $\{H_k\}$ are bounded in norm and that $H_k$ is sufficiently positive definite in $\Null(J_k)$ for all $k \in \N{}$ are typical for the constrained optimization literature.  In practice, one may choose $H_k$ to be (an approximation of) the Hessian of the Lagrangian at $(x_k,y_k)$ for some $y_k$, if such a matrix can be computed with reasonable effort in a manner that guarantees that Assumption~\ref{ass.H} holds.  A simpler alternative is that $H_k$ can be set to some positive definite diagonal matrix (independent of~$g_k)$.

Under Assumption~\ref{ass.H}, the tangential component $u_k$ solving \eqref{prob.tangential} can be obtained by solving
\bequation\label{eq.linsys}
  \bbmatrix H_k & J_k^T \\ J_k & 0 \ebmatrix \bbmatrix u_k \\ y_k \ebmatrix = -\bbmatrix g_k + H_kv_k \\ 0 \ebmatrix.
\eequation
Even if the constraint Jacobian $J_k$ does not have full row rank, the linear system \eqref{eq.linsys} is consistent since it represents sufficient optimality conditions (under Assumption~\ref{ass.H}) of the linearly constrained quadratic optimization problem in \eqref{prob.tangential}.  (Factorization methods that are popular in the context of solving symmetric indefinite linear systems of equations, such as the Bunch-Kaufman factorization, can fail when the matrix in \eqref{eq.linsys} is singular.  However, Krylov subspace methods provide a viable alternative, since for such methods singularity is benign as long as the system is known to be consistent, as is the case for \eqref{eq.linsys}.)  Under Assumption~\ref{ass.H}, the solution component~$u_k$ is unique, although the component~$y_k$ might not be unique (if $J_k$ does not have full row rank).

Upon computation of the search direction, our algorithm proceeds to determining a positive step size.  For this purpose, we employ the merit function $\phi : \R{n} \times \R{}_{\geq0} \to \R{}$ defined by
\bequation\label{eq.merit}
  \phi(x,\tau) = \tau f(x) + \|c(x)\|_2
\eequation
where $\tau$ is a merit parameter whose value is set dynamically.  The function~$\phi$ is a type of exact penalty function that is common in the literature \cite{Han77,HanMang79,Powe78}.  For setting the merit parameter value in each iteration, we employ a local model of $\phi$ denoted as $l : \R{n} \times \R{}_{\geq0} \times \R{n} \times \R{n} \to \R{}$ and defined~by
\bequationNN
  l(x,\tau,g,d) = \tau (f(x) + g^Td) + \|c(x) + J(x)d\|_2.
\eequationNN
Given the search direction vectors $v_k$, $u_k$, and $d_k \gets v_k + u_k$, the algorithm sets
\bequation\label{eq.merit_parameter_trial}
  \tautrial_k \gets \bcases \infty & \text{if $g_k^Td_k + u_k^TH_ku_k \leq 0$} \\[5pt] \displaystyle \frac{(1 - \sigma)(\|c_k\|_2 - \|c_k + J_kd_k\|_2)}{g_k^Td_k + u_k^TH_ku_k} & \text{otherwise,} \ecases
\eequation
where $\sigma \in (0,1)$ is user-defined.  The merit parameter value is then set as
\bequation\label{eq.tau_update}
  \tau_k \gets \bcases \tau_{k-1} & \text{if $\tau_{k-1} \leq \tautrial_k$} \\ \min\{(1 - \epsilon_\tau) \tau_{k-1}, \tautrial_k\} & \text{otherwise,} \ecases
\eequation
where $\epsilon_\tau \in (0,1)$ is user-defined.  This rule ensures that $\{\tau_k\}$ is monotonically nonincreasing, $\tau_k \leq \tautrial_k$ for all $k \in \N{}$, and, with the reduction function $\Delta l : \R{n} \times \R{}_{\geq0} \times \R{n} \times \R{n} \to \R{}$ defined by
\bequation\label{eq.model_reduction}
  \Delta l(x,\tau,g,d) = l(x,\tau,g,0) - l(x,\tau,g,d) = -\tau g^Td + \|c(x)\|_2 - \|c(x) + J(x)d\|_2
\eequation
and Assumption~\ref{ass.H}, it ensures the following fact that is critical for our analysis:
\bequation\label{eq.model_reduction_condition}
  \Delta l(x_k,\tau_k,g_k,d_k) \geq \tau_k u_k^TH_ku_k + \sigma(\|c_k\|_2 - \|c_k + J_kv_k\|_2).
\eequation

Similar to the algorithm in \cite{BeraCurtRobiZhou21}, our algorithm also adaptively sets other parameters that are used for determining an allowable range for the step size in each iteration.  (There exist constants that, if known in advance, could be used by the algorithm for determining the allowable range for each step size; see Lemma~\ref{lem.chi_upper} in our analysis later on.  However, to avoid the need to know these problem-dependent constants in advance, our algorithm generates these parameter sequences adaptively, which our analysis shows is sufficient to ensure convergence guarantees.)  For distinguishing between search directions that are dominated by the tangential component and others that are dominated by the normal component, the algorithm adaptively defines sequences $\{\chi_k\}$ and $\{\zeta_k\}$.  (These sequences were not present in the algorithm in~\cite{BeraCurtRobiZhou21}; they are newly introduced for the needs of our proposed algorithm.)  In particular, in iteration $k \in \N{}$, the algorithm employs the conditions
\bequation\label{eq.switch_condition}
  \|u_k\|_2^2 \geq \chi_{k-1} \|v_k\|_2^2\ \ \text{and}\ \ \thalf d_k^TH_kd_k < \tfrac14 \zeta_{k-1} \|u_k\|_2^2
\eequation
in order to set
\bequation\label{eq.switch}
  (\chi_k,\zeta_k) \gets \bcases ((1 + \epsilon_\chi)\chi_{k-1},(1 - \epsilon_\zeta)\zeta_{k-1}) & \text{if \eqref{eq.switch_condition} holds} \\ (\chi_{k-1},\zeta_{k-1}) & \text{otherwise,} \ecases
\eequation
where $\epsilon_\chi \in \R{}_{>0}$ and $\epsilon_\zeta \in (0,1)$ are user-defined.  It follows from \eqref{eq.switch} that $\{\chi_k\}$ is monotonically nondecreasing and $\{\zeta_k\}$ is monotonically nonincreasing.  It will be shown in our analysis that $\{\chi_k\}$ is bounded above by a positive real number and $\{\zeta_k\}$ is bounded below by a positive real number, where these bounds are uniform over all runs of the algorithm; i.e., these sequences are bounded \emph{deterministically}.  This means that despite the stochasticity of the algorithm iterates, these sequences have $(\chi_k,\zeta_k) = (\chi_{k-1},\zeta_{k-1})$ for all sufficiently large $k \in \N{}$ in any run of the algorithm.

Whether $\|u_k\|_2^2 \geq \chi_k \|v_k\|_2^2$ (i.e., the search direction is \emph{tangentially dominated}) or $\|u_k\|_2^2 < \chi_k \|v_k\|_2^2$ (i.e., the search direction is \emph{normally dominated}) influences two aspects of iteration $k \in \N{}$.  First, it influences a value that the algorithm employs to determine the range of allowable step sizes that represents a lower bound for the ratio between the reduction in the model $l$ of the merit function and a quantity involving the squared norm of the search direction.  (A similar, but slightly different sequence was employed for the algorithm in \cite{BeraCurtRobiZhou21}.)  In iteration $k \in \N{}$ of our algorithm, the estimated lower bound is set adaptively by first setting
\bequation\label{eq.ratio_trial}
  \xitrial_k \gets \bcases \displaystyle \frac{\Delta l(x_k,\tau_k,g_k,d_k)}{\tau_k \|d_k\|_2^2} & \text{if $\|u_k\|_2^2 \geq \chi_k \|v_k\|_2^2$} \\ \displaystyle \frac{\Delta l(x_k,\tau_k,g_k,d_k)}{\|d_k\|_2^2} & \text{otherwise,} \ecases
\eequation
then setting
\bequation\label{eq.ratio_update}
  \xi_k \gets \bcases \xi_{k-1} & \text{if $\xi_{k-1} \leq \xitrial_k$} \\ \min\{(1-\epsilon_{\xi}) \xi_{k-1}, \xitrial_k\} & \text{otherwise,} \ecases
\eequation
for some user-defined $\epsilon_\xi \in (0,1)$.  The procedure in \eqref{eq.ratio_update} ensures that $\{\xi_k\}$ is monotonically nonincreasing and $\xi_k \leq \xitrial_k$ for all $k \in \N{}$.  It will be shown in our analysis that $\{\xi_k\}$ is bounded away from zero \emph{deterministically}, even though in each iteration it depends on stochastic quantities.  (Like for $\{\chi_k\}$ and $\{\zeta_k\}$, there exists a constant that, if known in advance, could be used in place of $\xi_k$ for all $k \in \N{}$---see Lemma~\ref{lem.xi_lower}---but for ease of employment our algorithm generates $\{\xi_k\}$ instead.)  To achieve this property, it is critical that the denominator in \eqref{eq.ratio_trial} is different depending on whether the search direction is tangentially or normally dominated; see Lemma~\ref{lem.xi_lower} later on for details.  The second aspect of the algorithm that is affected by whether a search direction is tangentially or normally dominated is a rule for setting the step size; this will be seen in \eqref{eq.alpha} later on.

We are now prepared to present the mechanism by which a positive step size is selected in each iteration $k \in \N{}$ of our algorithm.  We present a strategy that allows for our convergence analysis in Section~\ref{sec.convergence} to be as straightforward as possible.  In Section~\ref{sec.numerical}, we remark on extensions of this strategy that are included in our software implementation for which our convergence guarantees also hold (as long as some additional cases are considered in one key lemma).

We motivate our strategy by considering an upper bound for the change in the merit function corresponding to the computed search direction, namely, $d_k \gets v_k + u_k$.  In particular, under Assumption~\ref{ass.main}, in iteration $k \in \N{}$, one has for any nonnegative step size $\alpha \in \R{}_{\geq0}$ that
\bequation\label{eq.merit_decrease}
    \baligned
      &\ \phi(x_k + \alpha d_k,\tau_k) - \phi(x_k,\tau_k) \\
      =&\ \tau_k f(x_k + \alpha d_k) - \tau_k f(x_k) + \|c(x_k + \alpha d_k)\|_2 - \|c_k\|_2 \\
      \leq&\ \alpha \tau_k \nabla f(x_k)^Td_k + \|c_k + \alpha J_k d_k\|_2 - \|c_k\|_2 + \thalf (\tau_k L + \Gamma) \alpha^2 \|d_k\|_2^2 \\
      \leq&\ \alpha \tau_k \nabla f(x_k)^Td_k + |1-\alpha|\|c_k\|_2 - \|c_k\|_2 + \alpha\| c_k + J_k d_k\|_2 + \thalf (\tau_k L + \Gamma) \alpha^2 \|d_k\|_2^2.
  \ealigned
\eequation
This upper bound is a convex, piecewise quadratic function in $\alpha$.  In a deterministic algorithm in which the gradient $\nabla f(x_k)$ is available, it is common to require that the step size $\alpha$ yields
\bequation\label{eq.suff_dec}
  \phi(x_k + \alpha d_k,\tau_k) - \phi(x_k,\tau_k) \leq -\eta \alpha \Delta l(x_k,\tau_k,\nabla f(x_k),d_k),
\eequation
where $\eta \in (0,1)$ is user-defined.  However, in our setting, \eqref{eq.suff_dec} cannot be enforced since our algorithm avoids the evaluation of $\nabla f(x_k)$ and in lieu of it only computes a stochastic gradient $g_k$.  The first main idea of our step size strategy is to determine a step size such that the upper bound in \eqref{eq.merit_decrease} is less than or equal to the right-hand side of \eqref{eq.suff_dec} when the true gradient $\nabla f(x_k)$ is replaced by its estimate~$g_k$.  Since~\eqref{eq.model_reduction_condition}, the orthogonality of $v_k \in \Range(J_k^T)$ and $u_k \in \Null(J_k)$, and the properties of the normal step (which, as shown in Lemma~\ref{lem.normal_component} later on, include that the left-hand side of \eqref{eq.Cauchy} is positive whenever $v_k \neq 0$) ensure that $\Delta l(x_k,\tau_k,g_k,d_k) > 0$ whenever $d_k \neq 0$, it follows that a step size satisfying this aforementioned property is given, for any $\beta_k \in (0,1]$, by
\bequation\label{eq.alpha_suff}
  \alphasuff_k \gets \min\left\{ \frac{2(1 - \eta)\beta_k\Delta l(x_k,\tau_k,g_k,d_k)}{(\tau_k L + \Gamma) \|d_k\|_2^2}, 1 \right\} \in \R{}_{>0}.
\eequation
The sequence $\{\beta_k\}$ referenced in \eqref{eq.alpha_suff} is chosen with different properties---namely, constant or diminishing---depending on the desired type of convergence guarantee.  We discuss details of the possible choices for $\{\beta_k\}$ and the consequences of these choices along with our convergence analysis.

Given that the step size $\alphasuff_k$ in \eqref{eq.alpha_suff} has been set based on a stochastic gradient estimate, a safeguard is needed for our convergence guarantees.  For this purpose, the second main idea of our step size selection strategy is to project the trial step size onto an interval that is appropriate depending on whether the search direction is tangentially dominated or normally dominated.  In particular, the step size is chosen as $\alpha_k \gets \proj_k(\alphasuff_k)$ where
\bequation\label{eq.alpha}
  \proj_k(\cdot) := \bcases \proj\(\cdot\ \Bigg| \left[ \displaystyle \frac{2(1 - \eta)\beta_k \xi_k \tau_k}{\tau_k L + \Gamma}, \frac{2(1 - \eta)\beta_k \xi_k \tau_k}{\tau_k L + \Gamma} + \theta\beta_k^2\right] \) & \text{if $\|u_k\|_2^2 \geq \chi_k \|v_k\|_2^2$} \\ \proj\(\cdot\ \Bigg| \left[ \displaystyle \frac{2(1 - \eta)\beta_k \xi_k}{\tau_k L + \Gamma}, \frac{2(1 - \eta)\beta_k \xi_k}{\tau_k L + \Gamma} + \theta\beta_k^2\right] \) & \text{otherwise.} \ecases
\eequation
Here, $\proj(\cdot | \Ical)$ denotes the projection onto the interval $\Ical \subset \R{}$.  In our analysis, the rules for $\{\beta_k\}$ (see Lemma~\ref{lem.key_decrease}) ensure that this projection only ever decreases the step size; hence, the overall motivation for the projection is to ensure that the step size is not too large compared to a conservative choice, namely, the lower end of the projection interval.  Motivation for the difference in the interval depending on whether the search direction is tangentially or normally dominated can be seen Lemma~\ref{lem.tau_zero_alpha_low} later on, where it is critical that the step size for a normally dominated search direction does not necessarily vanish if/when the merit parameter vanishes, i.e., $\{\tau_k\} \searrow 0$.

Overall, our step size selection mechanism can be understood as follows.  First, the algorithm adaptively sets the sequences $\{\chi_k\}$, $\{\zeta_k\}$, and $\{\xi_k\}$ in order to estimate bounds that are needed for the step size selection and are known to exist theoretically, but cannot be computed directly.  By the manner in which these sequences are set, our analysis shows that they remain constant for sufficiently large $k \in \N{}$ in any run of the algorithm.  With these values, our step size selection strategy aims to achieve a reduction in the merit function in expectation, with safeguards since the computed values are based on stochastic quantities.  One finds by the definition of the projection interval in \eqref{eq.alpha} that the step size \emph{for a tangentially dominated search direction} may decrease to zero if $\{\tau_k\} \searrow 0$; this is needed in cases when the problem is degenerate or infeasible, and the algorithm wants to avoid long steps in the tangential component that may ruin progress toward minimizing constraint violation.  Otherwise, \emph{for a normally dominated search direction}, the step size would remain bounded away from zero if $\beta_k = \beta \in (0,1]$ for all $k \in \N{}$; i.e., it can only decrease to zero if $\{\beta_k\}$ is diminishing.  If our algorithm did not make this distinction between the projection intervals for tangentially versus normally dominated search directions, then the algorithm would fail to have desirable convergence guarantees even in the deterministic setting.  (In particular, our proof in Appendix~\ref{sec.deterministic} of Theorem~\ref{th.deterministic}, which is upcoming in Section~\ref{sec.convergence}, would break down.)

Our complete algorithm is stated as Algorithm~\ref{alg.sqp} on page~\pageref{alg.sqp}.

\begin{algorithm}[ht]
  \caption{Stochastic SQP Algorithm}
  \label{alg.sqp}
  \begin{algorithmic}[1]
    \Require $L \in \R{}_{>0}$, a Lipschitz constant for $\nabla f$; $\Gamma \in \R{}_{>0}$, a Lipschitz constant for $c$; $\{\beta_k\} \subset (0,1]$; $x_0 \in \R{n}$; $\tau_{-1} \in \R{}_{>0}$; $\chi_{-1} \in \R{}_{>0}$; $\zeta_{-1} \in \R{}_{>0}$; $\omega \in \R{}_{>0}$; $\epsilon_v \in (0,1]$; $\sigma \in (0,1)$; $\epsilon_\tau \in (0,1)$; $\epsilon_\chi \in \R{}_{>0}$; $\epsilon_\zeta \in (0,1)$; $\epsilon_\xi \in (0,1)$; $\eta \in (0,1)$; $\theta \in \R{}_{\geq0}$
    \For{$k \in \N{}$}
    \If{$\|J_k^Tc_k\|_2 = 0$ and $\|c_k\|_2 > 0$}
        \State \textbf{terminate} and \textbf{return} $x_k$ (infeasible stationary point)
      \EndIf
      \State Compute a stochastic gradient $g_k$ \emph{at least} satisfying Assumption~\ref{ass.g}
      \State Compute $v_k \in \Range(J_k^T)$ that is feasible for problem~\eqref{prob.normal} and satisfies \eqref{eq.Cauchy}
      \State Compute $(u_k,y_k)$ as a solution of \eqref{eq.linsys}, and then set $d_k \gets v_k + u_k$ \label{line.linsys}
      \If{$d_k = 0$}
        \State Set $\tautrial_k \gets \infty$ and $\tau_k \gets \tau_{k-1}$
        \State Set $(\chi_k,\zeta_k) \gets (\chi_{k-1},\zeta_{k-1})$
        \State Set $\xitrial_k \gets \infty$ and $\xi_k \gets \xi_{k-1}$
        \State Set $\alphasuff_k \gets 1$ and $\alpha_k \gets 1$
      \Else
        \State Set $\tautrial_k$ by \eqref{eq.merit_parameter_trial} and $\tau_k$ by \eqref{eq.tau_update}
        \State Set $(\chi_k,\zeta_k)$ by \eqref{eq.switch_condition}--\eqref{eq.switch}
        \State Set $\xitrial_k$ by \eqref{eq.ratio_trial} and $\xi_k$ by \eqref{eq.ratio_update}
        \State Set $\alphasuff_k$ by \eqref{eq.alpha_suff} and $\alpha_k \gets \proj_k(\alphasuff_k)$ using \eqref{eq.alpha}
      \EndIf
      \State Set $x_{k+1} \gets x_k + \alpha_k d_k$
    \EndFor
  \end{algorithmic}
\end{algorithm}

%*********
% Section
%*********
\section{Convergence Analysis}\label{sec.convergence}

In this section, we prove convergence guarantees for Algorithm~\ref{alg.sqp}.  To understand the results that can be expected given our setting and the type of algorithm that we employ, let us first present a set of guarantees that can be proved if Algorithm~\ref{alg.sqp} were to be run with $g_k = \nabla f(x_k)$ and $\beta_k = \beta$ for all $k \in \N{}$, where $\beta \in \R{}_{>0}$ is sufficiently small.  For such an algorithm, we prove the following theorem in Appendix~\ref{sec.deterministic}.  The theorem is consistent with what can be proved for other deterministic algorithms in our context; e.g., see Theorem 3.3 in \cite{CurtNoceWaec09}.

\btheorem\label{th.deterministic}
  Suppose Algorithm~\ref{alg.sqp} is employed to solve problem~\eqref{prob.opt} such that Assumption~\ref{ass.main} holds, $g_k = \nabla f(x_k)$ for all $k \in \N{}$, $\{H_k\}$ satisfies Assumption~\ref{ass.H}, and $\beta_k = \beta$ for all $k \in \N{}$ where
  \bequation\label{eq.beta}
    \beta \in (0,1]\ \ \text{and}\ \ \frac{2 (1 - \eta) \beta \xi_{-1} \max\{\tau_{-1},1\}}{\Gamma} \in (0,1].
  \eequation
  If there exist $k_J \in \N{}$ and $\sigma_J \in \R{}_{>0}$ such that the singular values of $J_k$ are bounded below by $\sigma_J$ for all $k \geq k_J$, then the merit parameter sequence $\{\tau_k\}$ is bounded below by a positive real number and 
  \bequation\label{eq.det.limit1}
    0 = \lim_{k \to \infty} \left\|\bbmatrix \nabla f(x_k) + J_k^Ty_k \\ c_k \ebmatrix \right\|_2 = \lim_{k \to \infty} \left\| \bbmatrix Z_k^T \nabla f(x_k) \\ c_k \ebmatrix \right\|_2.
  \eequation
  Otherwise, if such $k_J$ and $\sigma_J$ do not exist, then it still follows that
  \bequation\label{eq.det.limit2}
    0 = \lim_{k \to \infty} \| J_k^Tc_k \|_2,
  \eequation
  and if $\{\tau_k\}$ is bounded below by a positive real number, then
  \bequation\label{eq.det.limit3}
    0 = \lim_{k \to \infty} \| \nabla f(x_k) + J_k^Ty_k \|_2 = \lim_{k \to \infty} \| Z_k^T\nabla f(x_k)\|_2.
  \eequation
\etheorem

Based on Theorem~\ref{th.deterministic}, the following aims---which are all achieved in certain forms in our analyses in Sections~\ref{sec.tau_low} and \ref{sec.tau_zero}---can be set for Algorithm~\ref{alg.sqp} in the stochastic setting.  First, if Algorithm~\ref{alg.sqp} is run and the singular values of the constraint Jacobians happen to remain bounded away from zero beyond some iteration, then (following~\eqref{eq.det.limit1}) one should aim to prove that a primal-dual stationarity measure (recall~\eqref{eq.KKT}) vanishes in expectation.  This is shown under certain conditions in Corollary~\ref{cor.stochastic_tau_finite_large_pi_finite} (and the subsequent discussion) on page~\pageref{cor.stochastic_tau_finite_large_pi_finite}.  Otherwise, a (sub)sequence of $\{J_k\}$ tends to singularity, in which case (following \eqref{eq.det.limit2}) one should at least aim to prove that $\{\|J_k^Tc_k\|_2\}$ vanishes in expectation, which would mean that a (sub)sequence of iterates converges in expectation to feasibility or at least stationarity with respect to the constraint infeasibility measure $\varphi$ (recall \eqref{eq.infeasible_stationary}).  Such a conclusion is offered under certain conditions by combining Corollary~\ref{cor.stochastic_tau_finite_large_pi_finite} (see page~\pageref{cor.stochastic_tau_finite_large_pi_finite}) and Theorem~\ref{th.tau_zero} (see page~\pageref{th.tau_zero}).  The remaining aim (paralleling~\eqref{eq.det.limit3}) is that one should aim to prove that even if a (sub)sequence of $\{J_k\}$ tends to singularity, if the merit parameter sequence~$\{\tau_k\}$ happens to remain bounded below by a positive real number, then $\{\|Z_k^T\nabla f(x_k)\|_2\}$ vanishes in expectation.  This can also be seen to occur under certain conditions in Corollary~\ref{cor.stochastic_tau_finite_large_pi_finite} on page~\pageref{cor.stochastic_tau_finite_large_pi_finite}.

In addition, due to its stochastic nature, there are events that one should consider in which the algorithm may exhibit behavior that cannot be exhibited by the deterministic algorithm.  One such event is when the merit parameter eventually remains fixed at a value that is not sufficiently small.  We show in Section~\ref{sec.tau_big}---with formal results stated and proved in Appendix~\ref{sec.totalprob}---that, under reasonable assumptions, the total probability of this event (over all possible runs of the algorithm) is zero.  We complete the picture of the possible behaviors of our algorithm by discussing remaining possible (practically irrelevant) events in Section~\ref{sec.diverge}.

Let us now commence our analysis of Algorithm~\ref{alg.sqp}.  If a run terminates finitely at iteration $k \in \N{}$, then an infeasible stationary point has been found.  Hence, without loss of generality throughout the remainder of our analysis and discussions, we assume that the algorithm does not terminate finitely, i.e., an infinite number of iterates are generated.  As previously mentioned, for much of our analysis, we merely assume that the stochastic gradient estimates satisfy Assumption~\ref{ass.g}.  This is done to show that many of our results hold under this general setting.  However, we will ultimately impose stronger conditions on $\{g_k\}$, as needed; see Sections~\ref{sec.tau_zero} and \ref{sec.tau_big} (and Appendix~\ref{sec.totalprob}).

We build to our main results through a series of lemmas.  Our first lemma has appeared for various deterministic algorithms in the literature.  It extends easily to our setting since the normal component computation is deterministic conditioned on the event that the algorithm reaches $x_k$.

\blemma\label{lem.normal_component}
  There exist $\kappa_v \in \R{}_{>0}$ and $\underline\omega \in \R{}_{>0}$ such that, in any run of the algorithm,
  \bsubequations
    \begin{align*}
      \|c_k\|_2(\|c_k\|_2 - \|c_k + J_k v_k\|_2) &\geq \kappa_v \|J_k^Tc_k\|_2^2 \\%\label{eq.normal_component.Cauchy} \\
      \text{and}\ \ \underline\omega \|J_k^Tc_k\|_2^2 \leq \|v_k\|_2 &\leq \omega \|J_k^Tc_k\|_2\ \ \text{for all}\ \ k \in \N{}\ \text{with}\ \|c_k\|_2 > 0. %\label{eq.normal_component.v_bound}
    \end{align*}
  \esubequations
\elemma
\proof{Proof.}
  The proof follows as for Lemmas~3.5 and 3.6 in \cite{CurtNoceWaec09}. \Halmos
\endproof

Our second lemma shows that the procedure for setting $\{\chi_k\}$ and $\{\zeta_k\}$ guarantees that these sequences are constant \emph{deterministically} for sufficiently large $k \in \N{}$.

\blemma\label{lem.chi_upper}
  There exist $(\chi_{\max},\zeta_{\min}) \in \R{}_{>0} \times \R{}_{>0}$ such that, in any run, there exists $k_{\chi,\zeta} \in \N{}$ such that $(\chi_k,\zeta_k) = (\chi_{k_{\chi,\zeta}},\zeta_{k_{\chi,\zeta}})$ for all $k \geq k_{\chi,\zeta}$, where $(\chi_{k_{\chi,\zeta}},\zeta_{k_{\chi,\zeta}}) \in (0,\chi_{\max}] \times [\zeta_{\min},\infty)$.
\elemma
\proof{Proof.}
  Consider arbitrary $k \in \N{}$ in any run.  If $d_k = 0$, then the algorithm sets $(\chi_k,\zeta_k) = (\chi_{k-1},\zeta_{k-1})$.  Otherwise, under Assumption~\ref{ass.H}, it follows for any $\chi \in \R{}_{>0}$ that $\|u_k\|_2^2 \geq \chi \|v_k\|_2^2$ implies
  \bequationNN
    \baligned
      \thalf d_k^TH_kd_k
        &= \thalf u_k^TH_ku_k + u_k^TH_kv_k + \thalf v_k^TH_kv_k \\
        &\geq \thalf \zeta \|u_k\|_2^2 - \|u_k\|_2 \|H_k\|_2 \|v_k\|_2 - \thalf \|H_k\|_2 \|v_k\|_2^2 \geq \(\frac{\zeta}{2} - \frac{\kappa_H}{\sqrt{\chi}} - \frac{\kappa_H}{2\chi}\) \|u_k\|_2^2.
    \ealigned
  \eequationNN
  Hence, for sufficiently large $\chi \in \R{}_{>0}$, one finds that $\|u_k\|_2^2 \geq \chi \|v_k\|_2^2$ implies $\thalf d_k^TH_kd_k \geq \tfrac14 \zeta \|u_k\|_2^2$.  The conclusion follows from this fact and the procedure for setting $(\chi_k,\zeta_k)$ in \eqref{eq.switch_condition}--\eqref{eq.switch}. \Halmos
\endproof

We now prove that the sequence $\{\xi_k\}$ is bounded below \emph{deterministically}.

\blemma\label{lem.xi_lower}
  There exists $\xi_{\min} \in \R{}_{>0}$ such that, in any run of the algorithm, there exists $k_{\xi} \in \N{}$ such that $\xi_k = \xi_{k_{\xi}}$ for all $k \geq k_{\xi}$, where $\xi_{k_{\xi}} \in [\xi_{\min},\infty)$.
\elemma
\proof{Proof.}
  Consider arbitrary $k \in \N{}$ in any run.  If $d_k = 0$, then the algorithm sets $\xi_k = \xi_{k-1}$.  If $d_k \neq 0$ and $\|u_k\|_2^2 \geq \chi_k \|v_k\|_2^2$, then it follows from \eqref{eq.model_reduction}--\eqref{eq.model_reduction_condition} and \eqref{eq.ratio_trial}--\eqref{eq.ratio_update} that either $\xi_k = \xi_{k-1}$ or
  \bequationNN
    \xi_k \geq (1 - \epsilon_\xi) \xitrial_k = (1 - \epsilon_\xi) \(\frac{\Delta l(x_k,\tau_k,g_k,d_k)}{\tau_k\|d_k\|_2^2}\) \geq (1 - \epsilon_\xi) \frac{\tau_k \zeta \|u_k\|_2^2}{\tau_k(1 + \chi_k^{-1})\|u_k\|_2^2} \geq (1 - \epsilon_\xi) \frac{\zeta}{(1 + \chi_{-1}^{-1})}.
  \eequationNN
  If $d_k \neq 0$ and $\|u_k\|_2^2 < \chi_k \|v_k\|_2^2$, then by \eqref{eq.model_reduction}--\eqref{eq.model_reduction_condition}, \eqref{eq.ratio_trial}--\eqref{eq.ratio_update}, and Lemmas~\ref{lem.normal_component}--\ref{lem.chi_upper} either $\xi_k = \xi_{k-1}$ or
  \bequationNN
    \xi_k \geq (1 - \epsilon_\xi) \xitrial_k = (1 - \epsilon_\xi) \(\frac{\Delta l(x_k,\tau_k,g_k,d_k)}{\|d_k\|_2^2}\) \geq (1 - \epsilon_\xi) \frac{\sigma \kappa_v \kappa_c^{-1} \|J_k^Tc_k\|_2^2}{(\chi_k + 1)\omega^2 \|J_k^Tc_k\|_2^2} \geq (1 - \epsilon_\xi) \frac{\sigma \kappa_v \kappa_c^{-1}}{(\chi_{\max} + 1)\omega^2}.
  \eequationNN
  Combining these results, the desired conclusion follows. \Halmos
\endproof

Our next two lemmas provide useful relationships between deterministic and stochastic quantities conditioned on the event that the algorithm has reached $x_k$ as the $k$th iterate.  The first result is similar to \cite[Lemma 3.6]{BeraCurtRobiZhou21}, although the proof presented here is different in order to handle potential rank deficiency of the constraint Jacobians.  Here and throughout the remainder of our analysis, conditioned on the event that the algorithm reaches $x_k$ as the $k$th iterate, we denote~$\utrue_k$ as the tangential component of the search direction that would be computed if $\nabla f(x_k)$ were used in place of $g_k$ in~\eqref{eq.linsys}, and similarly denote $\dtrue_k := v_k + \utrue_k$.

\blemma\label{lem.expectation}
  For all $k \in \N{}$ in any run, $\E_k[u_k] = \utrue_k$ and $\E_k[\|d_k - \dtrue_k\|_2] \leq \zeta^{-1} \sqrt{M}$.
\elemma
\proof{Proof.}
  Consider arbitrary $k \in \N{}$ in any run.  Under Assumption~\ref{ass.H}, it follows from \eqref{eq.linsys} that there exist $w_k$ and $\wtrue_k$ such that $u_k = Z_kw_k$ and $\utrue_k = Z_k\wtrue_k$ where $w_k = -(Z_k^TH_kZ_k)^{-1}Z_k^T(g_k + H_kv_k)$ and $\wtrue_k = -(Z_k^TH_kZ_k)^{-1}Z_k^T(\nabla f(x_k) + H_kv_k)$.  Since $(Z_k^TH_kZ_k)^{-1}Z_k^T$ and $Z_k$ are linear operators, it follows that $\E_k[w_k] = \wtrue_k$ and hence $\E_k[u_k] = \utrue_k$, as desired.  Then, it follows from consistency and submultiplicity of the spectral norm, orthonormality of $Z_k$, Jensen's inequality, concavity of the square root operator, and Assumptions~\ref{ass.g} and~\ref{ass.H} that
  \bequationNN
    \baligned
      \E_k[\|d_k - \dtrue_k\|_2] = \E_k[\|u_k - \utrue_k\|_2] &= \E_k[\|Z_k(w_k - \wtrue_k)\|_2] \\
      &= \E_k[\|Z_k(Z_k^TH_kZ_k)^{-1}Z_k^T(g_k - \nabla f(x_k))\|_2] \\
      &\leq \E_k[\|Z_k(Z_k^TH_kZ_k)^{-1}Z_k^T\|_2 \|g_k - \nabla f(x_k)\|_2] \\
      &= \|(Z_k^TH_kZ_k)^{-1}\|_2 \E_k[ \|g_k - \nabla f(x_k)\|_2] \\
      &\leq \zeta^{-1} \E_k[\|g_k - \nabla f(x_k)\|_2] \\
      &\leq \zeta^{-1} \sqrt{\E_k[\|g_k - \nabla f(x_k)\|_2^2]} \leq \zeta^{-1} \sqrt{M},
    \ealigned
  \eequationNN
  which is the final desired conclusion. \Halmos
\endproof

Our next result is part of \cite[Lemma 3.9]{BeraCurtRobiZhou21}; we provide a proof for completeness.

\blemma\label{lem.product_bounds}
  For all $k \in \N{}$ in any run, $\nabla f(x_k)^T\dtrue_k \geq \E_k[g_k^Td_k] \geq \nabla f(x_k)^T\dtrue_k - \zeta^{-1} M$.
\elemma
\proof{Proof.}
  Consider arbitrary $k \in \N{}$ in any run.  The arguments in the proof of Lemma~\ref{lem.expectation} give
  \bequationNN
    \baligned
      g_k^Tu_k &= -g_k^TZ_k(Z_k^TH_kZ_k)^{-1}Z_k^T(g_k + H_kv_k) \\ \text{and}\ \ \nabla f(x_k)^T\utrue_k &= -\nabla f(x_k)^TZ_k(Z_k^TH_kZ_k)^{-1}Z_k^T(\nabla f(x_k) + H_kv_k).
    \ealigned
  \eequationNN
  On the other hand, under Assumptions~\ref{ass.g} and \ref{ass.H}, it follows that
  \bequationNN
    \zeta^{-1}M \geq \E_k[\|Z_k^T(g_k - \nabla f(x_k))\|_{(Z_k^TH_kZ_k)^{-1}}^2] \geq 0,
  \eequationNN
  where
  \bequationNN
    \baligned
      &\ \E_k[\|Z_k^T(g_k - \nabla f(x_k))\|_{(Z_k^TH_kZ_k)^{-1}}^2] \\
      =&\ \E_k[\|Z_k^Tg_k\|_{(Z_k^TH_kZ_k)^{-1}}^2] - 2\E_k[g_k^TZ_k(Z_k^TH_kZ_k)^{-1}Z_k^T\nabla f(x_k)] + \|Z_k^T\nabla f(x_k)\|_{(Z_k^TH_kZ_k)^{-1}}^2 \\
      =&\ \E_k[\|Z_k^Tg_k\|_{(Z_k^TH_kZ_k)^{-1}}^2] - \|Z_k^T\nabla f(x_k)\|_{(Z_k^TH_kZ_k)^{-1}}^2.
    \ealigned
  \eequationNN
  Combining the facts above and again using Assumption~\ref{ass.g}, it follows that
  \bequationNN
    \baligned
      \nabla f(x_k)^T\dtrue_k - \E_k[g_k^Td_k]
        =&\ \nabla f(x_k)^Tv_k + \nabla f(x_k)^T\utrue_k - \E_k[g_k^Tv_k + g_k^Tu_k] \\
        =&\ \nabla f(x_k)^T\utrue_k - \E_k[g_k^Tu_k] \\
        =&\ -\nabla f(x_k)^TZ_k(Z_k^TH_kZ_k)^{-1}Z_k^T(\nabla f(x_k) + H_kv_k) \\
        &\ + \E_k[g_k^TZ_k(Z_k^TH_kZ_k)^{-1}Z_k^T(g_k + H_kv_k)] \\
        =&\ -\|Z_k^T\nabla f(x_k)\|_{(Z_k^TH_kZ_k)^{-1}}^2 + \E_k[\|Z_k^Tg_k\|_{(Z_k^TH_kZ_k)^{-1}}^2] \in [0,\zeta^{-1}M],
    \ealigned
  \eequationNN
  which gives the desired conclusion.
  \Halmos
\endproof

In the subsequent subsections, our analysis turns to offering guarantees conditioned on each of a few possible events that occur in a run of the algorithm, a few of which involve the merit parameter sequence eventually remaining constant in a run of the algorithm.  Before considering these events, let us first prove under certain circumstances that such behavior of the merit parameter sequence would occur.  As seen in Theorem~\ref{th.deterministic}, it is worthwhile to consider such an occurrence regardless of the properties of the sequence of constraint Jacobians.  That said, one might only be able to prove that it occurs when the constraint Jacobians are eventually bounded away from singularity.

Our first lemma here proves that if the constraint Jacobians are eventually bounded away from singularity, then the normal components of the search directions satisfy a useful upper bound.  The proof is essentially the same as that of~\cite[Lemma~3.15]{CurtNoceWaec09}, but we provide it for completeness.

\blemma\label{lem.det.v_upper}
  If, in a run of the algorithm, there exist $k_J \in \N{}$ and $\sigma_J \in \R{}_{>0}$ such that the singular values of $J_k$ are bounded below by $\sigma_J$ for all $k \geq k_J$, then there exists $\kappa_\omega \in \R{}_{>0}$ such that
  \bequationNN
    \|v_k\|_2 \leq \kappa_\omega (\|c_k\|_2 - \|c_k + J_kv_k\|_2)\ \ \text{for all}\ \ k \geq k_J.
  \eequationNN
\elemma
\proof{Proof.}
  Under the conditions of the lemma, $\|J_k^Tc_k\|_2 \geq \sigma_J \|c_k\|_2$ for all $k \geq k_J$.  Hence, along with Lemma~\ref{lem.normal_component}, it follows that $\|c_k\|_2(\|c_k\|_2 - \|c_k + J_kv_k\|_2) \geq \kappa_v \|J_k^Tc_k\|_2^2 \geq \kappa_v \sigma_J^2 \|c_k\|_2^{2}$ for all $k \geq k_J$.  Combining this again with Lemma~\ref{lem.normal_component}, it follows with the Cauchy-Schwarz inequality and \eqref{eq.bounds} that
  \bequationNN
    \|v_k\|_2 \leq \omega \|J_k^T\|_2 \|c_k\|_2 \leq \frac{\omega \kappa_J}{\kappa_v \sigma_J^2} (\|c_k\|_2 - \|c_k + J_kv_k\|_2)\ \ \text{for all}\ \ k \geq k_J,
  \eequationNN
  from which the desired conclusion follows. \Halmos
\endproof

We now prove that if the differences between the stochastic gradient estimates and the true gradients are bounded deterministically, then the sequence of tangential components is bounded.

\blemma\label{lem.det.u_bound}
  If, in a run of the algorithm, the sequence $\{\|g_k - \nabla f(x_k)\|_2\}$ is bounded by a positive real number $\kappa_g \in \R{}_{>0}$, then the sequence $\{\|u_k\|_2\}$ is bounded by a positive real number $\kappa_u \in \R{}_{>0}$.
\elemma
\proof{Proof.}
  Under Assumption~\ref{ass.main}, the sequence $\{\|\nabla f(x_k)\|_2\}$ is bounded; recall \eqref{eq.bounds}.  Hence, under the conditions of the lemma, $\{\|g_k\|_2\}$ is bounded.  The first block of~\eqref{eq.linsys} yields $u_k^TH_ku_k = -u_k^T(g_k + H_kv_k)$, which under Assumption~\ref{ass.H} yields $\zeta\|u_k\|_2^2 \leq -u_k^Tg_k - u_k^TH_kv_k \leq (\|g_k\|_2 + \|H_k\|_2 \|v_k\|_2) \|u_k\|_2$.  Hence, the conclusion follows from these facts, Assumption~\ref{ass.main}, Assumption~\ref{ass.H}, and Lemma~\ref{lem.normal_component}. \Halmos
\endproof

By combining the preceding two lemmas, the following lemma indicates certain circumstances under which the sequence of merit parameters will eventually remain constant.  We remark that it is possible in a run of the algorithm for the merit parameter sequence to remain constant eventually even if the conditions of the lemma do not hold, which is why our analyses in the subsequent subsections do not presume that these conditions hold.  That said, to prove that there exist settings in which the merit parameter is guaranteed to remain constant eventually, we offer the following.

\blemma\label{lem.det.pi_bounded}
  If, in a run, there exist $k_J \in \N{}$ and $\sigma_J \in \R{}_{>0}$ such that the singular values of $J_k$ are bounded below by $\sigma_J$ for all $k \geq k_J$ and $\{\|g_k - \nabla f(x_k)\|_2\}$ is bounded by a positive real number $\kappa_g \in \R{}_{>0}$, then there exist $k_\tau \in \N{}$ and $\tau_{\min} \in \R{}_{>0}$ such that $\tau_k = \tau_{\min}$ for all $k \in \N{}$ with $k \geq k_\tau$.
\elemma
\proof{Proof.}
  Observe that the algorithm terminates if $\|J_k^Tc_k\|_2 = 0$ while $\|c_k\|_2 > 0$.  Let us now show that if $\|c_k\|_2 = 0$, then the algorithm sets $\tau_k \gets \tau_{k-1}$.  Indeed, $\|c_k\|_2 = 0$ implies $v_k = 0$ by Lemma~\ref{lem.normal_component}.  If $u_k = 0$ as well, then $d_k = 0$ and the algorithm explicitly sets $\tau_k \gets \tau_{k-1}$.  Otherwise, if $v_k = 0$ and $u_k \neq 0$, then \eqref{eq.linsys} yields $0=g_k^Tu_k + u_k^TH_ku_k = g_k^Td_k + u_k^TH_ku_k$, in which case \eqref{eq.merit_parameter_trial}--\eqref{eq.tau_update} again yield $\tau_k \gets \tau_{k-1}$.  Overall, it follows that $\tau_k < \tau_{k-1}$ if and only if one finds $\|J_k^Tc_k\|_2 > 0$, $g_k^Td_k + u_k^TH_ku_k > 0$, and $\tau_{k-1} (g_k^Td_k + u_k^TH_ku_k) > (1 - \sigma)(\|c_k\|_2 - \|c_k + J_kv_k\|_2)$.  On the other hand, from the first equation in \eqref{eq.linsys}, the Cauchy-Schwarz inequality, \eqref{eq.bounds}, and Lemmas~\ref{lem.det.v_upper} and \ref{lem.det.u_bound}, it holds that
  \bequationNN
    \baligned
      g_k^Td_k + u_k^TH_ku_k = (g_k - H_ku_k)^Tv_k
        &= (g_k - \nabla f(x_k) + \nabla f(x_k) - H_ku_k)^Tv_k \\
        &\leq (\kappa_g + \kappa_{\nabla f} + \kappa_H \kappa_u) \|v_k\|_2 \\
        &\leq (\kappa_g + \kappa_{\nabla f} + \kappa_H \kappa_u) \kappa_\omega (\|c_k\|_2 - \|c_k + J_kv_k\|_2).
    \ealigned
  \eequationNN
  Combining these facts, the desired conclusion follows. \Halmos
\endproof

%************
% Subsection
%************
\subsection{Constant, Sufficiently Small Merit Parameter}\label{sec.tau_low}

Our goal in this subsection is to prove a convergence guarantee for our algorithm in the event $E_{\tau,\low}$, which is defined formally in the assumption below.  In the assumption, similar to our notation of $\utrue_k$ and $\dtrue_k$, we use $\tautruetrial_k$ to denote the value of $\tautrial_k$ that, conditioned on $x_k$ as the $k$th iterate, would be computed in iteration $k \in \N{}$ if the search direction were computed using the true gradient $\nabla f(x_k)$ in place of $g_k$ in \eqref{eq.linsys}.

\bassumption\label{ass.tau_low}
  In a run of the algorithm, event $E_{\tau,\low}$ occurs, i.e., there exists an iteration number $k_\tau \in \N{}$ and a merit parameter value $\tau_{\min} \in \R{}_{>0}$ such that
  \bequationNN
    \tau_k = \tau_{\min} \leq \tautruetrial_k,\ \ \chi_k = \chi_{k-1},\ \ \zeta_k = \zeta_{k-1},\ \ \text{and}\ \ \xi_k = \xi_{k-1}\ \ \text{for all}\ \ k \geq k_\tau.
  \eequationNN
  In addition, along the lines of Assumption~\ref{ass.g}, $\{g_k\}_{k \geq k_\tau}$ satisfies $\E_{k,\tau,\low}[g_k] = \nabla f(x_k)$ and $\E_{k,\tau,\low}[\|g_k - \nabla f(x_k)\|_2^2] \leq M$, where $\E_{k,\tau,\low}$ denotes expectation with respect to the distribution of~$\iota$ conditioned on the event that $E_{\tau,\low}$ occurs and the algorithm has reached $x_k$ in iteration $k \in \N{}$.
\eassumption

Recall from Lemmas~\ref{lem.chi_upper} and \ref{lem.xi_lower} that the sequences $\{\chi_k\}$, $\{\zeta_k\}$, and $\{\xi_k\}$ are guaranteed to be bounded deterministically, and in particular will remain constant for sufficiently large $k \in \N{}$.  Hence, one circumstance in which Assumption~\ref{ass.tau_low} may hold is under the conditions of Lemma~\ref{lem.det.pi_bounded}.  A critical distinction in Assumption~\ref{ass.tau_low} is that the value at which the merit parameter eventually settles is sufficiently small such that $\tau_k \leq \tautruetrial_k$ for all sufficiently large $k \in \N{}$.  This is the key distinction between the event $E_{\tau,\low}$ and some of the events we consider in Sections~\ref{sec.tau_big} and \ref{sec.diverge}.

For the sake of brevity in the rest of this subsection, let us temporarily redefine $\E_k := \E_{k,\tau,\low}$.

Our next lemma provides a key result that drives our analysis for this subsection.  It shows that as long as $\beta_k$ is sufficiently small for all $k \in \N{}$ (in a manner similar to \eqref{eq.beta}), the reduction in the merit function in each iteration is at least the sum of two terms: (1) the reduction in the model of the merit function corresponding to the \emph{true} gradient and its associated search direction, and (2) a pair of quantities that can be attributed to the error in the stochastic gradient estimate.

\blemma\label{lem.key_decrease}
  Suppose that $\{\beta_k\}$ is chosen such that
  \bequation\label{eq.key_decrease_beta}
    \beta_k \in (0,1]\ \ \text{and}\ \ \frac{2 (1 - \eta) \beta_k \xi_k \max\{\tau_k,1\}}{\tau_k L + \Gamma} \in (0,1]\ \ \text{for all}\ \ k \in \N{}.
  \eequation
  Then, for all $k \in \N{}$ in any such run of the algorithm, it follows that
  \bequationNN
    \baligned
      &\ \phi(x_k, \tau_k) - \phi(x_k + \alpha_k d_k, \tau_k) \\
      \geq&\ \alpha_k \Delta l(x_k,\tau_k,\nabla f(x_k),\dtrue_k) - (1 - \eta) \alpha_k \beta_k  \Delta l(x_k,\tau_k,g_k,d_k) - \alpha_k \tau_k \nabla f(x_k)^T (d_k - \dtrue_k).
    \ealigned
  \eequationNN
\elemma
\proof{Proof.}
  Consider arbitrary $k \in \N{}$ in any run.  From \eqref{eq.alpha_suff}--\eqref{eq.alpha} and the supposition about~$\{\beta_k\}$, one finds $\alpha_k \in (0,1]$.  Hence, with \eqref{eq.merit_decrease} and $J_kd_k = J_k\dtrue_k$ (since $J_ku_k = J_k\utrue_k = 0$ by \eqref{eq.linsys}), one has
  \begin{align}
    &\ \phi(x_k, \tau_k) - \phi(x_k + \alpha_k d_k, \tau_k) \nonumber \\
    \geq&\ -\alpha_k (\tau_k \nabla f(x_k)^T d_k - \|c_k\|_2 + \|c_k + J_kd_k\|_2) - \thalf (\tau_k L + \Gamma) \alpha_k^2 \|d_k\|_2^2 \nonumber \\
    =&\ -\alpha_k (\tau_k \nabla f(x_k)^T \dtrue_k - \|c_k\|_2 + \|c_k + J_k\dtrue_k\|_2) - \thalf (\tau_k L + \Gamma) \alpha_k^2 \|d_k\|_2^2 - \alpha_k \tau_k \nabla f(x_k)^T (d_k - \dtrue_k) \nonumber \\
    =&\ \alpha_k \Delta l(x_k,\tau_k,\nabla f(x_k),\dtrue_k) - \thalf (\tau_k L + \Gamma) \alpha_k^2 \|d_k\|_2^2 - \alpha_k \tau_k \nabla f(x_k)^T (d_k - \dtrue_k). \label{eq.merit_reduction_lower}
  \end{align}
  By \eqref{eq.alpha_suff}, it follows that $\alphasuff_k \leq \tfrac{2(1-\eta)\beta_k\Delta l(x_k,\tau_k,g_k,d_k)}{(\tau_k L + \Gamma) \|d_k\|_2^2}$. If $\|u_k\|_2^2 \geq \chi_k \|v_k\|_2^2$, then it follows from \eqref{eq.ratio_trial}--\eqref{eq.ratio_update} that $\xi_k \leq \xitrial_k = \frac{\Delta l(x_k,\tau_k,g_k,d_k)}{\tau_k\|d_k\|_2^2}$ and $\tfrac{2(1-\eta)\beta_k\Delta l(x_k,\tau_k,g_k,d_k)}{(\tau_k L + \Gamma) \|d_k\|_2^2} \geq \frac{2(1-\eta)\beta_k\xi_k\tau_k}{\tau_k L + \Gamma}$.  On the other hand, if $\|u_k\|_2^2 < \chi_k \|v_k\|_2^2$, then it follows from \eqref{eq.ratio_trial}--\eqref{eq.ratio_update} that $\xi_k \leq \xitrial_k = \frac{\Delta l(x_k,\tau_k,g_k,d_k)}{\|d_k\|_2^2}$ and $\tfrac{2(1-\eta)\beta_k\Delta l(x_k,\tau_k,g_k,d_k)}{(\tau_k L + \Gamma) \|d_k\|_2^2} \geq \frac{2(1-\eta)\beta_k\xi_k}{\tau_k L + \Gamma}$.  It follows from these facts and the supposition about $\{\beta_k\}$ that the projection in \eqref{eq.alpha} never sets $\alpha_k > \alphasuff_k$. Thus, $\alpha_k \leq \alphasuff_k \leq \tfrac{2(1-\eta)\beta_k\Delta l(x_k,\tau_k,g_k,d_k)}{(\tau_k L + \Gamma) \|d_k\|_2^2}$.  Hence, by \eqref{eq.merit_reduction_lower},
  \bequationNN
      \baligned
        &\ \phi(x_k, \tau_k) - \phi(x_k + \alpha_k d_k, \tau_k) \\
        \geq&\ \alpha_k \Delta l(x_k,\tau_k,\nabla f(x_k),\dtrue_k) - \thalf \alpha_k (\tau_k L + \Gamma) \(\tfrac{2(1-\eta)\beta_k \Delta l(x_k,\tau_k,g_k,d_k)}{(\tau_k L + \Gamma) \|d_k\|_2^2}\) \|d_k\|_2^2 - \alpha_k \tau_k \nabla f(x_k)^T (d_k - \dtrue_k) \\
        =&\ \alpha_k \Delta l(x_k,\tau_k,\nabla f(x_k),\dtrue_k) - (1 - \eta) \alpha_k \beta_k \Delta l(x_k,\tau_k,g_k,d_k) - \alpha_k \tau_k \nabla f(x_k)^T (d_k - \dtrue_k), \\
      \ealigned
    \eequationNN
  which completes the proof.\Halmos
\endproof

Our second result in this case offers a critical upper bound on the final term in the conclusion of Lemma~\ref{lem.key_decrease}.  The result follows in a similar manner as \cite[Lemma~3.11]{BeraCurtRobiZhou21}.

\blemma\label{lem.alphagd}
  For any run under Assumption~\ref{ass.tau_low}, it follows for any $k \geq k_\tau$ that
  \bequationNN
    \E_k[\alpha_k \tau_k \nabla f(x_k)^T (d_k - \dtrue_k)] \leq \beta_k^2 \theta \tau_{\min} \kappa_{\nabla f} \zeta^{-1} \sqrt{M}.
  \eequationNN
\elemma
\proof{Proof.}
  Consider $k \geq k_\tau$, where $k_\tau$ is defined in Assumption~\ref{ass.tau_low}.  We prove the desired conclusion under the assumption that the search direction in iteration $k$ is tangentially dominated, then argue that it also holds by a similar argument when this search direction is normally dominated.  Let $I_k$ be the event that $\nabla f(x_k)^T(d_k - \dtrue_k) \geq 0$ and let $I_k^c$ be the complementary event.  In addition, let $\P_k$ denote probability conditioned on the event that $E_{\tau,\low}$ occurs and $x_k$ is the $k$th iterate.  By the law of total expectation, Assumption~\ref{ass.tau_low}, and \eqref{eq.alpha}, one finds that
  \bequationNN
    \baligned
      &\ \E_k[\alpha_k \tau_k \nabla f(x_k)^T(d_k - \dtrue_k)] \\
      =&\ \E_k[\alpha_k \tau_{\min} \nabla f(x_k)^T(d_k - \dtrue_k) | I_k] \P_k[I_k] + \E_k[\alpha_k \tau_{\min} \nabla f(x_k)^T(d_k - \dtrue_k) | I_k^c] \P_k[I_k^c] \\
      \leq& \alpha_{k,\max} \tau_{\min} \E_k[\nabla f(x_k)^T(d_k - \dtrue_k) | I_k] \P_k[I_k] + \alpha_{k,\min} \tau_{\min} \E_k[\nabla f(x_k)^T(d_k - \dtrue_k) | I_k^c] \P_k[I_k^c],
    \ealigned
  \eequationNN
  where $\alpha_{k,\min} := \frac{2(1-\eta)\beta_k\hat\xi_{\min}\tau_{\min}}{\tau_{\min} L + \Gamma}$ and $\alpha_{k,\max} := \frac{2(1-\eta)\beta_k\hat\xi_{\min}\tau_{\min}}{\tau_{\min} L + \Gamma} + \theta \beta_k^2$ are, respectively, the lower and upper bounds for the step size for the tangentially dominated search direction from \eqref{eq.alpha} with $\hat\xi_{\min} \in [\xi_{\min},\infty)$ being the positive real number such that $\xi_k = \hat\xi_{\min}$ for all $k \geq k_\tau$ (see Lemma~\ref{lem.xi_lower} and Assumption~\ref{ass.tau_low}).  Thus, since $\E_k[d_k] = \dtrue_k$ by Lemma~\ref{lem.expectation}, the law of total expectation yields
  \bequationNN
    \baligned
      &\ \E_k[\alpha_k \tau_k \nabla f(x_k)^T(d_k - \dtrue_k)] \\
      \leq&\ \alpha_{k,\min} \tau_{\min} \E_k[\nabla f(x_k)^T(d_k - \dtrue_k) | I_k] \P_k[I_k] + \alpha_{k,\min} \tau_{\min} \E_k[\nabla f(x_k)^T(d_k - \dtrue_k) | I_k^c] \P_k[I_k^c] \\
      &\ + (\alpha_{k,\max} - \alpha_{k,\min}) \tau_{\min} \E_k[\nabla f(x_k)^T(d_k - \dtrue_k) | I_k] \P_k[I_k] \\
      =&\ (\alpha_{k,\max} - \alpha_{k,\min}) \tau_{\min} \E_k[\nabla f(x_k)^T(d_k - \dtrue_k) | I_k] \P_k[I_k].
    \ealigned
  \eequationNN
  Moreover, by the Cauchy-Schwarz inequality and law of total expectation, one finds
  \bequationNN
    \baligned
      &\ \E_k[\nabla f(x_k)^T(d_k - \dtrue_k) | I_k] \P_k[I_k] \\
      \leq&\ \E_k[\|\nabla f(x_k)\|_2 \|d_k - \dtrue_k\|_2 | I_k] \P_k[I_k] \\
      =&\ \E_k[\|\nabla f(x_k)\|_2 \|d_k - \dtrue_k\|_2] - \E_k[\|\nabla f(x_k)\|_2 \|d_k - \dtrue_k\|_2 | I_k^c] \P_k[I_k^c] \\
      \leq&\ \|\nabla f(x_k)\|_2 \E_k[\|d_k - \dtrue_k\|_2].
    \ealigned
  \eequationNN
  Combining the above results, \eqref{eq.bounds}, Lemma~\ref{lem.expectation}, and the fact that $\alpha_{k,\max} - \alpha_{k,\min} = \theta \beta_k^2$, the desired conclusion follows for tangentially dominated search directions.  Finally, using the same arguments---except with $\alpha_{k,\min} := \frac{2(1-\eta)\beta_k\hat\xi_{\min}}{\tau_{\min} L + \Gamma}$ and $\alpha_{k,\max} := \frac{2(1-\eta)\beta_k\hat\xi_{\min}}{\tau_{\min} L + \Gamma} + \theta \beta_k^2$, where again $\alpha_{k,\max} - \alpha_{k,\min} = \theta \beta_k^2$---the desired conclusion follows for normally dominated search directions as well. \Halmos
\endproof

Our next result in this case bounds the middle term in the conclusion of Lemma~\ref{lem.key_decrease}. 

\blemma\label{lem.Deltaq}
  For any run under Assumption~\ref{ass.tau_low}, it follows for any $k \geq k_\tau$ that 
  \bequationNN
    \E_k[\Delta l(x_k,\tau_{\min},g_k,d_k)] \leq \Delta l(x_k,\tau_{\min},\nabla f(x_k),\dtrue_k) + \tau_{\min} \zeta^{-1} M.
  \eequationNN
\elemma
\proof{Proof.}
  Consider arbitrary $k \geq k_\tau$ in any run under Assumption~\ref{ass.tau_low}.  By Assumption~\ref{ass.tau_low}, it follows from the model reduction definition~\eqref{eq.model_reduction}, Lemma~\ref{lem.product_bounds} and \eqref{eq.linsys} that
  \bequationNN
    \baligned
      \E_k[\Delta l(x_k,\tau_k,g_k,d_k)]
      =&\ \E_k[-\tau_{\min} g_k^Td_k + \|c_k\|_2 - \|c_k + J_kd_k\|_2] \\
      &\leq -\tau_{\min} \nabla f(x_k)^T\dtrue_k + \tau_{\min} \zeta^{-1} M + \|c_k\|_2 - \|c_k + J_k\dtrue_k\|_2 \\
      &= \Delta l(x_k,\tau_{\min},\nabla f(x_k),\dtrue_k) + \tau_{\min} \zeta^{-1} M,
    \ealigned
  \eequationNN
  as desired. \Halmos
\endproof

We now prove our main theorem of this subsection, where $\E_{\tau,\low}[\ \cdot\ ] := \E[\ \cdot\ | \ \text{Assumption~\ref{ass.tau_low} holds}]$.

\btheorem\label{th.stochastic_tau_constant_small}
  Suppose that Assumption~\ref{ass.tau_low} holds and the sequence $\{\beta_k\}$ is chosen such that \eqref{eq.key_decrease_beta} holds for all $k \in \N{}$.  For a given run of the algorithm, define $\hat\xi_{\min} \in \R{}_{>0}$ as the value in Assumption~\ref{ass.tau_low} such that $\xi_k = \hat\xi_{\min}$ for all $k \geq k_\tau$ and define
  \bequationNN
    \baligned
      \underline{A} &:= \min\left\{\tfrac{2(1-\eta)\hat\xi_{\min} \tau_{\min}}{\tau_{\min} L + \Gamma}, \tfrac{2(1-\eta)\hat\xi_{\min}}{\tau_{\min} L + \Gamma} \right\}, \quad \Abar := \max\left\{\tfrac{2(1-\eta)\hat\xi_{\min} \tau_{\min}}{\tau_{\min} L + \Gamma}, \tfrac{2(1-\eta)\hat\xi_{\min}}{\tau_{\min} L + \Gamma} \right\}, \\
      \text{and}\ \ \Mbar &:= \tau_{\min} \zeta^{-1} ((1-\eta) (\Abar + \theta) M + \theta \kappa_{\nabla f} \sqrt{M}).
    \ealigned
  \eequationNN
  If $\beta_k = \beta \in (0,\underline{A}/((1-\eta)(\Abar + \theta)))$ for all $k \geq k_\tau$, then for all $k \geq k_\tau$ one finds
  \bequation\label{eq.beta_fixed}
    \baligned
      &\ \E_{\tau,\low}\left[\frac{1}{k-k_\tau+1} \sum_{j=k_\tau}^{k} \Delta l(x_j,\tau_{\min},\nabla f(x_j),\dtrue_j) \right] \\
      \leq&\ \frac{\beta \Mbar}{\underline{A} - (1-\eta)(\Abar + \theta)\beta} + \frac{\E_{\tau,\low}[\phi(x_{k_\tau}, \tau_{\min})] - \phi_{\min}}{(k+1) \beta(\underline{A} - (1-\eta)(\Abar + \theta)\beta)} \xrightarrow{k\to\infty} \frac{\beta \Mbar}{\underline{A} - (1-\eta)(\Abar + \theta)\beta},
    \ealigned
  \eequation
  where, in the context of Assumption~\ref{ass.main}, $\phi_{\min} \in \R{}_{>0}$ is a lower bound for $\phi(\cdot,\tau_{\min})$ over $\Xcal$.  On the other hand, if $\sum_{j=k_\tau}^\infty \beta_j = \infty$ and $\sum_{j=k_\tau}^\infty \beta_j^2 < \infty$, then
  \bequation\label{eq.beta_diminishing}
    \lim_{k \geq k_\tau, k \to \infty} \E_{\tau,\low}\left[ \tfrac{1}{\(\sum_{j=k_\tau}^{k} \beta_j\)} \sum_{j=k_\tau}^{k} \beta_j\Delta l(x_j,\tau_{\min},\nabla f(x_j),\dtrue_j) \right] = 0.
  \eequation
\etheorem
\proof{Proof.}
  Consider arbitrary $k \geq k_\tau$ in any run under Assumption~\ref{ass.tau_low}.  From the definitions of $\underline{A}$ and $\Abar$ in the statement of the theorem, the manner in which the step sizes are set by \eqref{eq.alpha}, and the fact that $\beta_k \in (0,1]$, it follows that $\underline{A} \beta_k \leq \alpha_k \leq (\Abar + \theta) \beta_k$.  Hence, it follows from Lemmas~\ref{lem.key_decrease}--\ref{lem.Deltaq} and the conditions of the theorem that
  \bequationNN
    \baligned
      &\ \phi(x_k, \tau_{\min}) - \E_k[\phi(x_k + \alpha_k d_k, \tau_{\min})] \\
      \geq&\ \E_k[\alpha_k \Delta l(x_k,\tau_{\min},\nabla f(x_k),\dtrue_k) - (1-\eta) \alpha_k \beta_k  \Delta l(x_k,\tau_{\min},g_k,d_k) - \alpha_k \tau_{\min} \nabla f(x_k)^T (d_k - \dtrue_k)] \\
      \geq&\ \beta_k(\underline{A} - (1-\eta)(\Abar + \theta) \beta_k) \Delta l(x_k,\tau_{\min},\nabla f(x_k),\dtrue_k) - \beta_k^2 \Mbar.
    \ealigned
  \eequationNN
  If $\beta_k = \beta \in (0,\underline{A}/((1-\eta)(\Abar + \theta)))$ for all $k \geq k_\tau$, then total expectation under Assumption~\ref{ass.tau_low} yields
  \begin{multline*}
    \E_{\tau,\low}[\phi(x_k, \tau_{\min})] - \E_{\tau,\low}[\phi(x_k + \alpha_k d_k, \tau_{\min})] \\
    \geq \beta(\underline{A} - (1-\eta)(\Abar + \theta)\beta) \E_{\tau,\low}[\Delta l(x_k,\tau_{\min},\nabla f(x_k),\dtrue_k)] - \beta^2 \Mbar\ \ \text{for all}\ \ k \geq k_\tau.
  \end{multline*}
  Summing this inequality for $j \in \{k_\tau,\dots,k\}$, it follows under Assumption~\ref{ass.main} that
  \bequationNN
    \baligned
      &\ \E_{\tau,\low}[\phi(x_{k_\tau},\tau_{\min})] - \phi_{\min} \\
      \geq&\ \E_{\tau,\low}[\phi(x_{k_\tau},\tau_{\min})] - \E_{\tau,\low}[\phi(x_{k+1},\tau_{\min})] \\
      \geq&\ \beta(\underline{A} - (1-\eta)(\Abar + \theta)\beta) \E_{\tau,\low}\left[\sum_{j=k_\tau}^{k} \Delta l(x_j,\tau_{\min},\nabla f(x_j),\dtrue_j) \right] - (k-k_\tau+1) \beta^2 \Mbar,
    \ealigned
  \eequationNN
  from which \eqref{eq.beta_fixed} follows.  On the other hand, if $\{\beta_k\}$ satisfies $\sum_{j=k_\tau}^\infty \beta_j = \infty$ and $\sum_{j=k_\tau}^\infty \beta_j^2 < \infty$, then it follows for sufficiently large $k \geq k_\tau$ that $\beta_k \leq \eta \underline{A} / ((1 - \eta)(\Abar + \theta))$; hence, without loss of generality, let us assume that this inequality holds for all $k \geq k_\tau$, which implies that $\underline{A} - (1-\eta)(\Abar + \theta)\beta_k \geq (1-\eta) \underline{A}$ for all $k \geq k_\tau$.  As above, it follows that
  \begin{multline*}
    \E_{\tau,\low}[\phi(x_k, \tau_{\min})] - \E_{\tau,\low}[\phi(x_k + \alpha_k d_k, \tau_{\min})] \\ \geq (1 - \eta) \underline{A} \beta_k \E_{\tau,\low}[\Delta l(x_k,\tau_{\min},\nabla f(x_k),\dtrue_k)] - \beta_k^2 \Mbar\ \ \text{for all}\ \ k \geq k_\tau.
  \end{multline*}
  Summing this inequality for $j \in \{k_\tau,\dots,k\}$, it follows under Assumption~\ref{ass.main} that
  \bequationNN
    \baligned
      \E_{\tau,\low}[\phi(x_{k_\tau},\tau_{\min})] - \phi_{\min} 
      \geq&\ \E_{\tau,\low}[\phi(x_{k_\tau},\tau_{\min})] - \E_{\tau,\low}[\phi(x_{k+1},\tau_{\min})] \\
      \geq&\ (1-\eta) \underline{A} \E_{\tau,\low}\left[ \sum_{j=k_\tau}^{k} \beta_j \Delta l(x_j,\tau_{\min},\nabla f(x_j),\dtrue_j)\right] - \Mbar \sum_{j=k_\tau}^{k} \beta_j^2.
    \ealigned
  \eequationNN
  Rearranging this inequality yields
  \bequationNN
    \E_{\tau,\low} \left[ \sum_{j=k_\tau}^{k} \beta_j \Delta l(x_j,\tau_{\min},g_j,\dtrue_j)\right] \leq \frac{\E_{\tau,\low}[\phi(x_{k_\tau},\tau_{\min})] - \phi_{\min}}{(1-\eta)\underline{A}} + \frac{\Mbar}{(1-\eta)\underline{A}} \sum_{j=k_\tau}^{k} \beta_j^2,
  \eequationNN
  from which \eqref{eq.beta_diminishing} follows.\Halmos
\endproof

We end this subsection with a corollary in which we connect the result of Theorem~\ref{th.stochastic_tau_constant_small} to first-order stationarity measures (recall \eqref{eq.KKT}).  For this corollary, we require the following lemma.

\blemma\label{lem.utrue}
  For all $k \in \N{}$, it holds that $\|\utrue_k\|_2 \geq \kappa_H^{-1} \|Z_k^T(\nabla f(x_k) + H_kv_k)\|_2$.
\elemma
\proof{Proof.}
  Consider arbitrary $k \in \N{}$ in any run.  As in the proof of Lemma~\ref{lem.expectation}, $Z_k^TH_kZ_k\wtrue_k = -Z_k^T(\nabla f(x_k) + H_kv_k)$, meaning with Assumption~\ref{ass.H} that $\|\utrue_k\|_2 \geq \kappa_H^{-1} \|Z_k^T(\nabla f(x_k) + H_kv_k)\|_2$. \Halmos
\endproof

\bcorollary\label{cor.stochastic_tau_finite_large_pi_finite}
  Under the conditions of Theorem~\ref{th.stochastic_tau_constant_small}, the following hold true.
  \benumerate
    \item[(a)] If $\beta_k = \beta \in (0,\underline{A}/((1-\eta)(\Abar + \theta)))$ for all $k \geq k_\tau$, then for all $k \geq k_\tau$ one finds
    \bequationNN
      \baligned
      &\ \E_{\tau,\low}\left[\frac{1}{k-k_\tau+1} \sum_{j=k_\tau}^{k} \(\frac{\tau_{\min}\zeta\|Z_j^T(\nabla f(x_j) + H_jv_{j})\|_2^2}{\kappa_H^2} + \frac{\kappa_v \sigma \|J_j^Tc_j\|_2^2}{\kappa_c}\) \right] \\
      \leq&\ \frac{\beta \Mbar}{\underline{A} - (1-\eta)(\Abar + \theta)\beta} + \frac{\E_{\tau,\low}[\phi(x_{k_\tau}, \tau_{\min})] - \phi_{\min}}{(k+1) \beta(\underline{A} - (1-\eta)(\Abar + \theta)\beta)} 
      \xrightarrow{k\to\infty} \frac{\beta \Mbar}{\underline{A} - (1-\eta)(\Abar + \theta)\beta}.
      \ealigned
    \eequationNN
    \item[(b)] If $\sum_{j=k_\tau}^\infty \beta_j = \infty$ and $\sum_{j=k_\tau}^\infty \beta_j^2 < \infty$, then
    \bequationNN
      \lim_{k \geq k_\tau, k\to\infty} \E_{\tau,\low}\left[ \frac{1}{\(\sum_{j=k_\tau}^{k} \beta_j\)} \sum_{j=k_\tau}^{k} \beta_j \(\frac{\tau_{\min}\zeta \|Z_j^T(\nabla f(x_j) + H_jv_j)\|_2^2}{\kappa_H^2} + \frac{\kappa_v \sigma\|J_j^Tc_j\|_2^2}{\kappa_c}\) \right] = 0,
    \eequationNN
    from which it follows that
    \bequationNN
      \liminf_{k \geq k_\tau, k\to\infty}\ \E_{\tau,\low} \left[ \frac{\tau_{\min}\zeta\|Z_k^T(\nabla f(x_k) + H_kv_k)\|_2^2}{\kappa_H^2} + \frac{\kappa_v \sigma\|J_k^Tc_k\|_2^2}{\kappa_c} \right] = 0.
    \eequationNN
  \eenumerate
\ecorollary
\proof{Proof.}
  For all $k \in \N{}$, it follows under Assumption~\ref{ass.tau_low} that \eqref{eq.model_reduction_condition} holds with $\nabla f(x_k)$ in place of $g_k$ and $\utrue_k$ in place of $u_k$.  The result follows from this fact, Theorem~\ref{th.stochastic_tau_constant_small}, and Lemmas~\ref{lem.normal_component} and~\ref{lem.utrue}. \Halmos
\endproof

Observe that if the singular values of $J_k$ are bounded below by $\sigma_J \in \R{}_{>0}$ for all $k \geq k_J$ for some $k_J \in \N{}$, then (as in the proof of Lemma~\ref{lem.det.v_upper}) it follows that $\|J_k^Tc_k\|_2 \geq \sigma_J \|c_k\|_2$ for all $k \geq k_J$.  In this case, the results of Corollary~\ref{cor.stochastic_tau_finite_large_pi_finite} hold with $\sigma_J\|c_k\|_2$ in place of $\|J_k^Tc_k\|_2$.  Overall, Corollary~\ref{cor.stochastic_tau_finite_large_pi_finite} offers results for the stochastic setting that parallel the limits \eqref{eq.det.limit1} and \eqref{eq.det.limit3} for the deterministic setting.  The only difference is the presence of $Z_k^TH_kv_k$ in the term involving the reduced gradient $Z_k^T\nabla f(x_k)$ for all $k \in \N{}$.  However, this does not significantly weaken the conclusion.  After all, it follows from~\eqref{prob.normal} (see also Lemma~\ref{lem.normal_component}) that $\|v_k\|_2 \leq \omega \|J_k^Tc_k\|_2$ for all $k \in \N{}$.  Hence, since Corollary~\ref{cor.stochastic_tau_finite_large_pi_finite} shows that at least a subsequence of $\{\|J_k^Tc_k\|_2\}$ tends to vanish in expectation, it follows that $\{\|v_k\|_2\}$ vanishes in expectation along the same subsequence of iterations.  This, along with Assumption~\ref{ass.H} and the orthonormality of $Z_k$, shows that $\{\|Z_k^TH_kv_k\|_2\}$ exhibits this same behavior, which means that from the corollary one finds that a subsequence of $\{\|Z_k^T\nabla f(x_k)\|_2\}$ vanishes in expectation.

Let us conclude this subsection with a few remarks on how one should interpret its main conclusions.  First, one learns from the requirements on $\{\beta_k\}$ in Lemma~\ref{lem.key_decrease},  Theorem~\ref{th.stochastic_tau_constant_small}, and Corollary~\ref{cor.stochastic_tau_finite_large_pi_finite} that, rather than employ a prescribed sequence $\{\beta_k\}$, one should instead prescribe $\{\hat\beta_j\}_{j=0}^\infty \subset (0,1]$ and for each $k \in \N{}$ set $\beta_k$ based on whether or not an adaptive parameter changes its value.  In particular, anytime $k \in \N{}$ sees either $\tau_k < \tau_{k-1}$, $\chi_k > \chi_{k-1}$, $\zeta_k < \zeta_{k-1}$, or $\xi_k < \xi_{k-1}$, the algorithm should set $\beta_{k+j} \gets \lambda \hat\beta_j$ for $j = 0,1,2,\dots$ (continuing indefinitely or until $\khat \in \N{}$ with $\khat > k$ sees $\tau_{\khat} < \tau_{\khat-1}$, $\chi_{\khat} > \chi_{\khat-1}$, $\zeta_{\khat} < \zeta_{\khat-1}$, or $\xi_{\khat} < \xi_{\khat-1}$), where $\lambda \in \R{}_{>0}$ is chosen sufficiently small such that~\eqref{eq.key_decrease_beta} holds.  Since such a ``reset'' of $j \gets 0$ will occur only a finite number of times under event $E_{\tau,\low}$, one of the desirable results in Theorem~\ref{th.stochastic_tau_constant_small}/Corollary~\ref{cor.stochastic_tau_finite_large_pi_finite} can be attained if $\{\hat\beta_j\}$ is chosen as an appropriate constant or diminishing sequence.  Second, let us note that due to the generality of Assumption~\ref{ass.tau_low}, it is possible that for different runs of the algorithm the corresponding terminal merit parameter value, namely, $\tau_{\min}$, in Assumption~\ref{ass.tau_low} could become arbitrarily close to zero.  (This is in contrast to the conditions of Lemma~\ref{lem.det.pi_bounded}, which guarantee a \emph{uniform} lower bound for the merit parameter over all runs satisfying these conditions.)  Hence, while our main conclusions of this subsection hold under looser conditions than those in Lemma~\ref{lem.det.pi_bounded}, one should be wary in practice if/when the merit parameter sequence reaches small numerical values.

%************
% Subsection
%************
\subsection{Vanishing Merit Parameter}\label{sec.tau_zero}

Let us now consider the behavior of the algorithm in settings in which the merit parameter vanishes; in particular, we make Assumption~\ref{ass.tau_zero} below.

\bassumption\label{ass.tau_zero}
  In a run of the algorithm, event $E_{\tau,\zero}$ occurs, i.e., $\{\tau_k\} \searrow 0$.  In addition, along the lines of Assumption~\ref{ass.g}, the stochastic gradient sequence $\{g_k\}$ satisfies $\E_{k,\tau,\zero}[g_k] = \nabla f(x_k)$ and $\|g_k - \nabla f(x_k)\|_2^2 \leq M$, where $\E_{k,\tau,\zero}$ denotes expectation with respect to the distribution of $\iota$ conditioned on the event that $E_{\tau,\zero}$ occurs and the algorithm has reached $x_k$ in iteration $k \in \N{}$.
\eassumption

Recalling Theorem~\ref{th.deterministic} and Lemma~\ref{lem.det.pi_bounded}, one may conclude in general that the merit parameter sequence may vanish for one of two reasons: a (sub)sequence of constraint Jacobians tends toward rank deficiency or a (sub)sequence of stochastic gradient estimates diverges.  Our assumption here assumes that the latter event does not occur.  (In our remarks in Section~\ref{sec.diverge}, we discuss the obstacles that arise in proving convergence guarantees when the merit parameter vanishes and the stochastic gradient estimates diverge.)  Given our setting of constrained optimization, it is reasonable and consistent with Theorem~\ref{th.deterministic} to have convergence toward stationarity with respect to the constraint violation measure as the primary goal in these circumstances.

For the sake of brevity in the rest of this subsection, let us temporarily redefine $\E_k := \E_{k,\tau,\zero}$.

Our first result in this subsection is an alternative of Lemma~\ref{lem.key_decrease}.

\blemma\label{lem.key_decrease_zero}
  Under Assumption~\ref{ass.tau_zero} and assuming that $\{\beta_k\}$ is chosen such that \eqref{eq.key_decrease_beta} holds for all $k \in \N{}$, it follows for all $k \in \N{}$ that
  \bequationNN
    \baligned
      &\ \|c_k\|_2 - \|c(x_k + \alpha_k d_k)\|_2 \\
      \geq&\ \alpha_k (1 - (1 - \eta)\beta_k) \Delta l(x_k,\tau_k,g_k,d_k) - \tau_k (f_k - f(x_k + \alpha_k d_k)) - \alpha_k \tau_k (\nabla f(x_k) - g_k)^Td_k.
    \ealigned
  \eequationNN
\elemma
\proof{Proof.}
  Consider arbitrary $k \in \N{}$ in any run under Assumption~\ref{ass.tau_zero}.  As in the proof of Lemma~\ref{lem.key_decrease}, from \eqref{eq.alpha_suff}--\eqref{eq.alpha} and the supposition about~$\{\beta_k\}$, one finds $\alpha_k \in (0,1]$.  Hence, with \eqref{eq.merit_decrease}, one has
  \bequationNN
    \baligned
      \phi(x_k, \tau_k) - \phi(x_k + \alpha_k d_k, \tau_k)
      \geq&\ -\alpha_k (\tau_k \nabla f(x_k)^T d_k - \|c_k\|_2 + \|c_k + J_kd_k\|_2) - \thalf (\tau_k L + \Gamma) \alpha_k^2 \|d_k\|_2^2 \\
      =&\ \alpha_k \Delta l(x_k,\tau_k,g_k,d_k) - \thalf (\tau_k L + \Gamma) \alpha_k^2 \|d_k\|_2^2 - \alpha_k\tau_k(\nabla f(x_k) - g_k)^Td_k.
    \ealigned
  \eequationNN
  Now following the same arguments as in the proof of Lemma~\ref{lem.key_decrease}, it follows that $- \thalf (\tau_k L + \Gamma) \alpha_k^2 \|d_k\|_2^2 \geq - (1 - \eta) \alpha_k \beta_k \Delta l(x_k,\tau_k,g_k,d_k)$, which combined with the above yields the desired conclusion. \Halmos
\endproof

Our next result yields a bound on the final term in the conclusion of Lemma~\ref{lem.key_decrease_zero}.

\blemma\label{lem.alphagd_zero}
  For any run under Assumption~\ref{ass.tau_zero}, there exists $\kappa_\beta \in \R{}_{>0}$ such that
  \bequationNN
    \alpha_k \tau_k (\nabla f(x_k) - g_k)^T d_k \leq \bcases \beta_k \tau_k \kappa_\beta & \text{for all $k \in \N{}$ such that $\|u_k\|_2^2 < \chi_k \|v_k\|_2^2$} \\ \beta_k \tau_k \max\{\beta_k,\tau_k\} \kappa_\beta & \text{for all $k \in \N{}$ such that $\|u_k\|_2^2 \geq \chi_k \|v_k\|_2^2$.} \ecases
  \eequationNN
\elemma
\proof{Proof.}
  The existence of $\kappa_d \in \R{}_{>0}$ such that, in any run under Assumption~\ref{ass.tau_zero}, one finds $\|d_k\|_2 \leq \kappa_d$ for all $k \in \N{}$ follows from Assumption~\ref{ass.tau_zero}, the fact that $\|d_k\|_2^2 = \|v_k\|_2^2 + \|u_k\|_2^2$ for all $k \in \N{}$, Lemma~\ref{lem.det.u_bound}, Lemma~\ref{lem.normal_component}, and Assumption~\ref{ass.main}.  Now consider arbitrary $k \in \N{}$ in any run under Assumption~\ref{ass.tau_zero}.  If $(\nabla f(x_k) - g_k)^Td_k < 0$, then the desired conclusion follows trivially (for any $\kappa_\beta \in \R{}_{>0}$).  Hence, let us proceed under the assumption that $(\nabla f(x_k) - g_k)^Td_k \geq 0$.  If $\|u_k\| < \chi_k \|v_k\|_2^2$, then it follows from~\eqref{eq.alpha}, the facts that $0 \leq \tau_k$, $\xi_k \leq \xi_{-1}$, and $\beta_k \leq 1$ for all $k \in \N{}$, the Cauchy-Schwarz inequality, and Assumption~\ref{ass.tau_zero} that
  \bequationNN
    \baligned
      \alpha_k \tau_k (\nabla f(x_k) - g_k)^Td
      &\leq \(\frac{2(1-\eta)\beta_k\xi_k}{\tau_k L + \Gamma} + \theta \beta_k^2\) \tau_k \|\nabla f(x_k) - g_k\|_2 \|d_k\|_2 \\
      &\leq \(\frac{2(1-\eta)\xi_{-1}}{\Gamma} + \theta\) \beta_k \tau_k \sqrt{M} \kappa_d.
    \ealigned
  \eequationNN
  On the other hand, if $\|u_k\|_2^2 \geq \chi_k \|v_k\|_2^2$, then it follows under the same reasoning that
  \bequationNN
    \baligned
      \alpha_k \tau_k (\nabla f(x_k) - g_k)^Td
      &\leq \(\frac{2(1-\eta)\beta_k\xi_k\tau_k}{\tau_k L + \Gamma} + \theta \beta_k^2\) \tau_k \|\nabla f(x_k) - g_k\|_2 \|d_k\|_2 \\
      &\leq \(\frac{2(1-\eta)\xi_{-1}}{\Gamma} + \theta\) \beta_k \tau_k \max\{\beta_k,\tau_k\} \sqrt{M} \kappa_d.
    \ealigned
  \eequationNN
  Overall, the desired conclusion follows with $\kappa_{\beta} := \left(\tfrac{2 (1-\eta) \xi_{-1}}{\Gamma} + \theta \right) \sqrt{M} \kappa_d$. \Halmos
\endproof

Our third result in this subsection offers a formula for a positive lower bound on the step size that is applicable at points that are not stationary for the constraint infeasibility measure.  For this lemma and its subsequent consequences, we define for arbitrary $\gamma \in \R{}_{>0}$ the subset
\bequation\label{eq.X_gamma}
  \Xcal_\gamma := \{x \in \R{n} : \|J(x)^Tc(x)\|_2 \geq \gamma\}.
\eequation

\blemma\label{lem.tau_zero_alpha_low}
  There exists $\underline\alpha \in \R{}_{>0}$ such that $\alpha_k \geq \underline\alpha\beta_k$ for each $k \in \N{}$ such that $\|u_k\|_2^2 < \chi_k \|v_k\|_2^2$.  On the other hand, for each $\gamma \in \R{}_{>0}$, there exists $\epsilon_\gamma \in \R{}_{>0}$ $($proportional to $\gamma^2$$)$ such that
  \bequationNN
    x_k \in \Xcal_\gamma\ \ \text{implies}\ \ \alpha_k \geq \min\{\epsilon_\gamma \beta_k, \epsilon_\gamma \beta_k \tau_k + \theta \beta_k^2\}\ \ \text{whenever}\ \ \|u_k\|_2^2 \geq \chi_k \|v_k\|_2^2.
  \eequationNN
\elemma
\proof{Proof.}
  Define $\Kcal_\gamma := \{k \in \N{} : x_k \in \Xcal_\gamma\}$.  By Lemma~\ref{lem.normal_component}, it follows that $\|v_k\|_2 \geq \underline\omega \|J_k^Tc_k\|_2^2 \geq \underline\omega \gamma^2$ for all $k \in \Kcal_\gamma$.  Consequently, by Lemma~\ref{lem.det.u_bound}, it follows that
  \bequation\label{eq.fec}
    \|u_k\|_2 \leq \frac{\kappa_u}{\underline\omega \gamma^2} \|v_k\|_2\ \ \text{for all}\ \ k \in \Kcal_\gamma.
  \eequation
  It follows from \eqref{eq.alpha} that $\alpha_k \geq 2 (1 - \eta) \beta_k \xi_k / (\tau_k L + \Gamma)$ whenever $\|u_k\|_2^2 < \chi_k \|v_k\|_2^2$.  Otherwise, whenever $\|u_k\|_2^2 \geq \chi_k \|v_k\|_2^2$, it follows using the arguments in Lemma~\ref{lem.key_decrease} and \eqref{eq.alpha} that
  \bequationNN
    \alpha_k = \min\left\{ \frac{2 (1 - \eta) \beta_k \Delta l(x_k,\tau_k,g_k,d_k)}{(\tau_k L + \Gamma) \|d_k\|_2^2}, \frac{2 (1 - \eta) \beta_k \xi_k \tau_k}{\tau_k L + \Gamma} + \theta \beta_k^2, 1 \right\},
  \eequationNN
  which along with \eqref{eq.model_reduction_condition}, Lemma~\ref{lem.normal_component}, \eqref{eq.bounds}, and \eqref{eq.fec} imply that
  \bequationNN
    \baligned
      \alpha_k
        &\geq \min\left\{ \frac{ 2 ( 1- \eta) \beta_k \sigma (\|c_k\|_2 - \|c_k + J_kv_k\|_2)}{(\tau_k L + \Gamma) (\|u_k\|_2^2 + \|v_k\|_2^2)}, \frac{2 (1 - \eta) \beta_k \xi_k \tau_k}{\tau_k L + \Gamma} + \theta \beta_k^2, 1 \right\} \\
        &\geq \min\left\{ \frac{ 2 ( 1- \eta) \beta_k \sigma \kappa_v\|J_k^T c_k\|^2}{(\tau_k L + \Gamma) (\frac{\kappa_u^2}{\underline{\omega}^2\gamma^4}+1)\omega^2\|c_k\|\|J_k^T c_k\|^2}, \frac{2 (1 - \eta) \beta_k \xi_k \tau_k}{\tau_k L + \Gamma} + \theta \beta_k^2, 1 \right\} \\
        &\geq \min\left\{ \frac{2 (1 - \eta)  \beta_k \sigma \kappa_v \underline\omega^2 \gamma^4}{(\tau_k L + \Gamma) \kappa_c \omega^2 (\kappa_u^2 + \underline\omega^2 \gamma^4)}, \frac{2 (1 - \eta) \beta_k \xi_k \tau_k}{\tau_k L + \Gamma} + \theta \beta_k^2, 1 \right\}.
    \ealigned
  \eequationNN
  Combining the cases above with Lemma~\ref{lem.xi_lower} yields the desired conclusion. \Halmos
\endproof

We now prove our main theorem of this subsection.

\btheorem\label{th.tau_zero}
  Suppose that Assumption~\ref{ass.tau_zero} holds, the sequence $\{\beta_k\}$ is chosen such that \eqref{eq.key_decrease_beta} holds for all $k \in \N{}$, and either
  \benumerate
    \item[(a)] $\beta_k = \beta \in (0,1)$ for all $k \in \N{}$, or
    \item[(b)] $\sum_{k=0}^\infty \beta_k = \infty$, $\sum_{k=0}^\infty \beta_k^2 < \infty$, and either $|\{k \in \N{} : \|u_k\|_2^2 < \chi_k \|v_k\|_2^2\}| = \infty$ or $\sum_{k=0}^\infty \beta_k \tau_k = \infty$.
  \eenumerate
  Then, $\liminf_{k\to\infty} \|J_k^Tc_k\|_2 = 0$.
\etheorem
\proof{Proof.}  
  To derive a contradiction, suppose that there exists $k_\gamma \in \N{}$ and $\gamma \in \R{}_{>0}$ such that $x_k \in \Xcal_\gamma$ for all $k \geq k_\gamma$.  Our aim is to show that, under condition (a) or (b), a contradiction is reached.
  
  First, suppose that condition (a) holds.  By Lemmas~\ref{lem.key_decrease_zero}--\ref{lem.tau_zero_alpha_low}, \eqref{eq.bounds}, \eqref{eq.model_reduction_condition}, the fact that $\beta \in (0,1)$, Lemma~\ref{lem.normal_component}, and Assumption~\ref{ass.main}, it follows that there exists $\underline{\epsilon}_\gamma \in \R{}_{>0}$ such that
  \begin{align}
    \|c_k\|_2 - \|c_{k+1}\|_2
      \geq&\ \alpha_k (1 - (1 - \eta)\beta) \Delta l(x_k,\tau_k,g_k,d_k) - \tau_k (f_k - f_{k+1}) - \alpha_k \tau_k (\nabla f(x_k) - g_k)^Td_k \nonumber \\
      \geq&\ \underline{\epsilon}_\gamma \beta \eta \sigma (\|c_k\|_2 - \|c_k + J_kv_k\|_2) - \tau_k (f_{\sup} - f_{\inf}) - \beta \tau_k \max\{1,\tau_k\} \kappa_\beta \nonumber \\
      \geq&\ \underline{\epsilon}_\gamma \beta \eta \sigma \kappa_v \kappa_c^{-1} \|J_k^Tc_k\|_2^2 - \tau_k (f_{\sup} - f_{\inf} + \beta \max\{1,\tau_k\} \kappa_\beta)\ \ \text{for all}\ \ k \geq k_\gamma. \label{eq.ugh}
  \end{align}
  Since $\|J_k^Tc_k\|_2 \geq \gamma$ for all $k \geq k_\gamma$ and $\{\tau_k\} \searrow 0$ under Assumption~\ref{ass.tau_zero}, it follows that there exists $k_\tau \geq k_\gamma$ such that $\tau_k (f_{\sup} - f_{\inf} + \beta \max\{1,\tau_k\} \kappa_\beta) \leq \thalf \underline{\epsilon}_\gamma \beta \eta \sigma \kappa_v \kappa_c^{-1} \|J_k^Tc_k\|_2^2$ for all $k \geq k_\tau$.  Hence, summing~\eqref{eq.ugh} for $j \in \{k_\tau,\dots,k\}$, it follows with \eqref{eq.bounds} that
  \bequationNN
    \baligned
      \kappa_c \geq \|c_{k_\tau}\|_2 - \|c_{k+1}\|_2 \geq \half \underline{\epsilon}_\gamma \beta \eta \sigma \kappa_v \kappa_c^{-1} \sum_{j=k_\tau}^{k} \|J_j^Tc_j\|_2^{2}.
    \ealigned
  \eequationNN
  It follows from this fact that $\{J_k^Tc_k\}_{k\geq k_\tau} \to 0$, yielding the desired contradiction.
  
  Second, suppose that condition (b) holds.  Since $\sum_{k=0}^\infty \beta_k^2 < \infty$, it follows that there exists $k_\beta \in \N{}$ with $k_\beta \geq k_\gamma$ such that $(1 - (1-\eta)\beta_k) \geq \eta$ for all $k \geq k_\beta$.  Hence, for all $k \geq k_\beta$ with $\|u_k\|_2^2 < \chi_k \|v_k\|_2^2$, one finds from Lemmas~\ref{lem.key_decrease_zero}--\ref{lem.tau_zero_alpha_low}, \eqref{eq.bounds}, \eqref{eq.model_reduction_condition}, Lemma~\ref{lem.normal_component}, and Assumption~\ref{ass.main} that
  \begin{align*}
    %&\ \|c_k\|_2 - \|c_{k+1}\|_2 \\
      \|c_k\|_2 - \|c_{k+1}\|_2 \geq&\ \alpha_k (1 - (1 - \eta)\beta_k) \Delta l(x_k,\tau_k,g_k,d_k) - \tau_k (f_k - f_{k+1}) - \alpha_k \tau_k (\nabla f(x_k) - g_k)^Td_k \\
      \geq&\ \beta_k \underline\alpha \eta \sigma \kappa_v \kappa_c^{-1} \|J_k^Tc_k\|_2^2 - \tau_k (f_k - f_{\inf}) + \tau_k(f_{k+1} - f_{\inf}) - \beta_k \tau_k \kappa_\beta \\
      \geq&\ \beta_k \underline\alpha \eta \sigma \kappa_v \kappa_c^{-1} \|J_k^Tc_k\|_2^2 - \tau_{k-1} (f_k - f_{\inf}) + \tau_k(f_{k+1} - f_{\inf}) - \beta_k \tau_k \kappa_\beta.
  \end{align*}
  Similarly, for all sufficiently large $k \geq k_\beta$---specifically, $k \geq \kbar_\beta$, where $\kbar_\beta \in \N{}$ is sufficiently large such that $\kbar_\beta \geq k_\beta$ and $\epsilon_\gamma \geq \epsilon_\gamma \tau_k + \theta \beta_k$---with $\|u_k\|_2^2 \geq \chi_k \|v_k\|_2^2$, similar reasoning yields
  \begin{align*}
    &\ \|c_k\|_2 - \|c_{k+1}\|_2 \\
      \geq&\ \alpha_k (1 - (1 - \eta)\beta_k) \Delta l(x_k,\tau_k,g_k,d_k) - \tau_k (f_k - f_{k+1}) - \alpha_k \tau_k (\nabla f(x_k) - g_k)^Td_k \\
      \geq&\ \beta_k \max\{\beta_k,\tau_k\} \min\{\epsilon_\gamma,\theta\} \eta \sigma \kappa_v \kappa_c^{-1} \|J_k^Tc_k\|_2^2 - \tau_{k-1} (f_k - f_{\inf}) + \tau_k(f_{k+1} - f_{\inf}) - \beta_k \max\{\beta_k, \tau_k\} \tau_k \kappa_\beta.
  \end{align*}
  Since $\|J_k^Tc_k\|_2 \geq \gamma$ for all $k \geq \kbar_\beta \geq k_\beta \geq k_\gamma$ and $\{\tau_k\} \searrow 0$ under Assumption~\ref{ass.tau_zero}, it follows that there exists $k_\tau \geq \kbar_\beta$ such that $\tau_k \kappa_\beta \leq \thalf \underline\alpha \eta \sigma \kappa_v \kappa_c^{-1} \|J_k^Tc_k\|_2^2$ and $\tau_k \kappa_\beta \leq \thalf \min\{\epsilon_\gamma,\theta\} \eta \sigma \kappa_v \kappa_c^{-1} \|J_k^Tc_k\|_2^2$ for all $k \geq k_\tau$.  Hence, letting $\Kcal_u := \{k \in \N{} : \|u_k\|_2^2 \geq \chi_k \|v_k\|_2^2\}$ and $\Kcal_v:= \{k \in \N{} : \|u_k\|_2^2 < \chi_k \|v_k\|_2^2\}$, one can sum the inequalities above for $j \in \{k_\tau,...,k\}$ to find
  \begin{align}
    \kappa_c \geq \|c_{k_\tau}\|_2 - \|c_{k+1}\|_2
      \geq&\ - \tau_{k_\tau - 1} (f_{k_\tau} - f_{\inf}) + \tau_k (f_{k+1} - f_{\inf}) \nonumber \\
      &\ + \sum_{j=k_\tau,j\in\Kcal_v}^k \beta_j \(\underline\alpha \eta \sigma \kappa_v \kappa_c^{-1} \|J_j^Tc_j\|_2^2 - \tau_j \kappa_\beta\) \nonumber \\
      &\ + \sum_{j=k_\tau,j\in\Kcal_u}^k \beta_j \max\{\beta_j,\tau_j\} \(\min\{\epsilon_\gamma,\theta\} \eta \sigma \kappa_v \kappa_c^{-1} \|J_j^Tc_j\|_2^2 - \tau_k \kappa_\beta\) \nonumber \\
        \geq&\ - \tau_{k_\tau - 1} (f_{k_\tau} - f_{\inf}) \nonumber \\
      &\ + \sum_{j=k_\tau,j\in\Kcal_v}^k \beta_j \thalf \underline\alpha \eta \sigma \kappa_v \kappa_c^{-1} \|J_j^Tc_j\|_2^2 \nonumber \\
      &\ + \sum_{j=k_\tau,j\in\Kcal_u}^k \beta_j \max\{\beta_j,\tau_j\} \thalf \min\{\epsilon_\gamma,\theta\} \eta \sigma \kappa_v \kappa_c^{-1} \|J_j^Tc_j\|_2^2. \label{eq.toref}
  \end{align}
  It follows from this fact and the fact that either $|\Kcal_v| = \infty$ or at least $\sum_{j=k_\tau,j\in\Kcal_u} \beta_j \tau_j = \infty$ that $\{J_k^Tc_k\}_{k\geq k_\tau} \to 0$, yielding the desired contradiction. \Halmos
\endproof

There is one unfortunate case not covered by Theorem~\ref{th.tau_zero}, namely, the case when $\{\beta_k\}$ diminishes (as in condition (b)), the search direction is tangentially dominated for all sufficiently large $k \in \N{}$, and $\sum_{k=0}^\infty \beta_k \tau_k < \infty$.  One can see in the proof of the theorem why the desired conclusion, namely, that the limit inferior of $\{\|J_k^Tc_k\|_2\}$ is zero, does not necessarily follow in this setting: If, after some iteration, all search directions are tangentially dominated and $\sum_{k=0}^\infty \beta_k \tau_k < \infty$, then the coefficients on $\|J_k^Tc_k\|_2^2$ in \eqref{eq.toref} are summable, which means that there might not be a subsequence of $\{\|J_k^Tc_k\|_2^2\}$ that vanishes.  Fortunately, however, this situation is detectable in practice, in the sense that one can detect it using computed quantities.  In particular, if $\beta_k$ is below a small threshold, $\|J_k^Tc_k\|_2$ has remained above a threshold in all recent iterations, $\tau_k = \Ocal(\beta_k)$ in recent iterations, and the algorithm has computed tangentially dominated search directions in all recent iterations, then the algorithm may benefit by triggering a switch to a setting in which $\{\beta_k\}$ is kept constant in future iterations, in which case the desired conclusion follows under condition (a).  Such a trigger arguably does not conflict much with Section~\ref{sec.tau_low}, since the analysis in that section presumes that $\{\tau_k\}$ remains bounded away from zero, whereas here one has confirmed that $\tau_k \approx 0$.

%************
% Subsection
%************
\subsection{Constant, Insufficiently Small Merit Parameter}\label{sec.tau_big}

Our goal now is to consider the event that the algorithm generates a merit parameter sequence that eventually remains constant, but at a value that is too large in the sense that the conditions of Assumption~\ref{ass.tau_low} do not hold.  Such an event for the algorithm in~\cite{BeraCurtRobiZhou21} is addressed in Proposition~3.16 in that article, where under a reasonable assumption (paralleling \eqref{eq.daniel_special_label}, which we discuss later on) 
it is shown that, in a given run of the algorithm, the probability is zero of the merit parameter settling on too large of a value.  The same can be said of our algorithm, as discussed in this subsection.  That said, this does not address what might be the total probability, over all runs of the algorithm, of the event that the merit parameter remains too large.  We discuss in this section that, under reasonable assumptions, this total probability is zero, where a formal theorem and proof are provided in Appendix~\ref{sec.totalprob}.

For our purposes in this section, we make some mild simplifications.  First, as shown in Lemmas~\ref{lem.chi_upper} and~\ref{lem.xi_lower}, each of the sequences $\{\chi_k\}$, $\{\zeta_k\}$, and $\{\xi_k\}$ has a uniform bound that holds over any run of the algorithm.  Hence, for simplicity, we shall assume that the initial values of these sequences are chosen such that they are constant over $k \in \N{}$.  (Our discussions in this subsection can be generalized to situations when this is not the case; the conversation merely becomes more cumbersome, which we have chosen to avoid.)  Second, it follows from properties of the deterministic instance of our algorithm (recall Theorem~\ref{th.deterministic}) that if a subsequence of $\{\tautruetrial_k\}$ converges to zero, then a subsequence of the sequence of minimum singular values of the constraint Jacobians $\{J_k\}$ vanishes as well.  Hence, we shall consider in this subsection events in which there exists $\tautruetrial_{\min} \in \R{}_{>0}$ such that $\tautruetrial_k \geq \tautruetrial_{\min}$ for all $k \in \N{}$ in any run of the algorithm.  (We will remark on the consequences of this assumption further in Section~\ref{sec.diverge}.)  It follows from this and~\eqref{eq.tau_update} that if the cardinality of the set of iteration indices $\{k \in \N{} : \tau_k < \tau_{k-1}\}$ ever exceeds
\bequation\label{eq.sbar}
  \sbar(\tautruetrial_{\min}) := \left\lceil \frac{\log(\tautruetrial_{\min}/\tau_{-1})}{\log(1-\epsilon_\tau)} \right\rceil \in \N{},
\eequation
then for all subsequent $k \in \N{}$ one has $\tau_{k-1} \leq \tautruetrial_{\min} \leq \tautruetrial_k$.  This property of $\sbar(\tautruetrial_{\min})$ is relevant in our event of interest for this subsection, which we now define.

\bdefinition\label{def.tau_big}
  \emph{
  The event $E_{\tau,\big}(\tautruetrial_{\min})$ for some $\tautruetrial_{\min} \in \R{}_{>0}$ occurs in a run if and only if $\tautruetrial_k \geq \tautruetrial_{\min}$ for all $k \in \N{}$ and there exists an infinite index set $\Kcal \subseteq \N{}$ such that
  \bequation\label{eq.funky}
    \tautruetrial_k < \tau_{k-1}\ \ \text{for all}\ \ k \in \Kcal.
  \eequation
  }
\edefinition

Considering a given run of our algorithm in which it is presumed that $\tautruetrial_k \geq \tautruetrial_{\min}$ for some $\tautruetrial_{\min} \in \R{}_{>0}$ for all $k \in \N{}$, one has under a reasonable assumption (specifically, that \eqref{eq.daniel_special_label} in the lemma below holds for all $k \in \N{}$) that the probability is zero that $E_{\tau,\big}(\tautruetrial_{\min})$ occurs.  We prove this now using the same argument as in the proof of \cite[Proposition~3.16]{BeraCurtRobiZhou21}.  For this, we require the following lemma, proved here for our setting, which is slightly different than for the algorithm in \cite{BeraCurtRobiZhou21} (due to the slightly different formula for setting the merit parameter).

\blemma\label{lem:tautrialextra}
  For any $k \in \N{}$ in any run of the algorithm, it follows for any $p \in (0,1]$ that
  \bsubequations
    \begin{align}
      \P_k[g_k^Td_k + u_k^TH_ku_k \geq \nabla f(x_k)^T \dtrue_k + (\utrue_k)^TH_k\utrue_k] &\geq p \label{eq.daniel_special_label} \\
      \text{implies}\ \ \P_k[\tau_k < \tau_{k-1} | \tautruetrial_k < \tau_{k-1}] &\geq p. \label{eq.frank_special_label}
    \end{align}
  \esubequations
\elemma
\proof{Proof.}
  Consider any $k \in \N{}$ in any run of the algorithm such that $\tautruetrial_k < \tau_{k-1} \in \R{}_{>0}$. Then, it follows from \eqref{eq.merit_parameter_trial} that $\tautruetrial_k < \infty$, $\nabla f(x_k)^T \dtrue_k + (\utrue_k)^T H_k \utrue_k > 0$, and
  \bequationNN
    \tautruetrial_k = \frac{(1-\sigma)(\|c_k\|_2 - \|c_k + J_k\dtrue_k\|_2)}{\nabla f(x_k)^T \dtrue_k + (\utrue_k)^T H_k \utrue_k} < \tau_{k-1},
  \eequationNN
  from which it follows that
  \bequation\label{eq.mike1}
    (1-\sigma)(\|c_k\|_2 - \|c_k + J_k\dtrue_k\|_2) < (\nabla f(x_k)^T \dtrue_k + (\utrue_k)^T H_k \utrue_k) \tau_{k-1}.
  \eequation
  If, in addition, a realization of $g_k$ yields
  \bequation\label{eq.mike2}
    g_k^T d_k + u_k^TH_ku_k \geq \nabla f(x_k)^T \dtrue_k + (\utrue_k)^T H_k \utrue_k,
  \eequation
  then it follows from \eqref{eq.mike1} and the fact that $J_k\dtrue_k = J_kd_k$ that
  \bequationNN
    (1-\sigma)(\|c_k\|_2 - \|c_k + J_kd_k\|_2) < (g_k^T d_k + u_k^T H_k u_k) \tau_{k-1}.
  \eequationNN
  It follows from this inequality and Lemma~\ref{lem.normal_component} that $g_k^T d_k + u_k^T H_k u_k > 0$, and with \eqref{eq.tau_update} it holds that
  \bequationNN
    \tau_k \leq \tautrial_k = \frac{(1-\sigma)(\|c_k\|_2 - \|c_k + J_kd_k\|_2)}{g_k^T d_k + u_k^T H_k u_k} < \tau_{k-1}.
  \eequationNN
  Hence, conditioned on the event that $\tautruetrial_k < \tau_{k-1}$, one finds that \eqref{eq.mike2} implies that $\tau_k < \tau_{k-1}$.  Therefore, under the conditions of the lemma and the fact that, conditioned on the events leading up to iteration number $k$ one has that both $\tautruetrial_k$ and $\tau_{k-1}$ are deterministic, it follows that
  \begin{align*}
        &\ \P_k[\tau_k < \tau_{k-1} | \tautruetrial_k < \tau_{k-1}] \\
    \geq&\ \P_k[g_k^Td_k + u_k^TH_ku_k \geq \nabla f(x_k)^T\dtrue_k + (\utrue_k)^T H_k \utrue_k | \tautruetrial_k < \tau_{k-1}] \\
       =&\ \P_k[g_k^Td_k + u_k^TH_ku_k \geq \nabla f(x_k)^T\dtrue_k + (\utrue_k)^T H_k \utrue_k] \geq p,
  \end{align*}
  as desired. \Halmos
\endproof

We can now prove the following result for our algorithm.  (We remark that \cite{BeraCurtRobiZhou21} also discusses an illustrative example in which \eqref{eq.daniel_special_label} holds for all $k \in \N{}$; see Example~3.17 in that article.)

\begin{proposition}
  If, in a given run of our algorithm, there exist $\tautruetrial_{\min} \in \R{}_{>0}$ and $p \in (0,1]$ such that $\tautruetrial_k \geq \tautruetrial_{\min}$ and \eqref{eq.daniel_special_label} hold for all $k \in \N{}$, then the probability is zero that the event $E_{\tau,\big}(\tautruetrial_{\min})$ occurs in the run.
\end{proposition}
\proof{Proof.}
  Under the conditions of the proposition, the conclusion follows from Lemma~\ref{lem:tautrialextra} using the same argument as in the proof of \cite[Proposition~3.16]{BeraCurtRobiZhou21}. \Halmos
\endproof

The analysis above shows that if $\{\tautruetrial_k\}$ is bounded below uniformly by a positive real number, then the probability is zero that $E_{\tau,\big}(\tautruetrial_{\min})$ occurs in a given run.  From this property, it follows under this condition that the probability is zero that $E_{\tau,\big}(\tautruetrial_{\min})$ occurs in a countable number of runs.  However, this analysis does not address what may be the total probability, over all possible runs of the algorithm, that $E_{\tau,\big}(\tautruetrial_{\min})$ may occur.  (To understand this, recognize that a given run of the algorithm may be akin to a single realization from a continuous probability distribution.  Since the probability of any given realization is zero, one cannot simply take the fact that the probability of $E_{\tau,\big}(\tautruetrial_{\min})$ occurring in a given run is zero to imply that the probability of such an event is zero over all possible runs---since there may be an uncountable number of them.  Hence, an alternative approach needs to be taken.)  Proving that, under certain assumptions, the total probability is zero that this event occurs requires careful consideration of the stochastic process generated by the algorithm, and in particular consideration of the filtration defined by the initial conditions and the sequence of stochastic gradient estimates that are generated by the algorithm.  We prove in Appendix~\ref{sec.totalprob} a formal version of the following informally written theorem.

\btheorem[Informal version of Theorem~\ref{th.totalprob} in Appendix~\ref{sec.totalprob}]\label{th.tau_zero_informal}
  If the true trial merit parameter sequence is bounded below by a positive real number and there exists $p \in (0,1]$ such that a condition akin to \eqref{eq.daniel_special_label} always holds, then the total probability of the event that the merit parameter sequence eventually remains constant at too large of a value (as in Definition~\ref{def.tau_big}) is zero.
\etheorem

The key to our proof of Theorem~\ref{th.totalprob} is the construction of a tree to characterize the stochastic process generated by the algorithm in a manner that one can employ the multiplicative form of Chernoff's bound to capture the probability of having repeated \emph{missed} opportunities to decrease the merit parameter when it would have been reduced if the true gradients were computed.

%************
% Subsection
%************
\subsection{Complementary Events}\label{sec.diverge}

Our analyses in Sections~\ref{sec.tau_low}, \ref{sec.tau_zero}, and \ref{sec.tau_big} do not cover all possible events.  Ignoring events in which the stochastic gradients are biased and/or have unbounded variance, the events that complement $E_{\tau,\low}$, $E_{\tau,\zero}$, and $E_{\tau,\big}$ are the following:
\bitemize
  \item $E_{\tau,\zero,\text{bad}}$: $\{\tau_k\} \searrow 0$ and for all $M \in \R{}_{>0}$ there exists $k \in \N{}$ such that $\|g_k - \nabla f(x_k)\|_2^2 > M$;
  \item $E_{\tau,\big,\text{bad}}$: $\{\tautruetrial_k\} \searrow 0$ and there exists $\tau_{\big} \in \R{}_{>0}$ such that $\tau_k = \tau_{\big}$ for all $k \in \N{}$.
\eitemize

The event $E_{\tau,\zero,\text{bad}}$ represents cases in which the merit parameter vanishes while the stochastic gradient estimates do not remain in a bounded set.  The difficulty of proving a guarantee for this setting can be seen as follows.  If the merit parameter vanishes, then this is an indication that less emphasis should be placed on the objective over the course of the optimization process, which may indicate that the constraints are infeasible or degenerate.  However, if a subsequence of stochastic gradient estimates diverges at the same time, then each large (in norm) stochastic gradient estimate may suggest that a significant amount of progress can be made in reducing the objective function, despite the merit parameter having reached a small value (since it is vanishing).  This disrupts the balance that the merit parameter attempts to negotiate between the objective and the constraint violation terms in the merit function.  Our analysis of the event $E_{\tau,\zero}$ in Section~\ref{sec.tau_zero} shows that if the stochastic gradient estimates remain bounded, then the algorithm can effectively transition to solving the deterministic problem of minimizing constraint violation.  However, it remains an open question whether it is possible to obtain a similar guarantee if/when a subsequence of stochastic gradient estimates diverges.  Ultimately, one can argue that scenarios of unbounded noise, such as described here, might only be of theoretical interest rather than real, practical interest.  For instance, if $f$ is defined by a (large) finite sum of component functions whose gradients (evaluated at points in a set containing the iterates) are always contained in a ball of uniform radius about the gradient of $f$---a common scenario in practice---then $E_{\tau,\zero,\text{bad}}$ cannot occur.

Now consider the event $E_{\tau,\big,\text{bad}}$.  We have shown in Section~\ref{sec.tau_big} that under certain conditions, including if $\{\tautruetrial_k\}$ is bounded below by $\tautruetrial_{\min} \in \R{}_{>0}$, then $E_{\tau,\big}$ occurs with probability zero.  However, this does not account for situations in which $\{\tautruetrial_k\}$ vanishes while $\{\tau_k\}$ does not.  Nonetheless, we contend that $E_{\tau,\big,\text{bad}}$ can be ignored for practical purposes since the adverse effect that it may have on the algorithm is observable.  In particular, if the merit parameter remains fixed at a value that is too large, then the worst that may occur is that $\{\|J_k^Tc_k\|_2\}$ does not vanish.  A practical implementation of the algorithm would monitor this quantity in any case (since, by  Corollary~\ref{cor.stochastic_tau_finite_large_pi_finite}, even in $E_{\tau,\low}$ one only knows that the limit inferior of the expectation of $\{\|J_k^Tc_k\|_2\}$ vanishes) and reduce the merit parameter if progress toward reducing constraint violation is inadequate.  Hence, $E_{\tau,\big,\text{bad}}$ (and $E_{\tau,\big}$ for that matter) is an event that at most suggests practical measures of the algorithm that should be employed for $E_{\tau,\low}$ in any case.

%*********
% Section
%*********
\section{Numerical Experiments}\label{sec.numerical}

The goal of our numerical experiments is to compare the empirical performance of our proposed stochastic SQP method (Algorithm~\ref{alg.sqp}) against some alternative approaches on problems from a couple of test set collections.  We implemented our algorithm in Matlab.  Our code is publicly available.\footnote{\href{https://github.com/frankecurtis/StochasticSQP}{\texttt{https://github.com/frankecurtis/StochasticSQP}}}  We first consider equality constrained problems from the CUTEst collection \cite{GoulOrbaToin15}, then consider two types of constrained logistic regression problems with datasets from the LIBSVM collection \cite{chang2011libsvm}.  We compare the performance of our method versus a stochastic subgradient algorithm \cite{DaviDrusKaka20} employed to minimize the exact penalty function \eqref{eq.merit} and, in one set of our logistic regression experiments where it is applicable, versus a stochastic projected gradient method.  These algorithms were chosen since, like our method, they operate in the highly stochastic regime.  We do not compare against the aforementioned method from \cite{NaAnitKola21} since, as previously mentioned, that approach may refine stochastic gradient estimates during each iteration as needed by a line search.  Hence, that method offers different types of convergence guarantees and is not applicable in our regime of interest.

In all of our experiments, results are given in terms of feasibility and stationarity errors at the \emph{best} iterate, which is determined as follows.  If, for a given problem instance, an algorithm produced an iterate that was sufficiently feasible in the sense that $\|c_k\|_{\infty} \leq 10^{-6} \max \{1, \|c_0\|_{\infty}\}$ for some $k \in \N{}$, then, with the largest $k \in \N{}$ satisfying this condition, the feasibility error was reported as $\|c_k\|_\infty$ and the stationarity error was reported as $\|\nabla f(x_k) + J_k^Ty_k\|_{\infty}$, where $y_k$ was computed as a least-squares multiplier using the true gradient $\nabla f(x_k)$ and $J_k$.  (The multiplier $y_k$ and corresponding stationarity error are not needed by our algorithm; they are computed merely so that we could record the error for our experimental results.)  If, for a given problem instance, an algorithm did not produce a sufficiently feasible iterate, then the feasibility and stationarity errors were computed in the same manner at the least infeasible iterate (with respect to the measure of infeasibility $\|\cdot\|_\infty$).

%************
% Subsection
%************
\subsection{Implementation Details}

For all methods, Lipschitz constant estimates for the objective gradient and constraint Jacobian---playing the roles of $L$ and $\Gamma$, respectively---were computed using differences of gradients near the initial point.  Once these values were computed, they were kept constant for all subsequent iterations.  This procedure was performed in such a way that, for each problem instance, all algorithms used the same values for these estimates.

As mentioned in Section~\ref{sec.algorithm}, there are various extensions of our step size selection scheme with which one can prove, with appropriate modifications to our analysis, comparable convergence guarantees as are offered by our algorithm.  We included one such extension in our software implementation for our experiments.  In particular, in addition to $\alphasuff_k$ in \eqref{eq.alpha_suff}, one can directly consider the upper bound in \eqref{eq.merit_decrease} with the gradient $\nabla f(x_k)$ replaced by its estimate $g_k$, i.e.,
\bequationNN
  \baligned
    &\ \alpha \tau_k g_k^Td_k + |1-\alpha|\|c_k\|_2 - \|c_k\|_2 + \alpha\| c_k + J_k d_k\|_2 + \thalf (\tau_k L + \Gamma) \alpha^2 \|d_k\|_2^2 \\
    =&\ -\alpha \Delta l(x_k,\tau_k,g_k,d_k) + |1 - \alpha| \|c_k\|_2 - (1 - \alpha) \|c_k\|_2 + \thalf (\tau_k L + \Gamma) \alpha^2 \|d_k\|_2^2, \\
  \ealigned
\eequationNN
and consider the step size that minimizes this as a function of $\alpha$ (with scale factor $\beta_k$), namely,
\bequation\label{eq.alpha_min}
  \alphamin_k := \max\left\{\min\left\{ \frac{\beta_k \Delta l(x_k,\tau_k,g_k,d_k)}{(\tau_k L + \Gamma) \|d_k\|_2^2} , 1 \right\} , \frac{\beta_k \Delta l(x_k,\tau_k,g_k,d_k) - 2\|c_k\|_2}{(\tau_k L + \Gamma) \|d_k\|_2^2} \right\}.
\eequation
(Such a value is used in \cite{BeraCurtRobiZhou21}.)  The algorithm can then set a trial step size as any satisfying
\bequation\label{eq.alpha_trial}
  \alphatrial_k \in [ \min\{\alphasuff_k,\alphamin_k\},\max\{\alphasuff_k,\alphamin_k\}]
\eequation
and set $\alpha_k$ as the projection of this value, rather than $\alphasuff_k$, for all $k \in \N{}$.  (The projection interval in \eqref{eq.alpha} should be modified, specifically with each instance of $2(1-\eta)$ replaced by $\min\{2(1-\eta),1\}$, to account for the fact that the lower value in \eqref{eq.alpha_trial} may be smaller than $\alphasuff_k$.  A similar modification is needed in the analysis, specifically in the requirements for $\{\beta_k\}$ in Lemma~\ref{lem.key_decrease}.)

One can also consider rules that allow even larger step sizes to be taken.  For example, rather than consider the upper bound offered by the last expression in \eqref{eq.merit_decrease}, one can consider any step size that ensures that the penultimate expression in \eqref{eq.merit_decrease} is less than or equal to the right-hand side of \eqref{eq.suff_dec} with $\nabla f(x_k)$ replaced by $g_k$.  Such a value can be found with a one-dimensional search over~$\alpha$ with negligible computational cost.  Our analysis can be extended to account for this option as well.  However, for our experimental purposes here, we do not consider such an approach.

For our stochastic SQP method, we set $H_k \gets I$ and $\alphatrial_k \gets \max\{\alphasuff_k,\alphamin_k\}$ for all $k \in \N{}$.  Other parameters were set as $\tau_{-1} \gets 1$, $\chi_{-1} \gets 10^{-3}$, $\zeta_{-1} \gets 10^{3}$, $\xi_{-1} \gets 1$, $\omega \gets 10^2$, $\epsilon_v \gets 1$, $\sigma \gets 1/2$, $\epsilon_{\tau} \gets 10^{-2}$, $\epsilon_{\chi} \gets 10^{-2}$, $\epsilon_{\zeta} \gets 10^{-2}$, $\epsilon_{\xi} \gets 10^{-2}$, $\eta \gets 1/2$, and $\theta \gets 10^4$.  For the stochastic subgradient method, the merit parameter value and step size were tuned for each problem instance, and for the stochastic projected gradient method, the step size was tuned for each problem instance; details are given in the following subsections.  In all experiments, both the stochastic subgradient and stochastic projected gradient method were given many more iterations to find each of their best iterates for a problem instance; this is reasonable since the search direction computation for our method is more expensive than for the other methods.  Again, further details are given below.

%************
% Subsection
%************
\subsection{CUTEst problems}\label{sec.cutest}

In our first set of experiments, we consider equality constrained problems from the CUTEst collection.  Specifically, of the $136$ such problems in the collection, we selected those for which $(i)$ $f$ is not a constant function, and $(ii)$ $n + m + 1 \leq 1000$. This selection resulted in a set of $67$ problems.  In order to consider the context in which the LICQ does not hold, for each problem we duplicated the last constraint.  (This does not affect the feasible region nor the set of stationary points, but ensures that the problem instances are degenerate.)  Each problem comes with an initial point, which we used in our experiments.  To make each problem stochastic, we added noise to each gradient computation.  Specifically, for each run of an algorithm, we fixed a \emph{noise level} as $\epsilon_N \in \{10^{-8}, 10^{-4}, 10^{-2}, 10^{-1}\}$, and in each iteration set the stochastic gradient estimate as $g_k \gets \Ncal(\nabla f(x_k),\epsilon_N I)$.  For each problem and noise level, we ran $10$ instances with different random seeds.  This led to a total of $670$ runs of each algorithm for each noise level.

We set a budget of 1000 iterations for our stochastic SQP algorithm and a more generous budget of 10000 iterations for the stochastic subgradient method.  We followed the same strategy as in~\cite{BeraCurtRobiZhou21} to tune the merit parameter $\tau$ for the stochastic subgradient method, but also tuned the step sizes through the sequence $\{\beta_k\}$. Specifically, for each problem instance, we ran the stochastic subgradient method for 11 different values of $\tau$ and 4 different values of $\beta$, namely, $\tau \in \{10^{-10},10^{-9},\dots,10^{0}\}$ and $\beta \in \{10^{-3},10^{-2},10^{-1},10^0\}$, set the step size as $\tfrac{\beta \tau}{\tau L + \Gamma}$, and selected the combination of $\tau$ and $\beta$ for that problem instance that led to the best iterate overall.  (We found through this process that the selected $(\tau,\beta)$ pairs were relatively evenly distributed over their ranges, meaning that this extensive tuning effort was useful to obtain better results for the stochastic subgradient method.)  For our stochastic SQP method, we set $\beta_k \gets 1$ for all $k \in \N{}$.  Overall, between the additional iterations allowed in each run of the stochastic subgradient method, the different merit parameter values tested, and the different step sizes tested, the stochastic subgradient method was given $440$ times the number of iterations that were given to our stochastic SQP method for each problem.

\begin{figure}[ht]
  \centering
  \includegraphics[width=0.4\textwidth,clip=true,trim=20 5 90 50]{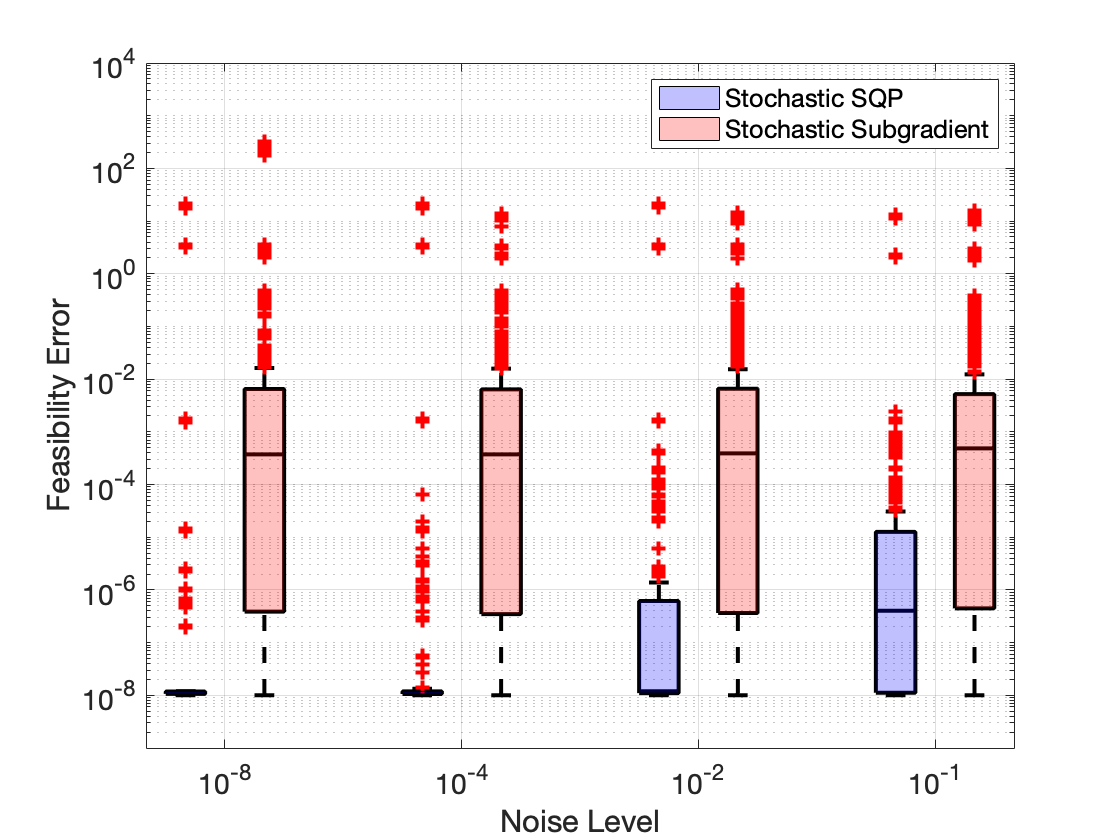}\quad
  \includegraphics[width=0.4\textwidth,clip=true,trim=20 5 90 50]{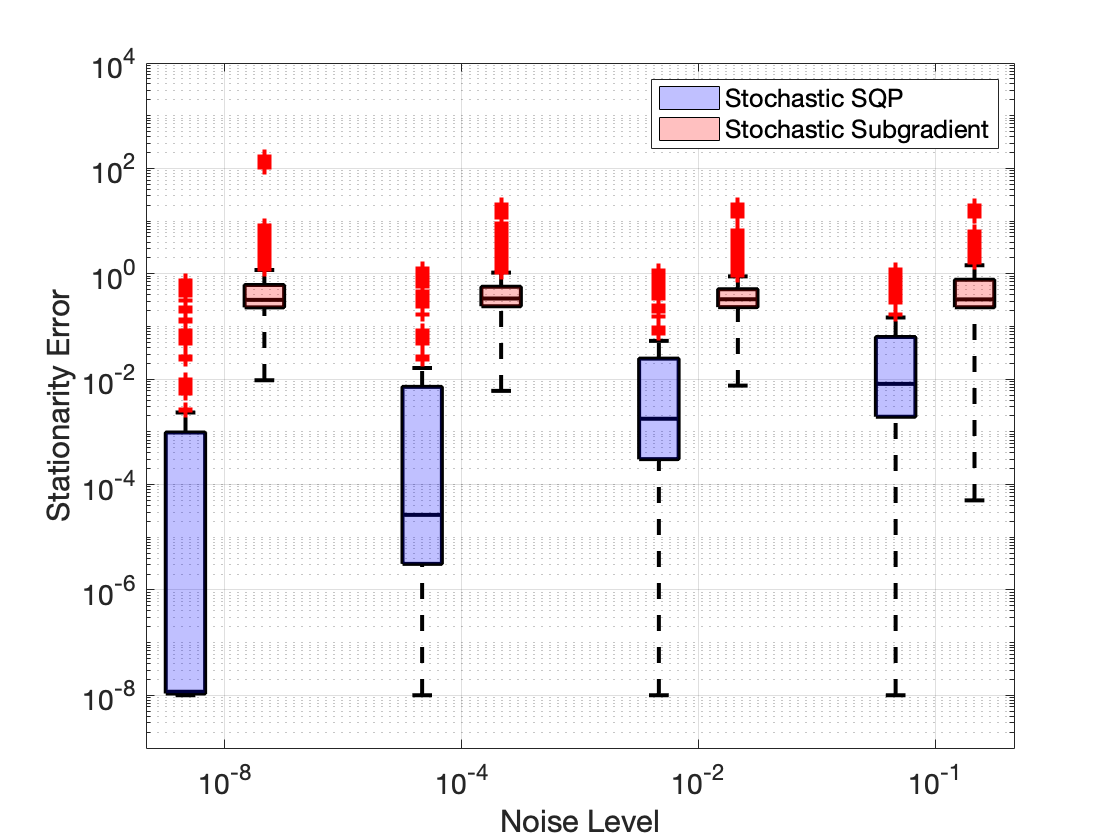}
  \caption{Box plots for feasibility errors (left) and stationarity errors (right) when our stochastic SQP method and a stochastic subgradient method are employed to solve equality constrained problems from the CUTEst collection.}
  \label{fig.perf_cutest}
\end{figure}

The results of this experiment are reported in the form of box plots in Figure~\ref{fig.perf_cutest}.  One finds that the best iterates from our stochastic SQP algorithm generally correspond to much lower feasibility and stationarity errors for all noise levels.  The stationarity errors for our method degrade as the noise level increases, but this is not surprising since these experiments are run with $\{\beta_k\}$ being a constant sequence.  It is interesting, however, that our algorithm typically finds iterates that are sufficiently feasible, even for relatively high noise levels.  This shows that our approach handles the deterministic constraints well despite the stochasticity of the objective gradient estimates.  Finally, we remark that for these experiments our algorithm found $\tau_{k-1} \leq \tautruetrial_k$ to hold in roughly 98\% of all iterations for all runs (across all noise levels), and found this inequality to hold in the last 50 iterations in 100\% of all runs.  This provides evidence for our claim that the merit parameter not reaching a sufficiently small value is not an issue of practical concern.

%************
% Subsection
%************
\subsection{Constrained Logistic Regression}\label{sec.logistic} 

In our next sets of experiments, we consider equality constrained logistic regression problems of the form
\bequation\label{eq.logistic}
  \min_{x \in \R{n}}\ f(x) = \frac{1}{N} \sum_{i=1}^N \log\left(1 + e^{-y_i(X_i^Tx)} \right)\ \ \st\ \ Ax = b,\ \ \|x\|_2^2 = 1,
\eequation
where $X \in \R{n\times N}$ contains feature data for $N$ data points (with $X_i$ representing the $i$th column of~$X$), $y \in \{-1,1\}^N$ contains corresponding label data, $A \in \R{(m+1) \times n}$ and $b \in \R{m+1}$.  For instances of $(X,y)$, we consider $11$ binary classification datasets from the LIBSVM collection \cite{chang2011libsvm}; specifically, we consider all of the datasets for which $12 \leq n \leq 1000$ and $256 \leq N \leq 100000$.  (For datasets with multiple versions, e.g., the $\{\texttt{a1a},\dots,\texttt{a9a}\}$ datasets, we consider only the largest version.)  The names of the datasets that we used and their sizes are given in Table~\ref{tab:data}.  For the linear constraints, we generated random $A$ and $b$ for each problem.  Specifically, the first $m = 10$ rows of $A$ and first $m$ entries in $b$ were set as random values with each entry being drawn from a standard normal distribution.  Then, to ensure that the LICQ was not satisfied (at any algorithm iterate), we duplicated the last constraint, making $m+1$ linear constraints overall.  For all problems and algorithms, the initial iterate was set to the vector of all ones of appropriate dimension.

\begin{table}[ht]
  \caption{Names and sizes of datasets. (Source: \cite{chang2011libsvm}.)}
  \label{tab:data}
  \centering
  {\footnotesize
\begin{tabular}{lcc}\toprule
\textbf{dataset}    & \textbf{dimension ($n$)} & \textbf{ datapoints ($N$)} \\ \hline
\texttt{a9a}             & $123$                          & $32,561$                                                                             \\ \hdashline
\texttt{australian}      & $14$                           & $690$                                                                                \\ \hdashline
%\texttt{breast-cancer}   & $10$                           & $683$                                                                                \\
%\texttt{cod-rna}         & $8$                            & $59.535$                                                                             \\
%\texttt{colon-cancer}    & $2,000$                        & $62$                                                                                 \\
%\texttt{diabetes}        & $8$                            & $768$                                                                                \\
%\texttt{gisette}         & $5,000$                        & $6,000$                                                                              \\ \hdashline
\texttt{heart}           & $13$                           & $270$                                                                                \\ \hdashline
\texttt{ijcnn1}          & $22$                           & $49,990$                                                                             \\ \hdashline
\texttt{ionosphere}      & $34$                           & $351$                                                                                \\ \hdashline
%\texttt{liver-disorders} & $5$                            & $145$                                                                                \\
\texttt{madelon}         & $500$                          & $2,000$                                                                              \\ \hdashline
\texttt{mushrooms}       & $112$                          & $8,124$                                                                              \\ \hdashline
\texttt{phising}         & $68$                           & $11,055$                                                                             \\ \hdashline
\texttt{sonar}           & $60$                           & $208$                                                                                \\ \hdashline
\texttt{splice}          & $60$                           & $1,000$                                                                              \\ \hdashline
\texttt{w8a}             & $300$                          & $49,749$       \\
\hline%\bottomrule                     
\end{tabular}}
\end{table}

For one set of experiments, we consider problems of the form~\eqref{eq.logistic} except without the norm constraint.  For this set of experiments, the performance of all three algorithms---stochastic SQP, subgradient, and projected gradient---are compared.  For each dataset, we considered two noise levels, where the level is dictated by the mini-batch size of each stochastic gradient estimate (recall~\eqref{eq.mini-batch}).  For the mini-batch sizes, we employed $b_k \in \{16,128\}$ for all problems.  For each dataset and mini-batch size, we ran $5$ instances with different random seeds.

A budget of $5$ epochs (i.e., number of effective passes over the dataset) was used for all methods.   For our stochastic SQP method, we used $\beta_k = 10^{-1}$ for all $k \in \N{}$.  For the stochastic subgradient method, the merit parameter and step size were tuned like in Section~\ref{sec.cutest} over the sets $\beta \in \{10^{-3},10^{-2},10^{-1},10^0\}$ and $\tau \in \{10^{-3},10^{-2},10^{-1},10^0\}$.  For the stochastic projected gradient method, the step size was tuned using the formula $\tfrac{\beta}{L}$ over $\beta \in \{10^{-8},10^{-7},\dots,10^{1},10^{2}\}$.  Overall, this meant that the stochastic subgradient and stochastic projected gradient methods were effectively run for $16$ and $11$ times the number of epochs, respectively, that were allowed for our method.

The results for this experiment are reported in Table~\ref{tab:logistic_regression1}.  For every dataset and mini-batch size, we report the average feasibility and stationarity errors for the best iterates of each run along with a $95\%$ confidence interval.  The results show that our method consistently outperforms the two alternative approaches despite the fact that each of the other methods were tuned with various choices of the merit and/or step size parameter.  For a second set of experiments, we consider problems of the form~\eqref{eq.logistic} with the norm constraint.  The settings for the experiment were the same as above, except that the stochastic projected gradient method is not considered.  The results are stated in Table~\ref{tab:logistic_regression2}.  Again, our method regularly outperforms the stochastic subgradient method in terms of the best iterates found.  For the experiments without the norm constraint, our algorithm found $\tau_{k-1} \leq \tautruetrial_k$ to hold in roughly 98\% of all iterations for all runs, and found this inequality to hold in all iterations in the last epoch in 100\% of all runs.  With the norm constraint, our algorithm found $\tau_{k-1} \leq \tautruetrial_k$ to hold in roughly 97\% of all iterations for all runs, and found this inequality to hold in all iterations in the last epoch in 99\% of all runs.

\begin{table}[tb]
  \caption{Average feasibility and stationarity errors, along with 95\% confidence intervals, when our stochastic SQP method, a stochastic subgradient method, and a stochastic projected gradient method are employed to solve logistic regression problems with linear constraints (only).  The results for the best-performing algorithm are shown in bold.}
  \label{tab:logistic_regression1}
\centering
{\tiny
\begin{tabular}{l|c|c|c|c|c|c} \toprule
                                &                                 & \multicolumn{2}{c|}{\begin{tabular}[c]{@{}c@{}}Stochastic\\ Subgradient\end{tabular}} & \multicolumn{1}{c|}{\begin{tabular}[c]{@{}c@{}}Stochastic \\ Projected Gradient\end{tabular}} & \multicolumn{2}{c}{\begin{tabular}[c]{@{}c@{}}Stochastic \\ SQP \end{tabular}}                               \\ \hline
\multicolumn{1}{l|}{dataset}  & \multicolumn{1}{c|}{batch} & \multicolumn{1}{c|}{Feasibility}           & \multicolumn{1}{c|}{Stationarity}          & \multicolumn{1}{c|}{Stationarity}              & \multicolumn{1}{c|}{Feasibility} & \multicolumn{1}{c}{Stationarity} \\ \hline
\texttt{a9a} & 16                     &       $8.30e-03   \pm 2.32e-03             $                    &     $      1.64e-01   \pm 3.55e-03          $                                                                 &     $ 3.64e-02   \pm 2.95e-03           $                            &          $     \pmb{1.22e-15   \pm 2.18e-16}       $            &        $\pmb{9.99e-03 \pm 6.92e-03 }$                      \\ \hdashline
\texttt{a9a} & 128                     &      $1.16e-02   \pm 4.60e-05           $                      &        $   1.69e-01   \pm 2.51e-02         $                                                                  &     $ 1.69e-02   \pm 2.79e-03            $                             &             $  \pmb{1.64e-15   \pm 4.00e-16 }     $             &      $   \pmb{7.33e-03   \pm 4.68e-05}     $                  \\ \hline
\texttt{australian} & 16                     &    $   7.94e-02   \pm 1.60e-05       $                          &       $    7.94e-02   \pm 1.60e-05       $                                                                     &  $    9.17e-02   \pm 4.32e-04   $                                    &     $         \pmb{ 5.72e-06   \pm 1.56e-06    }     $          &       $  \pmb{2.67e-02   \pm 6.43e-04   }  $                  \\ \hdashline
\texttt{australian} & 128                     &    $   5.02e-01   \pm 7.04e-05      $                           &       $    5.02e-01   \pm 7.04e-05    $                                                                        &    $  1.11e-02   \pm 7.19e-05       $                                 &      $       \pmb{  6.58e-05   \pm 7.90e-07   }   $             &     $   \pmb{ 5.50e-02   \pm 1.08e-03}   $                    \\ \hline
%\texttt{gisette} & 16                     &       XXXX   \pm XXXX                                 &           XXXX   \pm XXXX                               &    XXXX   \pm XXXX                                            &      XXXX   \pm XXXX                                        &               XXXX   \pm XXXX                   &         XXXX   \pm XXXX                       \\ \hdashline
%\texttt{gisette} & 128                     &       XXXX   \pm XXXX                                 &           XXXX   \pm XXXX                               &    XXXX   \pm XXXX                                            &      XXXX   \pm XXXX                                        &               XXXX   \pm XXXX                   &         XXXX   \pm XXXX                       \\ \hline
\texttt{heart} & 16                     &      $ 3.66e-01   \pm 4.37e-03           $                      &        $   3.28e+01   \pm 7.02e+00    $                                                                        &   $   \pmb{3.17e+01   \pm 6.72e+00   }        $                             &            $  \pmb{ 8.83e-03   \pm 2.77e-03       } $           &       $  3.39e+01   \pm 9.85e+00     $                  \\ \hdashline
\texttt{heart} & 128                     &    $   1.52e+00   \pm 4.96e-02        $                         &     $      \pmb{1.23e+01   \pm 1.40e+01 }   $                                                                        &   $   3.29e+01   \pm 3.21e+00        $                                &       $     \pmb{   1.26e-01   \pm 7.86e-04   }   $             &    $     3.24e+01   \pm 1.76e+00     $                  \\ \hline
\texttt{ijccn1} & 16                     &    $   3.58e-03   \pm 2.00e-05           $                      &        $   4.70e-02   \pm 6.45e-07           $                                                                 &   $   7.41e-02   \pm 3.33e-07          $                              &           $  \pmb{  3.03e-15   \pm 6.20e-16    }     $          &       $ \pmb{ 1.93e-03   \pm 4.07e-06 }        $              \\ \hdashline
\texttt{ijccn1} & 128                     &   $    3.90e-02   \pm 4.01e-06            $                     &        $   5.17e-02   \pm 1.65e-07            $                                                             &    $  3.88e-02   \pm 6.15e-07             $                           &             $ \pmb{ 2.16e-09   \pm 2.62e-09     }   $           &       $ \pmb{ 1.70e-02   \pm 5.19e-05 }          $            \\ \hline
\texttt{ionosphere} & 16                     &    $   5.41e-01   \pm 8.80e-05         $                        &         $  5.41e-01   \pm 8.80e-05        $                                                                   &      $  9.77e-01   \pm 8.55e-03        $                                &           $    \pmb{ 9.61e-07   \pm 2.77e-09    } $              &       $ \pmb{  4.17e-02   \pm 1.08e-03    }$                   \\ \hdashline
\texttt{ionosphere} & 128                     &    $   5.76e+00   \pm 3.76e-05        $                         &       $    5.76e+00   \pm 3.76e-05    $                                                                        &   $   5.98e+00   \pm 3.21e-03  $                                      &      $     \pmb{     1.31e-05   \pm 1.14e-09  }  $               &      $  \pmb{  1.55e-01   \pm 2.61e-03  } $                    \\ \hline
\texttt{madelon} & 16                     &     $  3.06e-02   \pm 1.85e-02           $                      &       $    5.46e+01   \pm 1.25e+01           $                                                                &      $2.11e+01   \pm 2.72e+00           $                             &            $   \pmb{2.88e-08   \pm 5.51e-08     }   $           &        $ \pmb{1.09e+01   \pm 3.00e+00   }     $               \\ \hdashline
\texttt{madelon} & 128                     &      $ 1.87e+00   \pm 7.62e-01         $                        &       $    2.21e+01   \pm 1.55e+01      $                                                                      &   $  \pmb{ 2.16e+01   \pm 4.17e+00  }   $                                    &        $     \pmb{  5.81e-01   \pm 1.63e-02 }     $             &      $   4.81e+01   \pm 4.75e+00     $                  \\ \hline
\texttt{mushrooms} & 16                     &    $   2.19e-01   \pm 6.55e-04        $                         &      $     2.19e-01   \pm 6.55e-04       $                                                                    &     $ 7.31e-03   \pm 3.21e-06     $                                   &         $     \pmb{  2.08e-15   \pm 3.28e-16     }   $           &       $  \pmb{5.95e-03   \pm 3.21e-05     }  $                \\ \hdashline
\texttt{mushrooms} & 128                     &   $    4.73e-01   \pm 4.37e-05     $                            &     $      4.73e-01   \pm 4.37e-05     $                                                                       & $     3.31e-02   \pm 7.13e-05 $                                       &     $     \pmb{      1.66e-09   \pm 6.20e-14 }   $               &    $   \pmb{  3.28e-02   \pm 9.15e-04}   $                    \\ \hline
\texttt{phishing} & 16                     &     $  2.67e-02   \pm 2.76e-07         $                        &         $  3.47e-02   \pm 1.39e-09         $                                                                   &    $   \pmb{2.20e-05   \pm 9.29e-06}   $                                    &            $   \pmb{ 4.26e-15   \pm 1.27e-15  }    $             &      $   3.37e-03   \pm 1.27e-06      $                 \\ \hdashline
\texttt{phishing} & 128                     &     $  3.06e-01   \pm 1.13e-06      $                           &      $     3.06e-01   \pm 1.13e-06   $                                                                        &   $    2.29e-01   \pm 8.88e-03   $                                    &          $    \pmb{  1.83e-15   \pm 4.99e-16 }    $              &     $   \pmb{ 2.20e-02   \pm 7.29e-03  }   $                  \\ \hline
\texttt{sonar} & 16                     &      $ 1.33e+00   \pm 1.08e-04          $                       &         $  1.33e+00   \pm 1.08e-04       $                                                                     &     $ 6.13e-01   \pm 2.22e-03     $                                   &            $    \pmb{7.02e-07   \pm 1.60e-07   }   $             &     $  \pmb{  2.34e-02   \pm 2.03e-04   }  $                  \\ \hdashline
\texttt{sonar} & 128                     &    $   1.33e+01   \pm 1.48e-04       $                          &         $  1.33e+01   \pm 1.48e-04     $                                                                     &    $  6.46e-02   \pm 4.73e-03   $                                     &           $   \pmb{  2.07e-06   \pm 6.70e-10   }   $             &      $  \pmb{ 2.98e-02   \pm 1.71e-03   } $                   \\ \hline
\texttt{splice} & 16                     &      $ 2.56e-03   \pm 3.39e-04        $                         &         $  4.56e-01   \pm 3.55e-02      $                                                                      &     $ 9.65e-01   \pm 3.19e-03    $                                    &           $   \pmb{  7.49e-14   \pm 1.03e-13   }   $             &      $  \pmb{ 2.19e-02   \pm 4.33e-03  }  $                   \\ \hdashline
\texttt{splice} & 128                     &    $   3.14e-01   \pm 1.09e-04      $                           &       $    4.83e-01   \pm 4.65e-05    $                                                                       &    $  1.23e+00   \pm 9.44e-05  $                                      &          $   \pmb{   3.54e-08   \pm 5.74e-09  }   $              &      $ \pmb{  1.07e-02   \pm 3.16e-04 } $                     \\ \hline
\texttt{w8a} & 16                     &      $ 2.38e-02   \pm 1.75e-03      $                           &         $  1.47e-01   \pm 1.89e-06       $                                                                  &      $ 9.85e-04   \pm 3.31e-05        $                                &              $  \pmb{7.35e-15   \pm 6.98e-16       }   $         &      $\pmb{   6.07e-05   \pm 6.46e-05    }  $                 \\ \hdashline
\texttt{w8a} & 128                     &     $  1.79e-02   \pm 1.25e-03    $                             &      $     1.49e-01   \pm 4.64e-03    $                                                                       &   $   3.41e-02   \pm 7.43e-03    $                                    &          $   \pmb{   5.96e-15   \pm 5.67e-16    }  $             &       $ \pmb{ 1.20e-03   \pm 1.85e-03}  $                     \\ \hline
\end{tabular}}
\end{table}

\begin{table}[tb]
  \caption{Average feasibility and stationarity errors, along with 95\% confidence intervals, when our stochastic SQP method and a stochastic subgradient method are employed to solve logistic regression problems with linear constraints and a squared $\ell_2$-norm constraint.  The results for the best-performing algorithm are shown in bold.}
  \label{tab:logistic_regression2}
\centering
{\tiny
\begin{tabular}{l|c|c|c|c|c} \toprule
                                &                                 & \multicolumn{2}{c|}{\begin{tabular}[c]{@{}c@{}}Stochastic\\ Subgradient\end{tabular}} & \multicolumn{2}{c}{\begin{tabular}[c]{@{}c@{}}Stochastic \\ SQP \end{tabular}}                               \\ \hline
\multicolumn{1}{l|}{dataset}  & \multicolumn{1}{c|}{batch} & \multicolumn{1}{c|}{Feasibility}                         & \multicolumn{1}{c|}{Stationarity}              & \multicolumn{1}{c|}{Feasibility} & \multicolumn{1}{c}{Stationarity} \\ \hline
\texttt{a9a} & 16                     &    $   4.62e-03   \pm 3.27e-04    $                                                                           &              $ 1.24e-01   \pm 7.52e-02       $            &         $\pmb{5.52e-05 \pm 5.04e-09}$        & $\pmb{6.07e-03 \pm 2.32e-05}$               \\ \hdashline
\texttt{a9a} & 128                     &      $ 4.27e-03   \pm 3.92e-04    $                                                                           &             $  1.90e-01   \pm 3.03e-03       $            &        $\pmb{6.38e-05 \pm 1.12e-08}$       & $\pmb{4.40e-03 \pm 1.41e-05}$              \\ \hline
\texttt{australian} & 16                     &    $   1.51e-01   \pm 1.07e-05      $                                                                         &        $       1.51e-01   \pm 1.07e-05   $                &      $\pmb{1.52e-04 \pm 5.58e-06}$       & $\pmb{5.65e-03 \pm 3.73e-05}$           \\ \hdashline
\texttt{australian} & 128                     &     $  3.96e-01   \pm 1.87e-04    $                                                                           &       $        3.96e-01   \pm 1.87e-04     $              &      $\pmb{3.83e-04 \pm 5.45e-05}$     & $\pmb{1.68e-02 \pm 3.29e-03}$            \\ \hline
%\texttt{gisette} & 16                     &       XXXX   \pm XXXX                                                                               &               XXXX   \pm XXXX                   &         XXXX   \pm XXXX        & XXXX   \pm XXXX               \\ \hdashline
%\texttt{gisette} & 16                     &       XXXX   \pm XXXX                                                                               &               XXXX   \pm XXXX                   &         XXXX   \pm XXXX        & XXXX   \pm XXXX               \\ \hline
\texttt{heart} & 16                     &    $   1.57e+00    \pm 5.76e-01 $                                                                              &           $    2.86e+01   \pm 1.00e+01       $            &       $\pmb{9.29e-01 \pm 3.47e-02}$        & $\pmb{2.65e+01 \pm 1.81e+01}$              \\ \hdashline
\texttt{heart} & 128                     &     $ \pmb{1.33e+00    \pm 6.69e-01} $                                                                         &              $\pmb{ 1.69e+01    \pm 2.23e+00  }$                 &         $1.88e+00 \pm 1.42e-01$        &$2.93e+00 \pm 1.26e+00$          \\ \hline
\texttt{ijcnn1} & 16                     &      $ 5.36e-02   \pm 9.37e-07         $                                                                      &           $    5.36e-02   \pm 9.37e-07   $                &         $\pmb{3.70e-02 \pm 9.24e-05}$        & $\pmb{4.60e-02 \pm 8.32e-03}$               \\ \hdashline
\texttt{ijcnn1} & 128                     &    $   5.41e-02   \pm 1.04e-06       $                                                                        &            $   5.41e-02   \pm 1.04e-06     $              &         $\pmb{3.64e-02 \pm 1.06e-04}$        & $\pmb{3.64e-02  \pm  1.06e-04} $              \\ \hline
\texttt{ionosphere} & 16                     &    $   3.35e-01   \pm 1.06e-03     $                                                                          &        $       3.35e-01   \pm 1.06e-03     $              &         $\pmb{5.79e-03 \pm 1.44e-04}$        & $\pmb{1.21e-02 \pm 4.96e-03}$              \\ \hdashline
\texttt{ionosphere} & 128                     &   $    8.70e-01   \pm 1.43e-03   $                                                                            &       $        8.70e-01   \pm 1.43e-03     $              &         $\pmb{5.92e-03 \pm 2.18e-05}$        & $\pmb{4.31e-02 \pm 3.52e-04}$               \\ \hline
\texttt{madelon} & 16                     &      $ 2.66e+00   \pm 6.84e-01             $                                                                 &              $ 3.86e+01   \pm 3.28e+01        $           &         $\pmb{3.74e-01 \pm 8.55e-02}$        & $\pmb{4.70e-01 \pm 3.27e-02}$              \\ \hdashline
\texttt{madelon} & 128                     &       $ \pmb{2.21e+01    \pm 4.90e-01} $                                                                              &               $ \pmb{4.77e+01    \pm 4.84e+00} $                 &         $ 7.21e+01 \pm 5.28e+00$        & $7.21e+01 \pm 5.28e+00$               \\ \hline
\texttt{mushrooms} & 16                     &     $  1.01e-01   \pm 5.79e-05        $                                                                       &           $    1.55e-01   \pm 8.22e-06        $           &         $\pmb{4.06e-04 \pm 8.76e-09}$        & $\pmb{4.65e-03 \pm 3.65e-05}$              \\ \hdashline
\texttt{mushrooms} & 128                     &   $    9.72e-01   \pm 9.94e-06     $                                                                          &           $    9.72e-01   \pm 9.94e-06    $               &         $\pmb{6.96e-04 \pm 1.52e-09}$        & $\pmb{3.34e-03 \pm 2.35e-07}$               \\ \hline
\texttt{phishing} & 16                     &     $  1.30e-01   \pm 1.61e-06            $                                                                   &         $      1.30e-01   \pm 1.61e-06       $            &         $\pmb{3.65e-05 \pm 2.44e-08}$        &  $\pmb{8.17e-03 \pm 2.43e-05}$              \\ \hdashline
\texttt{phishing} & 128                     &    $   1.53e-01   \pm 3.37e-08         $                                                                      &         $      1.53e-01   \pm 3.37e-08       $            &         $\pmb{1.26e-04 \pm 3.30e-09}$       & $\pmb{8.45e-04 \pm 2.73e-07}$               \\ \hline
\texttt{sonar} & 16                     &      $ 6.45e-01   \pm 5.62e-04              $                                                                 &              $ 6.45e-01   \pm 5.62e-04           $        &         $\pmb{3.38e-03 \pm 8.81e-06}$       & $\pmb{1.48e-02 \pm 2.58e-04}$             \\ \hdashline
\texttt{sonar} & 128                     &     $  5.04e+00   \pm 4.44e-03           $                                                                    &             $  5.04e+00   \pm 4.44e-03       $            &         $\pmb{5.71e-03 \pm 8.61e-06}$        & $\pmb{2.16e-02 \pm 8.48e-05}$               \\ \hline
\texttt{splice} & 16                     &       $\pmb{1.96e-03    \pm 1.78e-04}$                                                                             &            $   4.94e-01   \pm 7.35e-03       $            &         $3.96e-03 \pm 7.12e-07$        & $\pmb{1.03e-02 \pm 1.14e-05}$              \\ \hdashline
\texttt{splice} & 128                     &    $   1.40e+00   \pm 7.90e-05          $                                                                     &       $        1.40e+00   \pm 7.90e-05      $             &         $\pmb{5.52e-03    \pm 3.72e-06} $       & $\pmb{1.04e-02    \pm 1.06e-04}  $             \\ \hline
\texttt{w8a} & 16                     &      $ 1.32e-02   \pm 6.83e-04              $                                                                 &             $  1.15e-01   \pm 1.33e-02           $        &         $\pmb{2.15e-04    \pm 2.24e-09} $       & $\pmb{1.83e-03    \pm 8.90e-07}  $             \\ \hdashline
\texttt{w8a} & 128                     &    $   5.35e-02   \pm 7.79e-02           $                                                                    &           $    1.33e-01   \pm 1.74e-07        $           &         $\pmb{1.67e-04    \pm 6.01e-09} $       & $\pmb{1.00e-03    \pm 1.01e-06}   $            \\ \hline
\end{tabular}}
\end{table}

%*********
% Section
%*********
\section{Conclusion}\label{sec.conclusion}

We have proposed, analyzed, and tested a stochastic SQP method for solving equality constrained optimization problems in which the objective function is defined by an expectation of a stochastic function.  Our algorithm is specifically designed for cases when the LICQ does not necessarily hold in every iteration.  The convergence guarantees that we have proved for our method consider situations when the merit parameter sequence eventually remains fixed at a value that is sufficiently small, in which case the algorithm drives stationarity measures for the constrained optimization problem to zero, and situations when the merit parameter vanishes, which may indicate that the problem is degenerate and/or infeasible.  Numerical experiments demonstrate that our algorithm consistently outperforms alternative approaches in the highly stochastic regime.

\appendix
\section{Deterministic Analysis}\label{sec.deterministic}
In this appendix, we prove that Theorem~\ref{th.deterministic} holds, where in particular we consider the context when $g_k = \nabla f(x_k)$ and $\beta_k = \beta$ satisfy \eqref{eq.beta} for all $k \in \N{}$.  For this purpose, we introduce a second termination condition in Algorithm~\ref{alg.sqp}.  In particular, after line~\ref{line.linsys}, we terminate the algorithm if both $\|g_k + J_k^Ty_k\|_2 = 0$ and $\|c_k\|_2 = 0$.  In this manner, if the algorithm terminates finitely, then it returns an infeasible stationary point (recall~\eqref{eq.infeasible_stationary}) or primal-dual stationary point for problem~\eqref{prob.opt} and there is nothing left to prove.  Hence, without loss of generality, we proceed under the assumption that the algorithm runs for all $k \in \N{}$.

Throughout our analysis in this appendix, we simply refer to the tangential direction as $u_k$, the full search direction as $d_k = v_k + u_k$, etc., even though it is assumed throughout this appendix that these are the \emph{true} quantities computed using the true gradient $\nabla f(x_k)$ for all $k \in \N{}$.

It follows in this context that both Lemma~\ref{lem.normal_component} and Lemma~\ref{lem.chi_upper} hold.  In addition, Lemma~\ref{lem.xi_lower} holds, where, in the proof, the case that $d_k = 0$ can be ignored due to the following lemma.

\blemma\label{lem.det.d_nonzero}
  For all $k \in \N{}$, one finds that $d_k = v_k + u_k \neq 0$.
\elemma
\proof{Proof.}
  For all $k \in \N{}$, the facts that $v_k \in \Range(J_k^T)$ and $u_k \in \Null(J_k)$ imply $d_k = v_k + u_k = 0$ if and only if $v_k = 0$ and $u_k = 0$.  Since we suppose in our analysis that the algorithm does not terminate finitely with an infeasible stationary point, it follows for all $k \in \N{}$ that $\|J_k^Tc_k\|_2 > 0$ or $\|c_k\|_2 = 0$.  If $\|J_k^Tc_k\|_2 > 0$, then Lemma~\ref{lem.normal_component} implies that $v_k \neq 0$, and the desired conclusion follows.  Hence, we may proceed under the assumption that $\|c_k\|_2 = 0$.  In this case, it follows under Assumption~\ref{ass.H} that $g_k + J_k^Ty_k = 0$ if and only if $u_k = 0$, which under our supposition that the algorithm does not terminate finitely means that $u_k \neq 0$. \Halmos
\endproof

We now prove a lower bound on the reduction in the merit function that occurs in each iteration.  This is a special case of Lemmas~\ref{lem.key_decrease} and \ref{lem.key_decrease_zero} for the deterministic setting.

\blemma\label{lem.det.merit_red}
  For all $k \in \N{}$, it holds that $\phi(x_k,\tau_k) - \phi(x_k + \alpha_k d_k,\tau_k) \geq \eta \alpha_k \Delta l(x_k,\tau_k,g_k,d_k)$.
\elemma
\proof{Proof.}
  For all $k \in \N{}$, it follows by the definition of $\alphasuff_k$ that (recall \eqref{eq.suff_dec})
  \bequationNN
    \phi(x_k + \alpha d_k,\tau_k) - \phi(x_k,\tau_k) \leq -\eta \alpha \Delta l(x_k,\tau_k,g_k,d_k)\ \ \text{for all}\ \ \alpha \in [0,\alphasuff_k].
  \eequationNN
  If $\|u_k\|_2^2 \geq \chi_k \|v_k\|_2^2$, then the only way that $\alpha_k > \alphasuff_k$ is if
  \bequationNN
    \frac{2(1-\eta)\beta\xi_k\tau_k}{\tau_k L + \Gamma} > \min\left\{ \frac{2(1-\eta)\beta \Delta l(x_k,\tau_k,g_k,d_k)}{(\tau_k L + \Gamma) \|d_k\|_2^2}, 1 \right\}.
  \eequationNN
  By \eqref{eq.beta}, the left-hand side of this inequality is less than 1, meaning $\alpha_k > \alphasuff_k$ only if
  \bequationNN
    \frac{2(1-\eta)\beta\xi_k\tau_k}{\tau_k L + \Gamma} > \frac{2(1-\eta)\beta \Delta l(x_k,\tau_k,g_k,d_k)}{(\tau_k L + \Gamma) \|d_k\|_2^2} \iff \xi_k \tau_k > \frac{\Delta l(x_k,\tau_k,g_k,d_k)}{\|d_k\|_2^2}.
  \eequationNN
  However, this is not true since $\xi_k \leq \xitrial_k$ for all $k \in \N{}$.  Following a similar argument for the case when $\|u_k\|_2^2 < \chi_k \|v_k\|_2^2$, the desired conclusion follows. \Halmos
\endproof

For our purposes going forward, let us define the shifted merit function $\tilde\phi : \R{n} \times \R{}_{\geq0} \to \R{}$ by
\bequationNN
  \tilde\phi(x,\tau) = \tau (f(x) - f_{\inf}) + \|c(x)\|_2.
\eequationNN

\blemma\label{lem.det.phi_shifted}
  For all $k \in \N{}$, it holds that $\tilde\phi(x_k,\tau_k) - \tilde\phi(x_{k+1},\tau_{k+1}) \geq \eta \alpha_k \Delta l(x_k,\tau_k,g_k,d_k)$.
\elemma
\proof{Proof.}
  For arbitrary $k \in \N{}$, it follows from Lemma~\ref{lem.det.merit_red} that
  \bequationNN
    \baligned
      \tau_{k+1}(f(x_k + \alpha_k d_k) - f_{\inf}) + \|c(x_k + \alpha_k d_k)\|_2
      \leq&\ \tau_k(f(x_k + \alpha_k d_k) - f_{\inf}) + \|c(x_k + \alpha_k d_k)\|_2 \\
      \leq&\ \tau_k(f(x_k) - f_{\inf}) + \|c_k\|_2 - \eta \alpha_k \Delta l(x_k,\tau_k,g_k,d_k),
    \ealigned
  \eequationNN
  from which the desired conclusion follows. \Halmos
\endproof

We now prove our first main result of this appendix.

\blemma\label{lem.det.Jc}
  The sequence $\{\|J_k^Tc_k\|_2\}$ vanishes.  Moreover, if there exist $k_J \in \N{}$ and $\sigma_J \in \R{}_{>0}$ such that the singular values of $J_k$ are bounded below by $\sigma_J$ for all $k \geq k_J$, then $\{\|c_k\|_2\}$ vanishes.
\elemma
\proof{Proof.}
  Let $\gamma \in \R{}_{>0}$ be arbitrary.  Our aim is to prove that the number of iterations with $x_k \in \Xcal_\gamma$ (recall \eqref{eq.X_gamma}) is finite.  Since $\gamma$ has been chosen arbitrarily in $\R{}_{>0}$, the conclusion will follow.  By Lemma~\ref{lem.tau_zero_alpha_low} and the fact that $\{\beta_k\}$ is chosen as a constant sequence, it follows that there exists $\underline\alpha \in \R{}_{>0}$ such that $\alpha_k \geq \underline\alpha$ for all $k \in \Kcal_\gamma$ (regardless of whether the search direction is tangentially or normally dominated).  Hence, using Lemmas~\ref{lem.normal_component} and \ref{lem.det.phi_shifted}, it follows that
  \bequationNN
    \tilde\phi(x_k,\tau_k) - \tilde\phi(x_{k+1},\tau_{k+1}) \geq \eta \underline\alpha \Delta l(x_k,\tau_k,g_k,d_k) \geq \eta \underline\alpha \sigma (\|c_k\|_2 - \|c_k + J_kv_k\|_2) \geq \eta \underline\alpha \sigma \kappa_v \kappa_c^{-1} \gamma^2.
  \eequationNN
  Hence, the desired conclusion follows since $\{\tilde\phi(x_k,\tau_k)\}$ is monotonically nonincreasing by Lemma~\ref{lem.det.phi_shifted} and is bounded below under Assumption~\ref{ass.main}. \Halmos
\endproof

We now show a consequence of the merit parameter eventually remaining constant.

\blemma\label{lem.det.last}
  If there exists $k_\tau \in \N{}$ and $\tau_{\min} \in \R{}_{>0}$ such that $\tau_k = \tau_{\min}$ for all $k \geq k_\tau$, then
  \bequationNN
    0 = \lim_{k\to\infty} \|u_k\|_2 = \lim_{k\to\infty} \|d_k\|_2 = \lim_{k\to\infty} \|g_k + J_k^Ty_k\|_2 = \lim_{k\to\infty} \|Z_k^Tg_k\|_2.
  \eequationNN
\elemma
\proof{Proof.}
  Under Assumption~\ref{ass.main} and the conditions of the lemma, Lemmas~\ref{lem.tau_zero_alpha_low} and \ref{lem.det.phi_shifted} imply that $\Delta l(x_k,\tau_k,g_k,d_k)\} \to 0$, which with \eqref{eq.model_reduction_condition} and Lemma~\ref{lem.normal_component} implies that $\{\|u_k\|_2\} \to 0$, $\{\|v_k\|_2\} \to 0$, and $\{\|J_k^Tc_k\|_2\} \to 0$.  The remainder of the conclusion follows from Assumption~\ref{ass.H} and \eqref{eq.linsys}. \Halmos
\endproof

The proof of Theorem~\ref{th.deterministic} can now be completed.

\proof{Proof of Theorem~\ref{th.deterministic}.}
  The result follows from Lemmas~\ref{lem.det.pi_bounded}, \ref{lem.det.Jc}, and \ref{lem.det.last}.
\endproof
\section{Total Probability Result}\label{sec.totalprob}
\newcommand{\Prob}{\mathbb{P}}
\newcommand{\hdelta}{\hat{\delta}}
\newcommand{\Dtrue}{D^{\text{\rm true}}}
\newcommand{\Utrue}{U^{\text{\rm true}}}
\newcommand{\Ytrue}{Y^{\text{\rm true}}}

\newcommand{\smax}{s_{\max}}
\newcommand{\kmax}{k_{\max}}

\newcommand{\Lbad}{\Lcal_{\text{\rm bad}}}
\newcommand{\Lgood}{\Lcal_{\text{\rm good}}}
\newcommand{\Ebad}{E_{\text{\rm bad}}}
\newcommand{\Ebadkmax}{E_{\text{\rm bad},\kmax}}
\newcommand{\EbadkmaxJ}{E_{\text{\rm bad},\kmax,J}}
\newcommand{\EbadB}{E_{\text{\rm bad},B}}
\newcommand{\Cdec}{C_{\text{\rm dec}}}
\newcommand{\Idec}{\Ical_{\text{\rm dec}}}
\newcommand{\ktau}{k_{\tau}}
\newcommand{\Kctau}{\Kcal_{\tau}}
\newcommand{\tcalktritrue}{\Tcal_k^{\text{\rm trial},\text{\rm true}}}

\newcommand{\tktritrue}{\tau_k^{\text{\rm trial},\text{\rm true}}}
\newcommand{\Ttritrue}{\Tcal^{\text{\rm trial},\text{\rm true}}}
\newcommand{\Tktritrue}{\Tcal_k^{\text{\rm trial},\text{\rm true}}}
\newcommand{\Titritrue}{\Tcal_i^{\text{\rm trial},\text{\rm true}}}
\newcommand{\Tjtritrue}{\Tcal_j^{\text{\rm trial},\text{\rm true}}}
\newcommand{\dktrue}{d_k^{\text{\rm true}}}

\newcommand{\kapmax}{\kappa_{\max}}

\def\Ek#1{E_{k,#1}}

In this appendix, we prove a formal version of Theorem~\ref{th.tau_zero_informal}, which is stated at the end of this appendix as Theorem~\ref{th.totalprob}.  Toward this end, we formalize the quantities generated by Algorithm~\ref{alg.sqp} as a stochastic process, namely,
\bequationNN
  \{(X_k,G_k,V_k,U_k,\Utrue_k,D_k,\Dtrue_k,Y_k,\Ytrue_k,\Tcal_k,\Tktritrue,\Xcal_k,\Zcal_k,\Xi_k,\Acal_k)\},
\eequationNN
where, for all $k \in \N{}$, we denote the primal iterate as $X_k$, the stochastic gradient estimate as $G_k$, the normal search direction as $V_k$, the tangential search direction as $U_k$, the ``true'' tangential search direction as $\Utrue_k$, the search direction as $D_k$, the ``true'' search direction as $\Dtrue_k$, the Lagrange multiplier estimate as $Y_k$, the ``true'' Lagrange multiplier estimate as $\Ytrue_k$, the merit parameter as $\Tcal_k$, the ``true'' trial merit parameter as $\Tktritrue$, the curvature parameter as $\Xcal_k$, the curvature threshold parameter as $\Zcal_k$, the ratio parameter as $\Xi_k$, and the step size as $\Acal_k$.  A realization of the $k$th element of this process are the quantities that have appeared throughout the paper, namely, $(x_k,g_k,v_k,u_k,\utrue_k,d_k,\dtrue_k,y_k,\ytrue_k,\tau_k,\tautruetrial_k,\chi_k,\zeta_k,\xi_k,\alpha_k)$.  Algorithm~\ref{alg.sqp}'s behavior is dictated entirely by the initial conditions (i.e., initial point and parameter values) as well as the sequence of stochastic gradient estimates; i.e., assuming for simplicity that the initial conditions are predetermined, a realization of $\{G_0,\dots,G_{k-1}\}$ determines the realizations of
\bequationNN
  \{X_j\}_{j=1}^k\ \ \text{and}\ \ \{V_k,U_k,\Utrue_k,D_k,\Dtrue_k,Y_k,\Ytrue_k,\Tcal_k,\Tktritrue,\Xcal_k,\Zcal_k,\Xi_k,\Acal_k)\}_{j=0}^{k-1}.
\eequationNN

In the process of proving our main result (Theorem \ref{th.totalprob}), we prove a set of lemmas about the behavior of the merit parameter sequence after a finite number of iterations. Specifically, we consider the behavior of Algorithm~\ref{alg.sqp} when terminated at $k = \kmax \in \N{}$.  With this consideration, we define a tree with a depth bounded by $\kmax$, which will be integral to our arguments in this section.  The proof of Theorem \ref{th.totalprob} ultimately considers the behavior of the algorithm as $\kmax \to \infty$.

Let $\Ical[\cdot]$ denote the indicator function of an event and, for all $k \in \kmax$, define the random variables
\bequationNN
  Q_k := \Ical[\Tktritrue < \Tcal_{k-1}]\ \ \text{and}\ \ W_k := \sum_{i=0}^{k-1} \Ical[\Tcal_i < \Tcal_{i-1}].
\eequationNN
Accordingly, for any realization of a run of Algorithm \ref{alg.sqp} and any $k \in \kmax$, the realization $(q_k,w_k)$ of $(Q_k,W_k)$ is determined at the beginning of iteration $k$.  The \emph{signature} of a realization up to iteration $k$ is $(q_0,\dots,q_k,w_0,\dots,w_k)$, which encodes all of the pertinent information regarding the behavior of the merit parameter sequence and these indicators up to the start of iteration $k$.

We use the set of all signatures to define a tree whereby each node contains a subset of all realizations of the algorithm. To construct the tree, we denote the root node by $N(q_0,w_0)$, where~$q_0$ is the indicator of the event $\tau_0^{\rm \text{trial},\text{\rm true}} < \tau_{-1}$, which is deterministic based on the initial conditions of the algorithm, and $w_0 = 0$. All realizations of the algorithm follow the same initialization, so~$q_0$ and~$w_0$ are in the signature of every realization. Next, we define a node $N(q_{[k]},w_{[k]})$ at depth $k \in [\kmax]$ (where the root node has a depth of $0$) in the tree as the set of all realizations of the algorithm for which the signature of the realization up to iteration $k$ is $(q_0,\dots,q_k,w_0,\dots,w_k)$. We define the edges in the tree by connecting nodes at adjacent levels, where node $N(q_{[k]},w_{[k]})$ is connected to node $N(q_{[k]},q_{k+1},w_{[k]},w_{k+1})$ for any $q_{k+1} \in \{0,1\}$ and $w_{k+1} \in \{w_k,w_k+1\}$.

Notationally, since the behavior of a realization of the algorithm up to iteration $k \in \N{}$ is completely determined by the initial conditions and the realization of $G_{[k-1]}$, we say that a realization described by $G_{[k-1]}$ belongs in node $N(q_{[k]},w_{[k]})$ by writing that $G_{[k-1]} \in N(q_{[k]},w_{[k]})$.  The initial condition, denoted for consistency as $G_{[-1]} \in N(q_0,w_0)$, occurs with probability one.  Based on the description above, the nodes of our tree satisfy the property that for any $k \geq 2$, the event $G_{[k-1]} \in N(q_{[k]},w_{[k]})$ occurs if and only if
\bequation\label{eq.defdef}
  Q_k = q_k,\ \ W_k = w_k,\ \ \text{and}\ \ G_{[k-2]} \in N(q_{[k-1]},w_{[k-1]}).
\eequation

Similar to Section \ref{sec.tau_big}, we consider the following event, under which the true trial merit parameter sequence $\{\Tktritrue\}$ is bounded below by a positive real number.

\begin{definition}
  For some $\tautruetrial_{\min} \in \R{}_{>0}$, the event $E_{\tau}(\tautruetrial_{\min})$ occurs if and only if $\Tktritrue \geq \tautruetrial_{\min}$ for all $k \in \N{}$ in all realizations of a run of the algorithm.
\end{definition}

When this event occurs, it follows that if $W_{\kbar} = \sum_{i=0}^{\kbar-1} \Ical[\Tcal_i < \Tcal_{i-1}] \geq \bar{s}(\tautruetrial_{\min})$ (see \eqref{eq.sbar}) in a run, then for all $k \geq \kbar$ in the run it follows that $\Tcal_{k-1} \leq \tautruetrial_{\min} \leq \Tktritrue$.

We are now prepared to state the assumption under which Theorem~\ref{th.totalprob} is proved.

\bassumption \label{assum:eventEapp}
  For some $\tautruetrial_{\min} \in \R{}_{>0}$, the event $E_{\tau}(\tautruetrial_{\min})$ occurs.  In addition, $\{\Xcal_k\}$, $\{\Zcal_k\}$, and $\{\Xi_k\}$ are constant.  Lastly, there exists $p_{\tau} \in (0,1]$ such that, for all $k \in [\kmax]$, one finds
  \bequation \label{eq:ptau}
    \P\left[\Tcal_k < \Tcal_{k-1} | E_{\tau}(\tautruetrial_{\min}), G_{[k-1]} \in N(q_{[k]},w_{[k]}), \tcalktritrue < \Tcal_{k-1} \right] \geq p_{\tau}.
  \eequation
\eassumption
Intuitively, equation~\eqref{eq:ptau} states that conditioned on $E_{\tau}(\tautruetrial_{\min})$, the behavior of the algorithm up to the beginning of iteration $k$, and $Q_k = 1$, the probability that the merit parameter is decreased in iteration $k$ is at least $p_\tau$.  For simplicity of notation, henceforth we define $E := E_{\tau}(\tautruetrial_{\min})$ and $\smax := \bar{s}(\tautruetrial_{\min})$.  We remark that if \eqref{eq.daniel_special_label} holds for all realizations of Algorithm \ref{alg.sqp}---as is the case when the distribution of the stochastic gradients satisfies a mild form of symmetry (see \cite[Example~3.17]{BeraCurtRobiZhou21} for a simple example)---then \eqref{eq:ptau} holds by Lemma \ref{lem:tautrialextra}.

Our main result, Theorem~\ref{th.totalprob}, essentially shows that the probability that $\tcalktritrue < \Tcal_k$ occurs infinitely often is zero.  Toward proving this result, we first prove a bound on the probability that $\tcalktritrue < \Tcal_k$ occurs at least $J$ times for any $J \in \N{}$ such that $J > \frac{\smax}{p_{\tau}} + 1$. Given such a $J$, we can define a number of important sets of nodes in the tree.  First, let
\bequationNN
  \Lgood := \left\{N(q_{[k]},w_{[k]}) : \left(\sum_{i=0}^k q_i < J \right) \land (w_k = s_{\max} \lor k = \kmax) \right\}
\eequationNN
be the set of nodes at which the sum of the elements of $q_{[k]}$ is sufficiently small (less than $J$) and either $w_k$ has reached $\smax$ or $k$ has reached $\kmax$. Second, let
\bequationNN
  \Lbad := \left\{N(p_{[k]},w_{[k]}) : \sum_{i=0}^{k} q_i \geq  J \right\}
\eequationNN
be the nodes in the complement of $\Lgood$ at which the sum of the elements of $q_{[k]}$ is at least $J$.  Going forward, we restrict attention to the tree defined by the root node and all paths from the root node that terminate at a node contained in $\Lgood \cup \Lbad$.  From this restriction and the definitions of $\Lgood$ and $\Lbad$, the tree has finite depth with the elements of $\Lgood \cup \Lbad$ being leaves.

Let us now define relationships between nodes.  The parent of a node is defined as
\bequationNN
  P(N(q_{[k]},w_{[k]})) = N(q_{[k-1]},w_{[k-1]}).
\eequationNN
On the other hand, the children of node $N(q_{[k]},w_{[k]})$ are defined as
\bequationNN
  C(N(q_{[k]},w_{[k]})) = \bcases \{N(q_{[k]},q_{k+1},w_{[k]},w_{k+1})\} & \text{if $N(q_{[k]},w_{[k]}) \not\in \Lgood \cup \Lbad$} \\ \emptyset & \text{otherwise.} \ecases
\eequationNN
Under these definitions, the paths down the tree terminate at nodes in $\Lgood \cup \Lbad$, reaffirming that these nodes are the leaves of the tree.  For convenience in the rest of our discussions, let $C(\emptyset) = \emptyset$.

We define the height of node $N(q_{[k]},w_{[k]})$ as the length of the longest path from $N(q_{[k]},w_{[k]})$ to a leaf node, i.e., the height is denoted as
\bequationNN
  h(N(q_{[k]},w_{[k]})) := \left(\min \{j \in \N{} \setminus \{0\} : C^j(N(q_{[k]},w_{[k]})) = \emptyset\}\right)-1,
\eequationNN
where $C^j(N(q_{[k]},w_{[k]}))$ is shorthand for applying the mapping $C(\cdot)$ consecutively $j$ times. From this definition, $h(N(q_{[k]},w_{[k]})) = 0$ for all $N(q_{[k]},w_{[k]}) \in \Lgood \cup \Lbad$.

Finally, let us define the event $\EbadkmaxJ$ as the event that for some $j \in [\kmax]$ one finds
\bequation\label{eq:EbadkmaxJ}
  \sum_{i=0}^{j} Q_i = \sum_{i=0}^{j} \Ical[\Titritrue < \Tcal_{i-1}] \geq J.
\eequation
Our first goal in this section is to find a bound on the probability of this event occurring. We will then utilize this bound to prove Theorem \ref{th.totalprob}. As a first step towards bounding the probability of $\EbadkmaxJ$, we prove the following result about the leaf nodes of the tree.

\blemma \label{lem:leafnodes}
  For any $k \in [\kmax]$, $J \in \N{}$, and $(q_{[k]},w_{[k]})$ with $N(q_{[k]},w_{[k]}) \in \Lgood$, one finds
  \bequationNN
    \P[G_{[k-1]} \in N(q_{[k]},w_{[k]}) \land \EbadkmaxJ | E] = 0.
  \eequationNN
  On the other hand, for any $k \in [\kmax]$, $J \in \N{}$, and $(q_{[k]},w_{[k]})$ with $N(q_{[k]},w_{[k]}) \in \Lbad$, one finds
  \begin{align*}
    &\P[G_{[k-1]} \in N(q_{[k]},w_{[k]}) \land \EbadkmaxJ | E] \\
    &\leq \prod_{i=1}^k (\P[Q_i = q_i | E, W_i = w_i, G_{[i-2]} \in N(q_{[i-1]},w_{[i-1]})] \cdot \P[W_i = w_i | E, G_{[i-2]} \in N(q_{[i-1]},w_{[i-1]})]).
  \end{align*}
\elemma
\proof{Proof.}
  Consider arbitrary $k \in [\kmax]$ and $J \in \N{}$ as well as an arbitrary pair $(q_{[k]},w_{[k]})$ such that $N(q_{[k]},w_{[k]}) \in \Lgood$.  By the definition of $\Lgood$, it follows that $\sum_{i=0}^k q_i < J$.  Then, by \eqref{eq.defdef},
  \begin{align*}
    &\P\left[\sum_{i=0}^{k} Q_i \geq J \Big| E, G_{[k-1]} \in N(q_{[k]},w_{[k]})\right] = \P\left[\sum_{i=0}^k q_i \geq J \Big| E, G_{[k-1]} \in N(q_{[k]},w_{[k]})\right] = 0.
  \end{align*}
  Therefore, for any $j \in \{1,\dots,k\}$, one finds from conditional probability that
  \begin{align*}
    \P\left[G_{[j-1]} \in N(q_{[j]},w_{[j]}) \land \text{\eqref{eq:EbadkmaxJ} holds} | E\right] &=\P\left[\text{\eqref{eq:EbadkmaxJ} holds} \Big| E, G_{[j-1]} \in N(q_{[j]},w_{[j]})\right] \\
    &\quad \cdot \P\left[G_{[j-1]} \in N(q_{[j]},w_{[j]}) | E\right] = 0.
  \end{align*}
  In addition, \eqref{eq:EbadkmaxJ} cannot hold for $j=0$ since $\Ical[\tautruetrial_0 < \tau_{-1}] = q_0 < J$ by the definition of $\Lgood$. Hence, along with the conclusion above, it follows that $\EbadkmaxJ$ does not occur in any realization whose signature up to iteration $j \in \{1,\dots,k\}$ falls into a node along any path from the root to $N(q_{[k]},w_{[k]})$.  Now, by the definition of $\Lgood$, at least one of $w_k = \smax$ or $k = \kmax$ holds.  Let us consider each case in turn.  If $k = \kmax$, then it follows by the preceding arguments that 
  \begin{equation*}
    \P\left[\sum_{i=0}^{\kmax} Q_i < J \Big| E, G_{[k-1]} \in N(p_{[k]},w_{[k]})\right] = 1.
  \end{equation*}
  Otherwise, if $w_k = \smax$, then it follows by the definition of $\smax$ that $\Tcal_{k-1} \leq \tautruetrial_{\min}$ so that $Q_i = \Ical[\Titritrue < \Tcal_{i-1}] = 0$ holds for all $i \in \{k,\dots,\kmax\}$, and therefore the equation above again follows.  Overall, it follows that $\P[G_{[k-1]} \in N(q_{[k]},w_{[k]}) \land \EbadkmaxJ | E] = 0$, as desired.
   
  Now consider arbitrary $k \in [\kmax]$ and $J \in \N{}$ as well as an arbitrary pair $(q_{[k]},w_{[k]})$ with $N(q_{[k]},w_{[k]}) \in \Lbad$.  One finds that
  \begin{align*}
    &\ \P[G_{[k-1]} \in N(q_{[k]},w_{[k]}) \land \EbadkmaxJ | E] \\
    \leq&\ \P[G_{[k-1]} \in N(q_{[k]},w_{[k]}) | E] = \P[\text{\eqref{eq.defdef} holds} | E] \\
    =&\ \P[Q_k = q_k | E, W_k = w_k, G_{[k-2]} \in N(q_{[k-1]},w_{[k-1]})] \cdot \P[W_k = w_k \land G_{[k-2]} \in N(p_{[k-1]},w_{[k-1]}) | E] \\
    =&\ \P[Q_k = q_k | E, W_k = w_k, G_{[k-2]} \in N(q_{[k-1]},w_{[k-1]})] \\
    &\ \cdot \P[W_k = w_k | E, G_{[k-2]} \in N(p_{[k-1]},w_{[k-1]})]
    \cdot \P[G_{[k-2]} \in N(p_{[k-1]},w_{[k-1]}) | E] \\
    =&\ \P[G_{-1} \in N(q_0,w_0)] \\
    &\ \cdot \prod_{i=1}^k (\P[Q_i = q_i | E, W_i = w_i, G_{[i-2]} \in N(q_{[i-1]},w_{[i-1]})] \cdot \P[W_i = w_i | E, G_{[i-2]} \in N(q_{[i-1]},w_{[i-1]})]),
  \end{align*}
  which, since $\P[G_{[-1]} \in N(q_0,w_0)] = 1$, proves the remainder of the result. \Halmos
\endproof

Next, we show that the probability of the occurrence of $\EbadkmaxJ$ at any node in the tree can be bounded in terms of the probability of the sum of a set of independent Bernoulli random variables being less than a threshold defined by $\smax$.

\blemma \label{lem:appinduction}
  For any $k \in [\kmax]$, $J \in \N{}$, and $(q_{[k]},w_{[k]})$ with $N(q_{[k]},w_{[k]}) \not\in \Lgood$, let
  \bequation \label{eq:psidef}
    \psi_J(q_{[k]}) = J - 1 - \sum_{i=0}^k q_i.
  \eequation
  One finds that
  \begin{align}
    &\P[G_{[k-1]} \in N(q_{[k]},w_{[k]}) \land \EbadkmaxJ | E] \nonumber \\
    &\leq \prod_{i=1}^k (\P[Q_i = q_i | E, W_i = w_i, G_{[i-2]} \in N(q_{[i-1]},w_{[i-1]})] \nonumber \\
    &\hspace{35pt} \cdot \P[W_i = w_i | E, G_{[i-2]} \in N(q_{[i-1]},w_{[i-1]})]) \cdot \P\left[\sum_{j=1}^{\psi_J(q_{[k-1]})} Z_j \leq \smax - w_k\right], \label{eq:inductionbound}
  \end{align}
  where $\{Z_j\}$ are independent Bernoulli random variables with $\P[Z_j =  1] = p_{\tau}$ for all $j \in \N{}$.
\elemma
\proof{Proof.}
  Consider any $(q_{[k]},w_{[k]})$ with $h(N(q_{[k]},w_{[k]}))=0$. Since $N(q_{[k]},w_{[k]}) \not\in \Lgood$, it follows that $N(q_{[k]},w_{[k]}) \in \Lbad$. Then, by the definition of $\Lbad$, it follows that $\sum_{i=0}^{k} q_i \geq J$.  In addition, since $C(N(q_{[k]},w_{[k]})) = \emptyset$ for any node in $\Lbad$, it follows that $P(N(q_{[k]},w_{[k]})) \not\in \Lbad$, which implies that $\sum_{i=0}^{k-1} q_i < J$.  Thus, $\sum_{i=0}^{k} q_i = J$ and $\sum_{i=0}^{k-1} q_i = J - 1$, which implies that $\psi_J(q_{[k-1]}) = 0$. Therefore, overall, the result holds for any $(q_{[k]},w_{[k]})$ with $h(N(q_{[k]},w_{[k]}))=0$ by Lemma \ref{lem:leafnodes}.
    
  We prove the rest of the result by induction on the height of the node. We note that the base case, i.e., when $h(N(q_{[k]},w_{[k]}))=0$, holds by the above argument. Now, assume that \eqref{eq:inductionbound} holds for any $(q_{[k]},w_{[k]})$ with $N(q_{[k]},w_{[k]}) \not \in \Lgood$ such that $h(N(q_{[k]},w_{[k]})) \leq \hat{h}$.  Consider arbitrary $(q_{[k]},w_{[k]})$ such that $N(q_{[k]},w_{[k]}) \not \in \Lgood$ and $h(N(q_{[k]},w_{[k]})) = \hat{h}+1$. By the definition of $C$, one finds
  \begin{align*}
    &\P[G_{[k-1]} \in N(q_{[k]},w_{[k]}) \land \EbadkmaxJ | E] \\
    &= \sum_{\{(q_{k+1},w_{k+1}) : N(q_{[k+1]},w_{[k+1]}) \in C(N(q_{[k]},w_{[k]})) \}} \P[G_{[k]} \in N(q_{[k+1]},w_{[k+1]}) \land \EbadkmaxJ | E].
  \end{align*}
  Then, by the definition of $q_{[k]}$ and $w_{[k]}$, we can enumerate the children of $N(q_{[k]},w_{[k]})$ as
  \begin{align*}
    &\P[G_{[k-1]} \in N(q_{[k]},w_{[k]}) \land \EbadkmaxJ | E] \\
    &= \P[G_{[k]} \in N(q_{[k]},0,w_{[k]},w_k) \land \EbadkmaxJ | E] + \P[G_{[k]} \in N(q_{[k]},0,w_{[k]},w_k+1) \land \EbadkmaxJ | E] \\
    &\quad+ \P[G_{[k]} \in N(q_{[k]},1,w_{[k]},w_k) \land \EbadkmaxJ | E] + \P[G_{[k]} \in N(q_{[k]},1,w_{[k]},w_k+1) \land \EbadkmaxJ | E].
  \end{align*}
  Now, noting that all children of $N(q_{[k]},w_{[k]})$ have a height that is at most $\hat{h}$, we apply the induction hypothesis four times to obtain
  \begin{align*}
    &\P[G_{[k-1]} \in N(q_{[k]},w_{[k]}) \land \EbadkmaxJ | E] \\
    &\leq \Bigg(\P[Q_{k+1} = 0 | E, W_{k+1} = w_k, G_{[k-1]} \in N(q_{[k]},w_{[k]})] \\
    &\hspace{28pt} \cdot \P[W_{k+1} = w_k | E, G_{[k-1]} \in N(q_{[k]},w_{[k]})] \cdot \P\left[\sum_{j=1}^{\psi_J(q_{[k]})} Z_{j,1} \leq \smax - w_k\right] \\
    &\hspace{28pt}+ \P[Q_{k+1} = 0 | E, W_{k+1} = w_k+1, G_{[k-1]} \in N(q_{[k]},w_{[k]})] \\
    &\hspace{28pt} \cdot \P[W_{k+1} = w_k + 1 | E, G_{[k-1]} \in N(q_{[k]},w_{[k]})] \cdot \P\left[\sum_{j=1}^{\psi_J(q_{[k]})} Z_{j,2} \leq \smax - w_k - 1\right] \\
    &\hspace{28pt}+ \P[Q_{k+1} = 1 | E, W_{k+1} = w_k, G_{[k-1]} \in N(q_{[k]},w_{[k]})] \\
    &\hspace{28pt} \cdot \P[W_{k+1} = w_k | E, G_{[k-1]} \in N(q_{[k]},w_{[k]})] \cdot \P\left[\sum_{j=1}^{\psi_J(q_{[k]})} Z_{j,3} \leq \smax - w_k\right] \\
    &\hspace{28pt}+ \P[Q_{k+1} = 1 | E, W_{k+1} = w_k+1, G_{[k-1]} \in N(q_{[k]},w_{[k]})] \\
    &\hspace{28pt} \cdot \P[W_{k+1} = w_k + 1 | E, G_{[k-1]} \in N(q_{[k]},w_{[k]})] \cdot \P\left[\sum_{j=1}^{\psi_J(q_{[k]})} Z_{j,4} \leq \smax - w_k - 1\right]\Bigg) \\
    &\hspace{28pt} \cdot \prod_{i=1}^k (\P[Q_i = q_i | E, W_i = w_i, G_{[i-2]} \in N(q_{[i-1]},w_{[i-1]})] \cdot \P[W_i = w_i | E, G_{[i-2]} \in N(q_{[i-1]},w_{[i-1]})])
  \end{align*}
  where $Z_{j,p}$ for all $p \in \{1,\dots,4\}$ and $j \in \{1, \dots, \psi(q_{[k]})\}$ are four sets of independent Bernoulli random variables with $\P[Z_{j,p}=1]=p_{\tau}$.  Now, by the definitions of $Z_{j,1}$, $Z_{j,2}$, $Z_{j,3}$, and $Z_{j,4}$, it follows that
    \bequationNN
        \P\left[\sum_{j=1}^{\psi_J(q_{[k]})} Z_{j,1} \leq \smax - w_k\right] = \P\left[\sum_{j=1}^{\psi_J(q_{[k]})} Z_{j,3} \leq \smax - w_k\right],
    \eequationNN
    and
    \bequationNN
        \P\left[\sum_{j=1}^{\psi_J(q_{[k]})} Z_{j,2} \leq \smax - w_k-1\right] = \P\left[\sum_{j=1}^{\psi_J(q_{[k]})} Z_{j,4} \leq \smax - w_k-1\right].
    \eequationNN
    Therefore, it follows that
    \begin{align*}
        &\P[G_{[k-1]} \in N(q_{[k]},w_{[k]}) \land \EbadkmaxJ | E] \\
        &\leq \Bigg((\P[Q_{k+1} = 0 | E, W_{k+1} = w_k, G_{[k-1]} \in N(q_{[k]},w_{[k]})] \\
        &\hspace{28pt} + \P[Q_{k+1} = 1 | E, W_{k+1} = w_k, G_{[k-1]} \in N(q_{[k]},w_{[k]})]) \\
        &\hspace{28pt} \cdot \P[W_{k+1} = w_k | E, G_{[k-1]} \in N(q_{[k]},w_{[k]})] \cdot \P\left[\sum_{j=1}^{\psi_J(q_{[k]})} Z_{j,1} \leq \smax - w_k\right] \\
        &\hspace{28pt} + (\P[Q_{k+1} = 0 | E, W_{k+1} = w_k+1, G_{[k-1]} \in N(q_{[k]},w_{[k]})] \\
        &\hspace{28pt}+ \P[Q_{k+1} = 1 | E, W_{k+1} = w_k+1, G_{[k-1]} \in N(q_{[k]},w_{[k]})]) \\
        &\hspace{28pt} \cdot \P[W_{k+1} = w_k + 1 | E, G_{[k-1]} \in N(q_{[k]},w_{[k]})] \cdot \P\left[\sum_{j=1}^{\psi_J(q_{[k]})} Z_{j,2} \leq \smax - w_k - 1\right] \Bigg) \\
        &\hspace{28pt} \cdot \prod_{i=1}^k (\P[Q_i = q_i | E, W_i = w_i, G_{[i-2]} \in N(q_{[i-1]},w_{[i-1]})] \cdot \P[W_i = w_i | E, G_{[i-2]} \in N(q_{[i-1]},w_{[i-1]})]).
    \end{align*}
    Now, by the law of total probability, it follows that
    \begin{align*}
        1 &= \P[Q_{k+1} = 0 | E, W_{k+1} = w_k, G_{[k-1]} \in N(q_{[k]},w_{[k]})] \\
        &\quad + \P[Q_{k+1} = 1 | E, W_{k+1} = w_k, G_{[k-1]} \in N(q_{[k]},w_{[k]})],
    \end{align*}
    and
    \begin{align*}
        1 &= \P[Q_{k+1} = 0 | E, W_{k+1} = w_k+1, G_{[k-1]} \in N(q_{[k]},w_{[k]})] \\
        &\quad + \P[Q_{k+1} = 1 | E, W_{k+1} = w_k+1, G_{[k-1]} \in N(q_{[k]},w_{[k]})].
    \end{align*}
    Thus,
    \begin{align}
        &\P[G_{[k-1]} \in N(q_{[k]},w_{[k]}) \land \EbadkmaxJ | E] \nonumber \\
        &\leq \Bigg(\P[W_{k+1} = w_k | E, G_{[k-1]} \in N(q_{[k]},w_{[k]})] \cdot \P\left[\sum_{j=1}^{\psi_J(q_{[k]})} Z_{j,1} \leq \smax - w_k\right] \nonumber \\
        &\hspace{28pt}+ \P[W_{k+1} = w_k + 1 | E, G_{[k-1]} \in N(q_{[k]},w_{[k]})] \cdot \P\left[\sum_{j=1}^{\psi_J(q_{[k]})} Z_{j,2} \leq \smax - w_k - 1\right] \Bigg) \nonumber \\
        &\hspace{28pt} \cdot \prod_{i=1}^k (\P[Q_i = q_i | E, W_i = w_i, G_{[i-2]} \in N(q_{[i-1]},w_{[i-1]})] \cdot \P[W_i = w_i | E, G_{[i-2]} \in N(q_{[i-1]},w_{[i-1]})]). \label{eq:appmainlemeq1}
    \end{align}
    
    We proceed by considering two cases.  First, suppose $q_k = 1$. By Assumption \ref{assum:eventEapp}, it follows that
    \begin{align*}
        &\P[W_{k+1} = w_k + 1 | E, G_{[k-1]} \in N(q_{[k]},w_{[k]})] \\
        &= \P[\Tcal_k < \Tcal_{k-1} | E, G_{[k-1]} \in N(q_{[k]},w_{[k]}), \tcalktritrue < \Tcal_{k-1}] \geq p_{\tau}.
    \end{align*}
    Additionally, using the law of total probability, we have
    \bequationNN
      1 = \P[W_{k+1} = w_k | E, G_{[k-1]} \in N(q_{[k]},w_{[k]})] + \P[W_{k+1} = w_k + 1 | E, G_{[k-1]} \in N(q_{[k]},w_{[k]})].
    \eequationNN
    Therefore, it follows that
    \begin{align}
        &\P[G_{[k-1]} \in N(q_{[k]},w_{[k]}) \land \EbadkmaxJ | E] \nonumber \\
        &\leq \underset{p \in [p_{\tau},1]}{\max} \Bigg((1-p) \P\left[\sum_{j=1}^{\psi_J(q_{[k]})} Z_{j,1} \leq \smax - w_k\right] + p \P\left[\sum_{j=1}^{\psi_J(q_{[k]})} Z_{j,2} \leq \smax - w_k - 1\right] \Bigg) \nonumber \\
        &\hspace{28pt} \cdot \prod_{i=1}^k (\P[Q_i = q_i | E, W_i = w_i, G_{[i-2]} \in N(q_{[i-1]},w_{[i-1]})] \cdot \P[W_i = w_i | E, G_{[i-2]} \in N(q_{[i-1]},w_{[i-1]})]). \label{eq:maxptau}
    \end{align}
    In addition, by the definition of $Z_{j,1}$ and $Z_{j,2}$, one finds that
    \begin{equation*}
        \P\left[\sum_{j=1}^{\psi_J(q_{[k]})} Z_{j,2} \leq \smax - w_k - 1\right] \leq \P\left[\sum_{j=1}^{\psi_J(q_{[k]})} Z_{j,1} \leq \smax - w_k\right].
    \end{equation*}
    Therefore, it follows that the max in \eqref{eq:maxptau} is given by $p = p_{\tau}$.
    Thus,
    \begin{align*}
        &\P[G_{[k-1]} \in N(q_{[k]},w_{[k]}) \land \EbadkmaxJ | E] \\
        &\leq \Bigg((1-p_{\tau}) \P\left[\sum_{j=1}^{\psi_J(q_{[k]})} Z_{j,1} \leq \smax - w_k\right] + p_{\tau} \P\left[\sum_{j=1}^{\psi_J(q_{[k]})} Z_{j,1} \leq \smax - w_k - 1\right] \Bigg) \\
        &\hspace{28pt} \cdot \prod_{i=1}^k (\P[Q_i = q_i | E, W_i = w_i, G_{[i-2]} \in N(q_{[i-1]},w_{[i-1]})] \cdot \P[W_i = w_i | E, G_{[i-2]} \in N(q_{[i-1]},w_{[i-1]})]),
    \end{align*}
    where, by the definitions of $Z_{j,1}$ and $Z_{j,2}$, we have used the fact that
    \bequationNN
        \P\left[\sum_{j=1}^{\psi_J(q_{[k]})} Z_{j,1} \leq \smax - w_k - 1\right] = \P\left[\sum_{j=1}^{\psi_J(q_{[k]})} Z_{j,2} \leq \smax - w_k - 1\right].
    \eequationNN
    Now, for all $j \in \{1,\dots,\psi_J(q_{[k]})\}$, define $Z_j = Z_{j,1}$ and let $Z_{\psi_J(q_{[k]})+1}$ be a Bernoulli random variable with $\P[Z_{\psi_J(q_{[k]})+1} = 1] = p_{\tau}$. Then, it follows that
    \begin{align*}
        &\P[G_{[k-1]} \in N(q_{[k]},w_{[k]}) \land \Ebadkmax | E] \\
        &\leq \P\left[\sum_{j=1}^{\psi_J(q_{[k]})+1} Z_j \leq \smax - w_k\right] \\
        &\qquad \cdot \prod_{i=1}^k (\P[Q_i = q_i | E, W_i = w_i, G_{[i-2]} \in N(q_{[i-1]},w_{[i-1]})] \cdot \P[W_i = w_i | E, G_{[i-2]} \in N(q_{[i-1]},w_{[i-1]})]).
    \end{align*}
    This proves the result in this case by noting that $q_k=1$ implies
    \bequationNN
        \psi_J(q_{[k]})+1 =  J - 1 - \sum_{i=0}^k q_i + 1 = J - 1 - \sum_{i=0}^{k-1} q_i = \psi_J(q_{[k-1]}).
    \eequationNN

    Next, consider the case where $q_k = 0$. Recalling that
    \bequationNN
      1 = \P[W_{k+1} = w_k | E, G_{[k-1]} \in N(q_{[k]},w_{[k]})] + \P[W_{k+1} = w_k + 1 | E, G_{[k-1]} \in N(q_{[k]},w_{[k]})],
    \eequationNN
    it follows from \eqref{eq:appmainlemeq1} that
    \begin{align}
        &\P[G_{[k-1]} \in N(q_{[k]},w_{[k]}) \land \EbadkmaxJ | E] \nonumber \\
        &\leq \underset{p\in[0,1]}{\max} \Bigg((1-p) \P\left[\sum_{j=1}^{\psi_J(q_{[k]})} Z_{j,1} \leq \smax - w_k\right] + p \P\left[\sum_{j=1}^{\psi_J(q_{[k]})} Z_{j,2} \leq \smax - w_k - 1\right] \Bigg) \nonumber \\
        &\hspace{28pt} \cdot \prod_{i=1}^k (\P[Q_i = q_i | E, W_i = w_i, G_{[i-2]} \in N(q_{[i-1]},w_{[i-1]})] \cdot \P[W_i = w_i | E, G_{[i-2]} \in N(q_{[i-1]},w_{[i-1]})]). \label{eq:newlabel}
    \end{align}
    Similar to before, noting that
    \begin{equation*}
        \P\left[\sum_{j=1}^{\psi_J(q_{[k]})} Z_{j,2} \leq \smax - w_k - 1\right] \leq \P\left[\sum_{j=1}^{\psi_J(q_{[k]})} Z_{j,1} \leq \smax - w_k\right],
    \end{equation*}
    it follows that the max in \eqref{eq:newlabel} is given by $p=0$, so, with $Z_j = Z_{j,1}$ for all $j \in \{1,\dots,\psi_J(q_{[k]})\}$,
    \begin{align*}
        &\P[G_{[k-1]} \in N(q_{[k]},w_{[k]}) \land \EbadkmaxJ | E] \\
        &\leq \P\left[\sum_{j=1}^{\psi_J(q_{[k]})} Z_j \leq \smax - w_k\right] \\
        &\quad \cdot \prod_{i=1}^k (\P[Q_i = q_i | E, W_i = w_i, G_{[i-2]} \in N(q_{[i-1]},w_{[i-1]})] \cdot \P[W_i = w_i | E, G_{[i-2]} \in N(q_{[i-1]},w_{[i-1]})]).
    \end{align*}
    The result follows from this inequality and the fact that $\psi_J(q_{[k]}) = \psi_J(q_{[k-1]})$ since $q_k=0$. \Halmos
\endproof

We now apply Lemma \ref{lem:appinduction} to obtain a high probability bound.

\begin{lemma} \label{lem:highprobapp}
  For any $J > \frac{\smax}{p_{\tau}}+1$, one finds that
  \bequation \label{eq:highprobbound}
    \P\left[\sum_{k=0}^{\kmax} \Ical[\Tktritrue < \Tcal_k] \geq J\right] \leq e^{- \frac{p_{\tau}(J-1)}{2} \left(1-\frac{\smax}{p_{\tau}(J-1)}\right)^2}.
  \eequation
\end{lemma}

\proof{Proof.}
    Recalling that the initial condition for the tree, $G_{[-1]} \in N(q_0,w_0)$, occurs with probability one, by Lemma \ref{lem:appinduction}, it follows that there exist $J-1$ independent Bernoulli random variables $Z_j$ with $\P[Z_j=1]=p_{\tau}$ for all $j \in \{1,\dots,J-1\}$ such that
    \bequationNN
        \P[\EbadkmaxJ | E] = \P[G_{[-1]} \in N(q_0,w_0) \land \Ebadkmax | E] \leq \P\left[\sum_{j=1}^{J-1} Z_j \leq \smax\right]. 
    \eequationNN
    Let
    \bequationNN
        \mu := \sum_{j=1}^{J-1} \P[Z_j = 1] = p_{\tau} (J-1) \ \text{ and } \ \rho := 1 - \smax/\mu.
    \eequationNN
    Noting that $\rho \in (0,1)$ by the definition of $J$, by the multiplicative form of Chernoff's bound,
    \bequationNN
        \P\left[\sum_{j=1}^{J-1} Z_j \leq \smax \right] \leq e^{-\thalf \mu \rho^2} = e^{- \thalf \mu (1-\smax/\mu)^2} = e^{- \frac{p_{\tau}(J-1)}{2} \left(1-\frac{\smax}{p_{\tau}(J-1)}\right)^2}.
    \eequationNN
    For all $k \in [\kmax]$, we have $\Tcal_k \leq \Tcal_{k-1}$. Thus, by the definition of $\EbadkmaxJ$, it follows that
    \begin{align*}
        \P\left[\sum_{k=0}^{\kmax} \Ical[\Tktritrue < \Tcal_{k}] \geq J \Bigg| E \right] &\leq \P\left[\sum_{k=0}^{\kmax} Q_k \geq J \Bigg| E \right] \leq \P[\EbadkmaxJ | E] \leq e^{- \frac{p_{\tau}(J-1)}{2} \left(1-\frac{\smax}{p_{\tau}(J-1)}\right)^2},
    \end{align*}
    as desired.
    \Halmos
\endproof

We now prove the main result of this appendix.

\btheorem\label{th.totalprob}
  Under Assumption \ref{assum:eventEapp}, it follows that
  \bequation
    \P\left[\sum_{k=0}^{\infty} \Ical[\Tktritrue < \Tcal_k] < \infty \Bigg| E\right] = 1.
  \eequation
\etheorem
\proof{Proof.}
    By Lemma \ref{lem:highprobapp}, for any $\kmax \in \N{} \setminus \{0\}$ and $J > \frac{\smax}{p_{\tau}}+1$, it follows that
    \bequationNN
        \P\left[\sum_{k=0}^{\kmax} \Ical[\Tktritrue < \Tcal_k] \geq J \Bigg| E \right] \leq e^{- \frac{p_{\tau}(J-1)}{2} \left(1-\frac{\smax}{p_{\tau}(J-1)}\right)^2}.
    \eequationNN
    Let $A_{\kmax}$ denote the event that
    \bequationNN
        \sum_{k=0}^{\kmax} \Ical[\Tktritrue < \Tcal_k] \geq J.
    \eequationNN
    It follows from this definition that $A_{\kmax} \subseteq A_{\kmax+1}$ for any $\kmax \in \N{} \setminus \{0\}$. Therefore, by the properties of an increasing sequence of events (see, for example \cite[Section 1.5]{DStirzaker_2003}), it follows that
    \begin{align*}
        \P\left[\sum_{k=0}^{\infty} \Ical[\Tktritrue < \Tcal_k] \geq J \Bigg| E\right]
        &= \P\left[\lim_{\kmax \rightarrow \infty} \sum_{k=0}^{\kmax} \Ical[\Tktritrue < \Tcal_k] \geq J \Bigg| E\right] \\
        &= \lim_{\kmax \rightarrow \infty} \P\left[\sum_{k=0}^{\kmax} \Ical[\Tktritrue < \Tcal_k] \geq J \Bigg| E\right] \leq e^{- \frac{p_{\tau}(J-1)}{2} \left(1-\frac{\smax}{p_{\tau}(J-1)}\right)^2}.
    \end{align*}
    Next, let $A_J$ denote the event that
    \bequationNN
        \sum_{k=0}^{\infty} \Ical[\Tktritrue < \Tcal_k] < J.
    \eequationNN
    From the definition of $A_J$, it follows that $A_J \subseteq A_{J+1}$ for any $J > \frac{\smax}{p_{\tau}}+1$.  Thus, as above,
    \begin{align*}
        \P\left[\sum_{k=0}^{\infty} \Ical[\Tktritrue < \Tcal_k] < \infty \Bigg| E\right]
        &= \P\left[\lim_{J \rightarrow \infty} \sum_{k=0}^{\infty} \Ical[\Tktritrue < \Tcal_k] < J \Bigg| E\right] \\
        &= \lim_{J \rightarrow \infty} \P\left[\sum_{k=0}^{\infty} \Ical[\Tktritrue < \Tcal_k] < J \Bigg| E\right] \\
        &\geq \lim_{J \rightarrow \infty} 1 - e^{- \frac{p_{\tau}(J-1)}{2} \left(1-\frac{\smax}{p_{\tau}(J-1)}\right)^2} = 1,
    \end{align*}
    which is the desired conclusion.
    \Halmos
\endproof

\section*{Acknowledgments.}

This material is based upon work supported by the U.S.~National Science Foundation's Division of Computing and Communication Foundations under award number CCF-1740796 and by the Office of Naval Research under award number N00014-21-1-2532.

%**************
% Bibliography
%**************
\bibliographystyle{plain}
\bibliography{references}

\begin{thebibliography}{10}

\bibitem{BeraCurtRobiZhou21}
Albert~S. Berahas, Frank~E. Curtis, Daniel~P. Robinson, and Baoyu Zhou.
\newblock {Sequential Quadratic Optimization for Nonlinear Equality Constrained
  Stochastic Optimization}.
\newblock {\em {SIAM Journal on Optimization}}, 31(2):1352--1379, 2021.

\bibitem{BottCurtNoce18}
L{\'e}on Bottou, Frank~E Curtis, and Jorge Nocedal.
\newblock Optimization methods for large-scale machine learning.
\newblock {\em SIAM Review}, 60(2):223--311, 2018.

\bibitem{chang2011libsvm}
Chih-Chung Chang and Chih-Jen Lin.
\newblock {LIBSVM}: a library for support vector machines.
\newblock {\em ACM Transactions on Intelligent Systems and Technology (TIST)},
  2(3):1--27, 2011.

\bibitem{ChenTungVeduMori18}
Changan Chen, Frederick Tung, Naveen Vedula, and Greg Mori.
\newblock {Constraint-aware deep neural network compression}.
\newblock In {\em Proceedings of the European Conference on Computer Vision
  (ECVC)}, pages 400--415, 2018.

\bibitem{CurtNoceWaec09}
Frank~E. Curtis, Jorge Nocedal, and Andreas W\"{a}chter.
\newblock A matrix-free algorithm for equality constrained optimization
  problems with rank deficient {J}acobians.
\newblock {\em {SIAM Journal on Optimization}}, 20(3):1224--1249, 2009.

\bibitem{DaviDrusKaka20}
Damek Davis, Dmitriy Drusvyatskiy, and Sham Kakade.
\newblock Stochastic subgradient method converges on tame functions.
\newblock {\em {Foundations of Computational Mathematics}}, 20:119--154, 2020.

\bibitem{GoulLuciRomaToin99}
N.~I.~M. Gould, S.~Lucidi, M.~Roma, and Ph.~L. Toint.
\newblock {Solving the trust-region subproblem using the Lanczos method}.
\newblock {\em SIAM Journal on Optimization}, 9(2):504--525, 1999.

\bibitem{GoulOrbaToin15}
Nicolas I.~M. Gould, Dominique Orban, and Philippe~L. Toint.
\newblock {CUTEst}: a constrained and unconstrained testing environment with
  safe threads for mathematical optimization.
\newblock {\em Computational Optimization and Applications}, 60:545--557, 2015.

\bibitem{Han77}
S.~P. Han.
\newblock A globally convergent method for nonlinear programming.
\newblock {\em Journal of Optimization Theory and Applications},
  22(3):297--309, 1977.

\bibitem{HanMang79}
S.~P. Han and O.~L. Mangasarian.
\newblock Exact penalty functions in nonlinear programming.
\newblock {\em Mathematical Programming}, 17:251--269, 1979.

\bibitem{hazan2016variance}
Elad Hazan and Haipeng Luo.
\newblock Variance-reduced and projection-free stochastic optimization.
\newblock In {\em Proceedings of International Conference on Machine Learning
  (ICML)}, pages 1263--1271, 2016.

\bibitem{KumaSoumMhamHara18}
Soumava Kumar~Roy, Zakaria Mhammedi, and Mehrtash Harandi.
\newblock {Geometry aware constrained optimization techniques for deep
  learning}.
\newblock In {\em Proceedings of Computer Vision and Pattern Recognition
  (CVPR)}, pages 4460--4469, 2018.

\bibitem{locatello2019stochastic}
Francesco Locatello, Alp Yurtsever, Olivier Fercoq, and Volkan Cevher.
\newblock Stochastic {F}rank-{W}olfe for composite convex minimization.
\newblock In {\em Proceedings of Neural Information Processing Systems
  (NeurIPS)}, pages 14269--14279, 2019.

\bibitem{lu2020generalized}
Haihao Lu and Robert~M Freund.
\newblock Generalized stochastic {F}rank-{W}olfe algorithm with stochastic
  “substitute” gradient for structured convex opt.
\newblock {\em Math. Prog}, pages 1--33, 2020.

\bibitem{NaAnitKola21}
Sen Na, Mihai Anitescu, and Mladen Kolar.
\newblock An adaptive stochastic sequential quadratic programming with
  differentiable exact augmented {L}agrangians.
\newblock {\em arXiv preprint arXiv:2102.05320}, 2021.

\bibitem{NandPathAbhiSing19}
Yatin Nandwani, Abhishek Pathak, and Parag Singla.
\newblock {A primal-dual formulation for deep learning with constraints}.
\newblock In {\em Proceedings of Neural Information Processing Systems
  (NeurIPS)}, pages 12157--12168, 2019.

\bibitem{NoceWrig06}
Jorge Nocedal and Stephen Wright.
\newblock {\em Numerical optimization}.
\newblock Springer Series in Operations Research and Financial Engineering.
  Springer-Verlag New York, 2006.

\bibitem{Omoj89}
E.~O. Omojokun.
\newblock {\em {Trust Region Algorithms for Optimization with Nonlinear
  Equality and Inequality Constraints}}.
\newblock PhD thesis, University of Colorado, Boulder, CO, USA, 1989.

\bibitem{Powe78}
M.~J.~D. Powell.
\newblock A fast algorithm for nonlinearly constrained optimization
  calculations.
\newblock In {\em Numerical Analysis}, Lecture Notes in Mathematics, pages
  144--157. Springer, Berlin, 1978.

\bibitem{RaviDinhLokhSing19}
Sathya~N Ravi, Tuan Dinh, Vishnu~Suresh Lokhande, and Vikas Singh.
\newblock Explicitly imposing constraints in deep networks via conditional
  gradients gives improved generalization and faster convergence.
\newblock In {\em Proceedings of the AAAI Conference on Artificial
  Intelligence}, volume~33, pages 4772--4779, 2019.

\bibitem{reddi2016stochastic}
Sashank~J Reddi, Suvrit Sra, Barnab{\'a}s P{\'o}czos, and Alex Smola.
\newblock Stochastic {F}rank-{W}olfe methods for nonconvex optimization.
\newblock In {\em 2016 54th Annual Allerton Conference}, pages 1244--1251.
  IEEE, 2016.

\bibitem{RobbMonr51}
Herbert Robbins and Sutton Monro.
\newblock A stochastic approximation method.
\newblock {\em The Annals of Mathematical Statistics}, 22(3):400--407, 1951.

\bibitem{RobbSieg71}
Herbert Robbins and David Siegmund.
\newblock A convergence theorem for nonnegative almost supermartingales and
  some applications.
\newblock In Jagdish~S. Rustagi, editor, {\em Optimizing Methods in
  Statistics}. Academic Press, 1971.

\bibitem{Stei83}
T.~Steihaug.
\newblock The conjugate gradient method and trust regions in large scale
  optimization.
\newblock {\em SIAM Journal on Numerical Analysis}, 20(3):626--637, 1983.

\bibitem{DStirzaker_2003}
David Stirzaker.
\newblock {\em Elementary probability}.
\newblock Cambridge University Press, 2003.

\bibitem{zhang2020one}
Mingrui Zhang, Zebang Shen, Aryan Mokhtari, Hamed Hassani, and Amin Karbasi.
\newblock One sample stochastic {F}rank-{W}olfe.
\newblock In {\em AISTATS}, pages 4012--4023, 2020.

\end{thebibliography}

\end{document}